\documentclass[12pt,leqno]{book}
\usepackage{graphicx}
\usepackage{amsfonts, color, amsmath, amscd, amsfonts}
\usepackage{makeidx, endnotes, wrapfig}
\usepackage{mystyle}

\begin{document}

\title{Planar Brownian Motion and Complex Analysis}

\author{Greg Markowsky}

\maketitle

\tableofcontents

%\frontmatter
%\tableofcontents
%\include{preface}
\setcounter{chapter}{-1}

\mainmatter
\chapter{Preface}

This was originally intended as a graduate level textbook on the connections between planar Brownian motion and complex analysis. Then back problems forced my retirement as a writer. I still have hope that this partial draft may be of use to someone, somewhere, which is why I have placed it online. I am deeply indebted to Maher Boudabra and Andrew Buttigieg for their helpful comments, enthusiasm, and proofreading. The book, such as it is, is dedicated to my ever patient Juyoun, without whom it would never have been written.
\chapter{Preliminaries} \lll{chzero}

\section{Basic notions from probability theory} \label{zeroBM}

We begin with a few basic notions from probability theory, and define Brownian motion. This chapter is largely meant for review and in order to fix notation; a reader unfamiliar with the topics discussed here is advised to consult an introductory probability textbook. A probability space is a triple $(\Om, \FF, P)$. Here $\Om$ is any set whatsoever; $\FF$ is a $\sigma$-field: a collection of subsets of $\Om$, containing $\Om$ itself, and which is closed under complements, finite intersections, and countable unions; $P$ is a probability measure. To say that $X$ is a random variable (abbreviated {\it r.v.}) on a probability space $(\Om, \FF, P)$ means that $X=X(\om)$ is a function on $\Om$, with $X^{-1}(A) \in \FF$ for all Borel sets $A$.  We will also say $X$ is {\it $\FF$-measurable}. If $X$ is a random variable, there will be a unique smallest $\si$-field upon which $X$ is measurable; this is referred to as the {\it $\si$-field generated by $X$}, and is denoted $\si(X)$. Two $\si$-fields $\FF$ and $\GG$ are {\it independent} if $P(A \cap B) = P(A)P(B)$ for all $A \in \FF, B \in \GG$. Two random variables $X,Y$ are independent if the $\si$-fields $\si(X), \si(Y)$ are independent.

\vski

Suppose that $X$ is a r.v. taking values in $\RR^n$. If there is a function $f_X(t)$ on $\RR^n$ such that, for any Borel set $A \subseteq \RR^n$, we have

\begin{equation} \label{}
P(X \in A) = \int_{A} f_X(t) dt,
\end{equation}

then we will call $f(x)$ the {\it probability density function} (abbreviated as {\it p.d.f}) for $X$. If $g$ is a measurable function on $\RR^n$, then we define the {\it expectation} of $g(X)$ to be

\begin{equation} \label{}
E[g(X)] = \int_{\RR^n} g(t) f_X(t)dt.
\end{equation}

If $X$ takes values on $\RR$, the {\it mean} of $X$ is then $E[X]$, and the {\it variance} is $\Var(X) = E[X^2]-E[X]^2$. %It is not hard to see in general that if $X$ and $Y$ are r.v.'s which have p.d.f.'s and which take values on $\RR^n$ and $\RR^m$ respectively, then the p.d.f. for the r.v. $(X,Y)$ taking values in $\RR^{n+m}$ satisfies $f_{(X,Y)}((t_t,t_2)) = f_X(t_1)f_Y(t_2)$.
If $X$ and $Y$ are independent, then $E[XY]=E[X]E[Y]$ and $\Var(X+Y) = \Var(X) + \Var(Y)$. The {\it normal} or {\it Gaussian distribution} with mean $\mu$ and variance $\si^2$ is a distribution on $\RR$ characterized by the p.d.f.

\begin{equation} \label{}
\rho(x) = \frac{1}{\si \sqrt{2\pi}} e^{-\frac{(x-\mu)^2}{2 \si^2}}.
\end{equation}

We will use the notation $X \sim N(\mu, \si^2)$ to signify that $X$ is a normally distributed r.v. with mean $\mu$ and variance $\si^2$.

\vski

A {\it stochastic process} is an indexed family of r.v.'s $\{X_t: t \geq 0\}$ and an increasing collection of $\si$-fields $\{\FF_t: t \geq 0\}$ (which means that $\FF_s \subseteq \FF_t$ whenever $s \leq t$) so that each $X_t$ is $\FF_t$-measurable. The collection $\{\FF_t\}$ is called the {\it filtration} of the process, and in some sense represents the information available to an observer at any time. In general, we take $\FF_t$ to be the $\si$-field generated by $\{X_s: 0 \leq s \leq t\}$. %Let $X$ be a $\FF$-measurable random variable on a space $\Om$, and let $\GG$ be a $\si$-field on $\Om$ with $\GG \subseteq \FF$, so that $X$ is not necessarily $\GG$-measurable. There is a $\GG$-measurable random variable, denoted $E[X|\GG]$ and referred to as the {\it conditional expectation of X with respect to $\GG$}, such that $E[X 1_A] = E[E[X|\GG]1_A]$ for all $A \in \GG$. There are a few rules for this:

%\begin{itemize}
%    \item[$\bullet$] $E[aX+bY|\GG] = aE[X|\GG]+bE[Y|\GG]$.%

%    \item[$\bullet$] If $\GG = \{\emptyset,\Om\}$, then $E[X|\GG]=E[X]$.

%    \item[$\bullet$] If $X$ is $\GG$-measurable, then $E[X|\GG]=X$, and more generally $E[XY|\GG]=XE[Y|\GG]$

%    \item[$\bullet$] If $\GG_1 \subseteq \GG_2$, then $E[E[X|\GG_2]|\GG_1] = E[X|\GG_1]$.

%    \item[$\bullet$] If $\si(X)$ and $\GG$ are independent, then $E[X|\GG]=E[X]$.
%\end{itemize}

%Stochastic processes for which $E[X_t|\FF_s] = X_s$ whenever $t>s$ are called martingales.

\vski

\section{Brownian motion} \label{zeroBM}

A stochastic process $X_t$ represents a system of some sort changing over time. Brownian motion is the stochastic process which can most aptly be described as "purely random motion". That is, it represents a continuous motion with no memory or momentum. In probabilistic terms, this requires $X_t - X_s$ to be independent of $X_s$ for $s<t$: the process should have no memory of what it did to get to $X_s$ in order to decide what to do to get to $X_t$. The process should also move in a uniform way over different time intervals; that is, $X_t - X_s$ should have the same distribution as $X_{t-s}-X_0$, for instance. However, if these properties hold then we see that we can write $X_t - X_0= (X_{t/n} - X_0) + (X_{2t/n} - X_{t/n}) + \ldots (X_{nt/n} - X_{(n-1)t/n})$ for any $n \in \NN$ and thereby express $X_t$ as the sum of a large number of i.i.d. random variables. This brings us to the realm of the Central Limit Theorem, which states that the sum of a large number of i.i.d. random variables with variances will be approximately normally distributed. It is therefore natural (and necessary) that our Brownian motion be normally distributed.

{\it One dimensional Brownian motion} $B_t$ is defined to be a stochastic process with the following properties.

\begin{itemize}

\item[(i)] $B_t: [0,\ff) \lar \reals$ is continuous a.s.

\item[(ii)]  For $t_1 < \ldots < t_n$ the random variables $B_{t_1}, B_{t_2}-B_{t_1}, \ldots , B_{t_n}-B_{t_{n-1}}$ are independent.

\item[(iii)]  For $t'<t''$, $B_{t''}-B_{t'}$ is a normally distributed variable with mean 0 and variance $t''-t'$.

\end{itemize}

In general, we will assume that $B_0 = a$ a.s. for some constant $a$, and we will use the notations $P_a$ and $E_a$ to mean probability and expectation taken with respect to this initial condition. That Brownian motion even exists is nontrivial, but well-established, and its study is a fascinating and rich endeavor practiced by many. In this text we will focus on properties of Brownian motion relevant to complex analysis, but so that the reader may get a better feel for the process we briefly describe the trajectories (the image of $B_t$ for a fixed $\om$). With probability 1, Brownian paths are nondifferentiable for every $t$ and unbounded but recurrent, which means that they return infinitely often to every point. Brownian motion satisfies {\it self-similarity}, which means roughly that the path looks identical regardless of the scale at which it is viewed; this can be seen by noting that $\frac{1}{c}B_{c^2t}$ for $c \in (0,\ff)$ satisfies the requirements above and is therefore a Brownian motion. The process is therefore revealed as possessing non-smooth, crinkled-up paths which travel slowly throughout $\RR$ but always return to their starting point.

A {\it stopping time} $\tau$ is a random variable taking values in $[0,\ff]$ such that

\be \label{}
\{ \tau \leq t \} \in \FF_t.
\ee

The most commonly used stopping times are {\it hitting times}, and we will make extensive use of them in this book. For a closed set $K$ in $\CC$, let $T_K = \inf\{t \geq 0: B_t \in K\}$. $T_K$ is the hitting time of $K$, and is a stopping time. In fact, strictly speaking the set $K$ need not be closed in order for $T_K$ to be a stopping time; however, we choose to require this because for an open set $\Om$ we define $T_\Om = T_{\dd \Om} =T_{\Om^c}$, the hitting time of the complement of $\Om$ (or the {\it exit time} of $\Om$). This is a bit of an abuse of notation; however, the point is that we will be generally running Brownian motion on open sets and stopping it when it exits the set or hits a particular closed set, and we tailor our definitions accordingly.

\vski

The importance of stopping times can be seen in the following theorem, as well as in It\=o's Theorem in the next chapter.

\begin{theorem}[Strong Markov property] \label{StrongMP}
If $\tau$ is a finite stopping time, and $\hat{B}_t = B_{\tau+t}-B_\tau$, then $\hat{B}_t$ is a Brownian motion independent of $\{B_t: 0 \leq t \leq \tau\}$.
\end{theorem}

The proof of this theorem is also nontrivial. The reader interested in seeing the finer details of the result, of the construction of Brownian motion, or of any of the other properties given here is referred to \ref{rogers} as well as any number of other excellent texts on the subject available in any academic library.

\vski

The primary object of interest in this book will be {\it planar Brownian motion}, which is simply the process $(R_t,I_t)$ in $\RR^2$, where $R_t$ and $I_t$ are independent one dimensional Brownian motions. Henceforth, any reference to Brownian motion will be to the planar variety, unless otherwise specified. The following important observation leads to the development in the rest of the book.

\begin{theorem} \label{trublu}
\begin{itemize} \label{}

\item[(i)] Brownian motion is rotation invariant. That is, if $\Theta$ is a rotation of the plane, then $Y_t=\Theta(B_t)$ is again a Brownian motion.

\item[(ii)] If $c \in [0,\ff)$, then $Z_t = \frac{1}{\sqrt{c}} B_{ct} = (\frac{1}{\sqrt{c}} R_{ct}, \frac{1}{\sqrt{c}} I_{ct})$ is again a Brownian motion.

\end{itemize}
\end{theorem}

{\bf Proof:} That $Y_t$ and $Z_t$ are continuous with independent increments are clear. It is also clear that the components of $Z_t$ are independent, and each component can be shown to be normally distributed with the correct mean and variance by a simple calculation. Since the components of Brownian motion are independent, the p.d.f. for $B_t$ starting at 0 is simply the product of the p.d.f.'s of the components, namely $\rho_t(x,y) = \Big(\frac{1}{\sqrt{2 \pi t}} e^{-x^2/2t}\Big)\Big(\frac{1}{\sqrt{2 \pi t}} e^{-y^2/2t}\Big) = \frac{1}{2 \pi t} e^{-(x^2+y^2)/2t}$. Converting to polar coordinates, this is $\rho_t(r,\th) = \frac{1}{2 \pi t} e^{-r^2/2t}$, which is rotation invariant (i.e. does not depend on $\th$). Thus, $\Theta(B_t)$ has the same distribution as $B_t$, namely a bivariate normal with independent components which are normalized to have mean $0$ and variance $t$. A similar argument holds if the Brownian motion starts at a point other than 0. \qed

Let us examine this theorem more closely. A rotation by $\th$ takes the point $(r,\th_0)$ in polar coordinates to $(r,\th_0 + \th)$, or in rectangular coordinates takes $(r\cos \th_0, r \sin \th_0)$ to $(r\cos (\th_0+ \th), r \sin (\th_0 + \th))$. Using the angle sum formulas for $\cos, \sin$, this is $(r\cos \th_0 \cos \th - r \sin \th_0 \sin \th, r \sin \th_0 \cos \th+ r \cos \th_0 \sin \th)$. If we then multiply both components by $c \in \RR^+$, we arrive at the point $(r c \cos \th_0 \cos \th - r c \sin \th_0 \sin \th, r c \sin \th_0 \cos \th+ r c \cos \th_0 \sin \th)$, and we know that this composition of maps takes Brownian motion to Brownian motion (with the clock being changed according to Theorem \ref{trublu}). Associating this map with the point $(c \cos \th, c \sin \th)$ in $\RR^2$, we are led naturally to the operation $(x,y) \dot (x',y') = (xx' -yy', xy' + x'y)$, and we will call this {\it complex multiplication}. Associating the point $(1,0)$ with the real number 1, and using the symbol $i$ for $(0,1)$, we can write any vector in $\RR^2$ as $x+yi$, and complex multiplication becomes $(x+yi)(x'+y'i) = (xx' - yy') + (xy' + x'y)i$. Note that $i^2 = -1$, which is often taken as the starting point of an investigation into complex numbers. We will refer to $\RR^2$ with this multiplication as the {\it complex plane}, and write $\CC$ to denote it.

%\begin{proposition} \label{trublu}
%Suppose that $W_t, W'_t$ are independent Brownian motions, and $a,b \in \RR$. Set $Y_t = \aa W_t + \bb W'_t$. Then $Y_{t/(\aa^2 + \bb^2)}$ is again a Brownian motion. Furthermore, if we set $Z_t = \bb W_t - \aa W'_t$, then $Y_t$ and $Z_t$ are independent.
%\end{proposition}

\section{Complex numbers} \label{zerocompnum}

Let $z = x + yi \in \CC$. We will refer to the real numbers $x$ and $y$ as the {\it real} and {\it imaginary parts} of $z$, and write $Re(z) = x, Im(z) = y$. Define the {\it modulus} of $z$ to be $|z| = \sqrt{x^2 + y^2}$,  the ordinary Euclidean metric. We define $\bar z = x- iy$ and refer to $\bar z$ as the {\it conjugate} of z; one of the primarily uses of $\bar z$ lies in the identity $z \bar z = |z|^2$. This allows us to define quotients on $\CC$ by reducing to division by real numbers, since if $ w \neq 0$ then we may write

\begin{equation} \label{}
\frac{z}{w} = \frac{z\bar w}{|w|^2} = \frac{Re(z \bar w)}{|w|^2} + i \frac{Im(z \bar w)}{|w|^2}.
\end{equation}

Armed with these definitions, it is straightforward to verify the following important facts.

\begin{proposition} \label{modproprop}

\begin{itemize} \label{}

\item[(i)] $|zw| = |z||w|$.

\item[(ii)] $|\frac{z}{w}| = \frac{|z|}{|w|}$.

\item[(iii)] (Triangle inequality) $|z+w| \leq |z| + |w|$.

\item[(iv)] (Reverse triangle inequality) $|z \pm w| \geq ||z| - |w||$.

\item[(v)] $|z \pm w|^2 = |z|^2 \pm 2 Re (z \bar w) + |w|^2$.
\end{itemize}
\end{proposition}

The proof of this proposition is left to the reader. Fundamental to the study of the complex numbers is {\it Euler's formula}:

\begin{equation} \label{}
e^{i \th} = \cos \th + i \sin \th.
\end{equation}

This can be verified by considering the well-known Taylor series for the real-valued functions $e^x, \cos x, \sin x$. This allows us to extend the domain of the exponential function to the entire complex plane by $e^{x+yi} = e^x(\cos y + i \sin y)$. The exponential so extended satisfies virtually every property one would like; the following proposition lists some of the more basic ones.

\begin{proposition} \label{exproprop}
\begin{itemize} \label{}

\item[(i)] If $z=x + yi$, then $|e^z| = e^x$, and $\frac{Im (e^z)}{Re (e^z)} = \tan y$.

\item[(ii)] The function $e^z$ is periodic with period $2\pi i$. That is, $e^z = e^{z+2\pi i}$.

\item[(iii)] $e^{z_1+z_2} = e^{z_1} e^{z_2}$.

\item[(iv)] $e^{ \bar z} = \seg{10}{e^z}$.

\end{itemize}
\end{proposition}

Again, we leave the proof to the reader. Part $(i)$ includes an important new fact, that the angle which the vector $e^z$ makes with the positive real axis is controlled by $y$. We will define the {\it argument} of any $z \in \CC$ to be $\mbox{arg }(z) = \tan^{-1}\frac{Im(z)}{Re(z)}$; this is clearly an ambiguous definition, since the value of $\tan^{-1}$ is only defined up to multiples of $2\pi$, and so for exactitude we define the {\it principal argument} of $z$ to be $\Arg(z) = \Tan^{-1}\frac{Im(z)}{Re(z)}$, where here $\Tan^{-1}$ is the branch of $\tan^{-1}$ which takes values in $(-\pi, \pi]$. We see that $\Arg z$ is well defined for all $z \neq 0$, although the function $\Arg z$ is discontinuous on the negative real axis. We see that for any $z =x+yi \neq 0$ we can write $z = re^{i \th}$, where $r=|z|$ and $\th = \Arg z$ (or, for that matter, any branch of $\arg z$). We will consider this to be the polar form of $z$. Note that the properties in Propositions \ref{modproprop} and \ref{exproprop} show that $(r_1e^{i\th_1})(r_2e^{i\th_2}) = r_1r_2 e^{i(\th_1 + \th_2)}$.

\vski

Let us now consider a number of simple transformations of $\CC$. We have seen that a rotation of the plane can be realized by multiplication by a complex number with modulus 1, which can be expressed as $e^{i\th}$ for some $\th \in \RR$. That is, the function $f(z) = e^{i\th} z$ is simply a rotation of $\CC$ by $\th$ radians counterclockwise. By an {\it expansion} of $\CC$, we mean a function of the form $f(z) = rz$, where $r \in [0,\ff)$; this fixes the arguments of all points in $\CC$, and simply scales the moduli of all points by the multiplicative factor $r$. A {\it dilation} is the composition of a rotation and expansion, and can therefore be expressed as $f(z) = re^{i\th} z$. A {\it translation} will be a map of the form $f(z) = z+c$, which simply shifts the entire plane by $c$. Composing a translation with a dilation gives a general {\it linear function} $f(z) = az + c$. These are essentially the most basic functions from $\CC$ to itself, and it would be well for the reader to entirely familiarize themselves with the geometric properties of these transformations, as they will be important throughout.

\vski

The {\it stereographic projection} is a method of adding a "point at infinity" to $\CC$ in order to make $\CC$ compact, which helps greatly in visualizing various transformations
\endnote{It is also common to perform the stereographic projection with the sphere tangent to the complex plane, for instance by taking the sphere $\{x_1^2 +x_2^2 + (x_3-1)^2 = 1\}$ rather than $\SS$, or by embedding $\CC$ in $\RR^2$ as $\{x_3 = -1\}$. We prefer the setup given in the text, however, as the important sets $\DD=\{|z|<1\}$, $\dd \DD$, and $\{|z|>1\}$ appear naturally on the sphere as the intersections of $\SS$ with $\{x_3 < 0\}, \{x_3 = 0\},$ and $\{x_3 > 0\}$, respectively. }.
Imagine $\CC$ embedded in $\RR^3 = \{(x_1,x_2,x_3)\}$ as the plane $x_3=0$, and consider also the unit sphere $\SS = \{x_1^2+x_2^2+x_3^2 = 1\}$. Associate any point $z=(x, y, 0) \in \CC$ with the point of intersection of $\SS$ and the line passing through $(0,0,1)$ and $z$. We will let $\Theta$ denote this map from $\SS \bsh \{(0,0,1)\}$ to $\CC$, which is in fact a bijection. This picture may help visualize the projection.

\hspace{1.2in} \includegraphics[width=110mm,height=90mm]{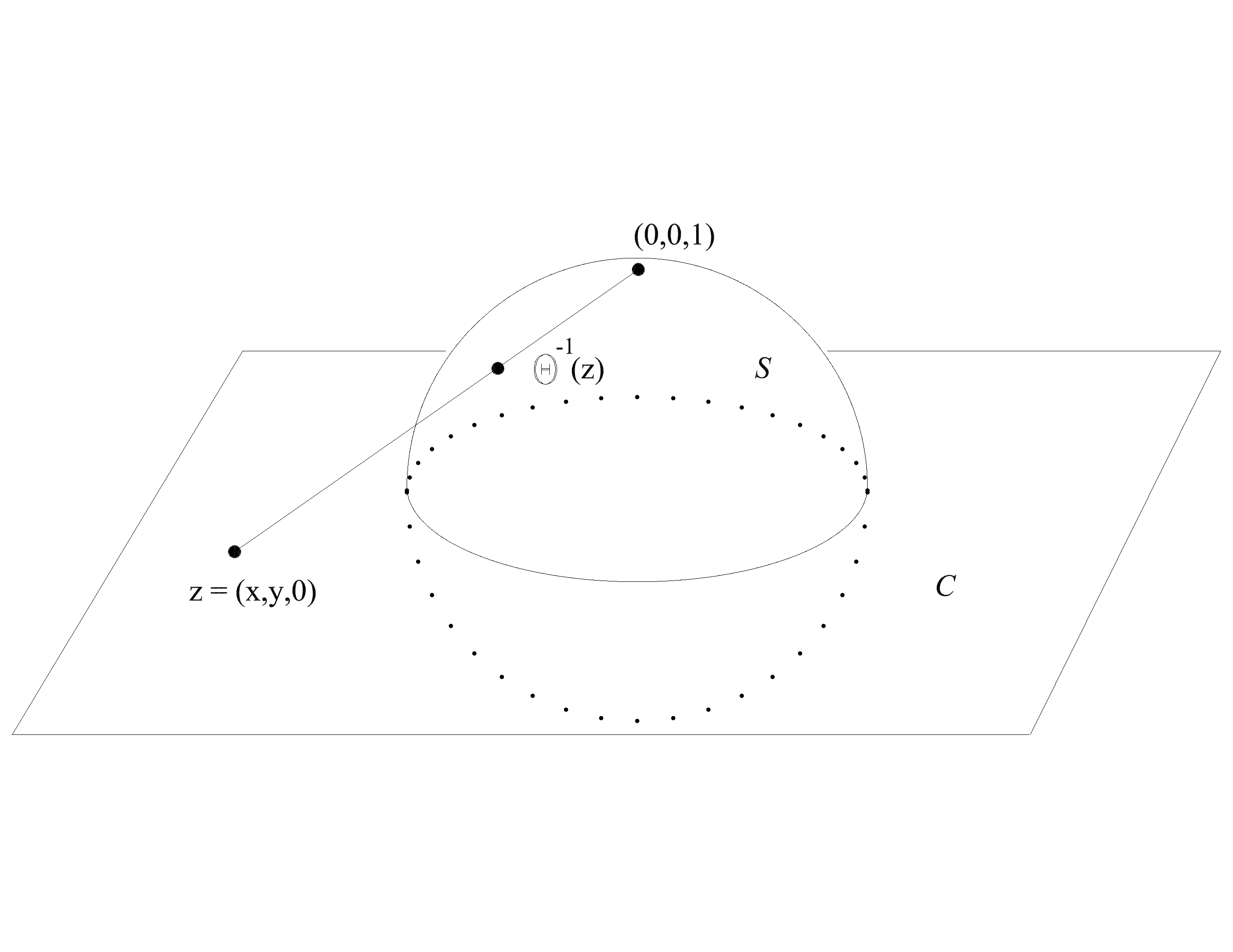}

In this way, $\CC$ is essentially wrapped onto the surface of a sphere, with only the point $(0,0,1)$ not matched to a corresponding point in $\CC$. $(0,0,1)$ therefore plays a special role in this construction, and is considered to be the {\it point at $\ff$}. Note that this is an appropriate label, since as $|z| \lar \ff$ in $\CC$, $\Theta^{-1}(z) \lar (0,0,1)$ on $\SS$. The complex plane $\CC$ together with the point at $\ff$ is commonly referred to as the {\it Riemann sphere}, and denoted $\hat \CC$. The line which passes through $(0,0,1)$ and $z=(x,y,0)$ can be realized by $\ga(t) = (xt,yt,1-t)$, and this intersects $\SS$ when $x^2t^2 +y^2t^2 + (1-t)^2 = 1$, which happens precisely when $(x^2+y^2+1) t^2 - 2t = 0$. $t=0$ is a solution to this quadratic equation, which yields the point of projection $(0,0,1)$. This means that the point $\Theta^{-1}(z)$ is given by the other solution, which occurs when $t=\frac{2}{x^2+y^2+1}$. Thus, $\Theta^{-1}(z) = (\frac{2x}{x^2+y^2+1},\frac{2y}{x^2+y^2+1}, \frac{x^2+y^2-1}{x^2+y^2+1})$. This formula allows us to deduce the following important fact.

\begin{proposition} \label{circles}
If $\CCC$ is a circle on $\SS$, then $\Theta(\CCC)$ is either a circle or line in $\CC$. To be specific, if $\CCC$ passes through $(0,0,1)$ then $\Theta(\CCC)$ is a straight line in $\CC$, while if $\CCC$ does not pass through $(0,0,1)$ then $\Theta(\CCC)$ is a circle in $\CC$.
\end{proposition}

{\bf Proof:} Several proofs exist; an elegant geometric one which might have pleased the ancient Greeks is outlined in Exercise \ref{stereogeo} below. However, we can give a simple analytic one as follows. A circle $\CCC$ on $\SS$ can be realized as the intersection of $\SS$ with a plane with equation $ax_1 + bx_2 + cx_3 = l$, where $a^2 + b^2 + c^2 =1$ and $ l \in (-1,1)$. Substituting the formula for $\Theta^{-1}(z)$ given above into this equation and simplifying leads to

\begin{equation} \label{}
(c-l)(x^2+y^2) + 2ax + 2by = c+l.
\end{equation}
If $c \neq l$, then this is the equation of a circle in $\CC$ (the bijective nature of $\Theta$ shows that the circle does not degenerate). If, however, $c=l$, then it is the equation of a line, and note in this case that the point $(0,0,1)$ lies on $\CC$. \qed

A set $U \subseteq \CC$ is {\it open} if for any $z \in U$ there is a small disk centered at $z$ contained in $U$. A set is {\it closed} if its complement is open. The {\it boundary} of $U$, denoted $\dd U$, is the set of all points $z$ for which there is a sequence $z_n$ in $U$ converging to $z$ and also a sequence $w_n$ in $U^c$ converging to $z$. A set $U$ is {\it connected} if, for any $z,w \in U$ there is a path in $U$ connecting them. Equivalently, an open set is connected if it cannot be expressed as the union of two nonempty, disjoint, open sets. A {\it domain} is an open, connected subset of $\CC$. The generic picture of a domain therefore looks something like

\hspace{1.2in} \includegraphics[width=80mm,height=65mm]{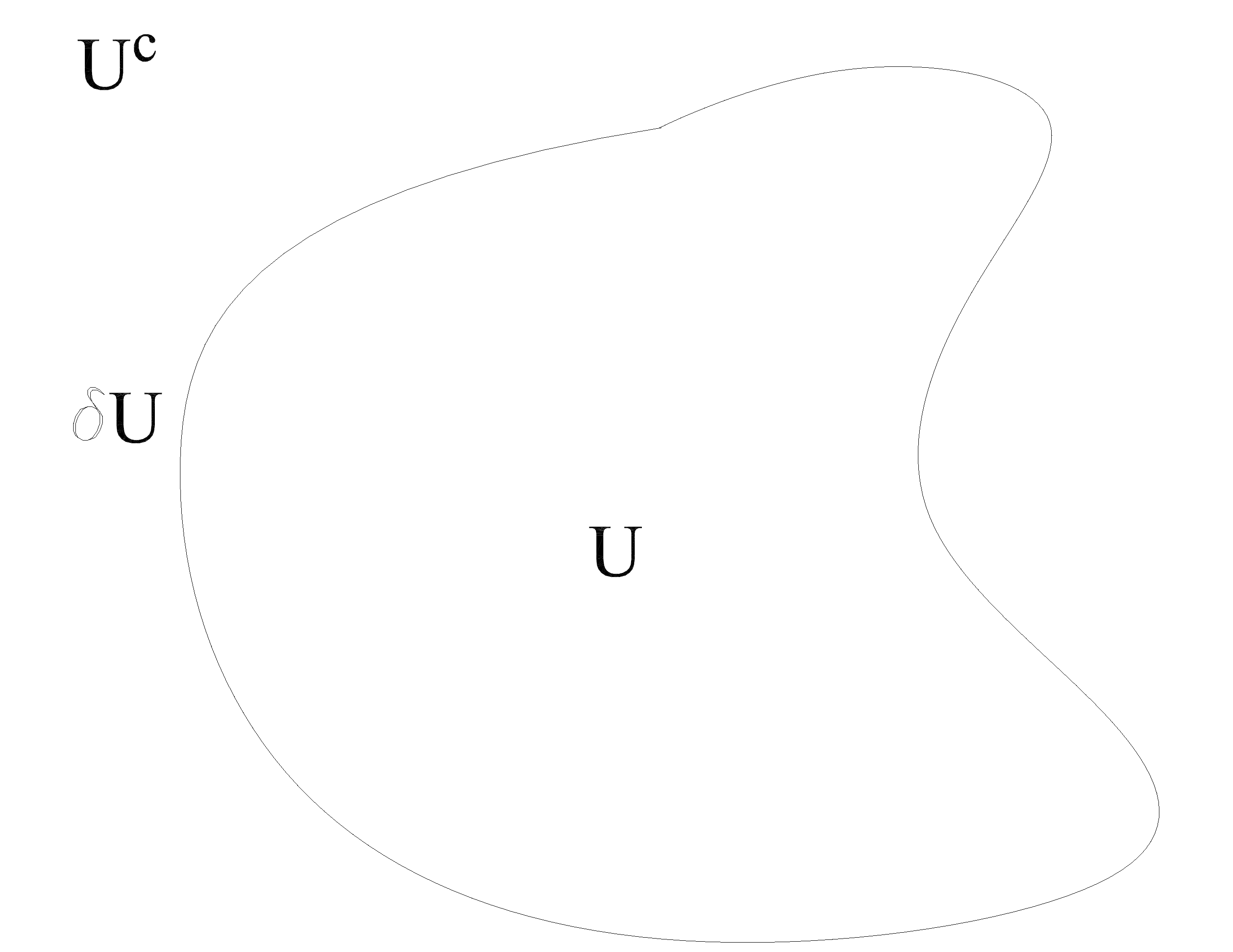}

Our work will take place predominantly on domains: we will generally be concerned with functions defined upon domains, and will run Brownian motions on domains. Several domains are prominent enough to be labeled: the entire plane $\CC$; the unit disk $\DD = \{|z|<1\}$; the upper half-plane $\HH = \{Im(z)>0\}$. We will also use the notation $r\DD$ to denote  $\{|z|<r\}$.

\section{Exercises}

\ccases{gausprop} Verify that the p.d.f. given for normal r.v.'s in the text really is a p.d.f.; that is, show that

\be \label{}
\frac{1}{\si \sqrt{2 \pi}} \int_{-\ff}^{\ff} e^{\frac{-(x-\mu)^2}{2\si^2}}dx = 1.
\ee
Show that if $X$ has this p.d.f., then indeed $E[X] = \mu$ and $\Var(X) = \si^2$. (Hint: For the first part, consider $(\int_{-\ff}^{\ff}e^{-x^2}dx)(\int_{-\ff}^{\ff}e^{-y^2}dy)$ as an integral over the plane $\RR^2$, and convert to polar coordinates. For the third, make use of the first and integration by parts.)

\vski

\ccases{gausmoment} Let $X$ be an $N(0,1)$ random variable. Show that

\be E[X^n] = \left \{ \begin{array}{ll}
0 & \qquad  \mbox{if  $n>0$ is odd}  \\
\frac{(2m)!}{2^m m!} & \qquad \mbox{if $n=2m\geq 0$ is even(Recall $0! = 1$)}
\end{array} \right. \ee

by applying integration by parts and induction to the integral

\be \label{}
\frac{1}{\sqrt{2\pi}} \int_{-\ff}^{\ff} x^n e^{-\frac{x^2}{2}}dx.
\ee

\ccases{transBM} Let $B_t$ be a one dimensional Brownian motion. Show that the following processes are also Brownian motions on $[0,T]$, where $T \in (0,\ff)$.

a) $X_t = -B_t$

b) $X_t = B_T-B_{T-t}$

c) $X_t = \frac{1}{c} B_{c^2 t}$

d) $X_t = tB_{1/t}$ for $t>0$, $X_0=0$ (showing continuity of paths at $t=0$ here is difficult, so it may be assumed).

\vski

\ccases{infts} Let $B_t$ be a one dimensional Brownian motion. If $t,s > 0$, show that $E[B_tB_s]=\min(t,s)$.

\vski

\ccases{mppp} Prove the properties of the modulus listed in Proposition \ref{modproprop}.

\vski

\ccases{eppp} Prove the properties of the exponential listed in Proposition \ref{exproprop}.

%\ccases{hfd} For any $z, w \in \CC$, show that

%\begin{equation} \label{zeromodsquare}
%|z \pm w|^2 = |z|^2 \pm 2 Re z \bar w + |w|^2.
%\end{equation}

\vski

\ccases{stereo77} Show that if $(x_1, x_2, x_3) \in \SS$ with $x_3 \neq 1$ then $\Theta((x_1, x_2, x_3)) = (\frac{x_1}{1-x_3}, \frac{x_2}{1-x_3},0)$.

\vski

\ccases{stereogeo} The purpose of this exercise is to prove Proposition \ref{circles} by purely geometric arguments.

i) Let $l$ be a straight line in $\CC$. Show that under the stereographic projection, $l$ corresponds to a circle on $\SS$ which passes through the point at $\ff$ using the fact that, in $\RR^3$, a line and any point not on the line are contained in a common plane.

ii) Let $c$ be a circle in $\CC$. Show that under the stereographic projection, $c$ corresponds to a circle on $\SS$ which does not pass through the point at $\ff$ by following the ensuing steps.

\vski

a) We will consider $\CCC$ to be the intersection of $\SS$ with a plane $\PP$, and by rotating if necessary we can assume that the maximal $x_3$ values on $\CCC$ occur when $x_2=0$. We then take the $x_2=0$ slice of the configuration, and obtain something like the following picture.

\hspace{1.1in} \includegraphics[width=120mm,height=90mm]{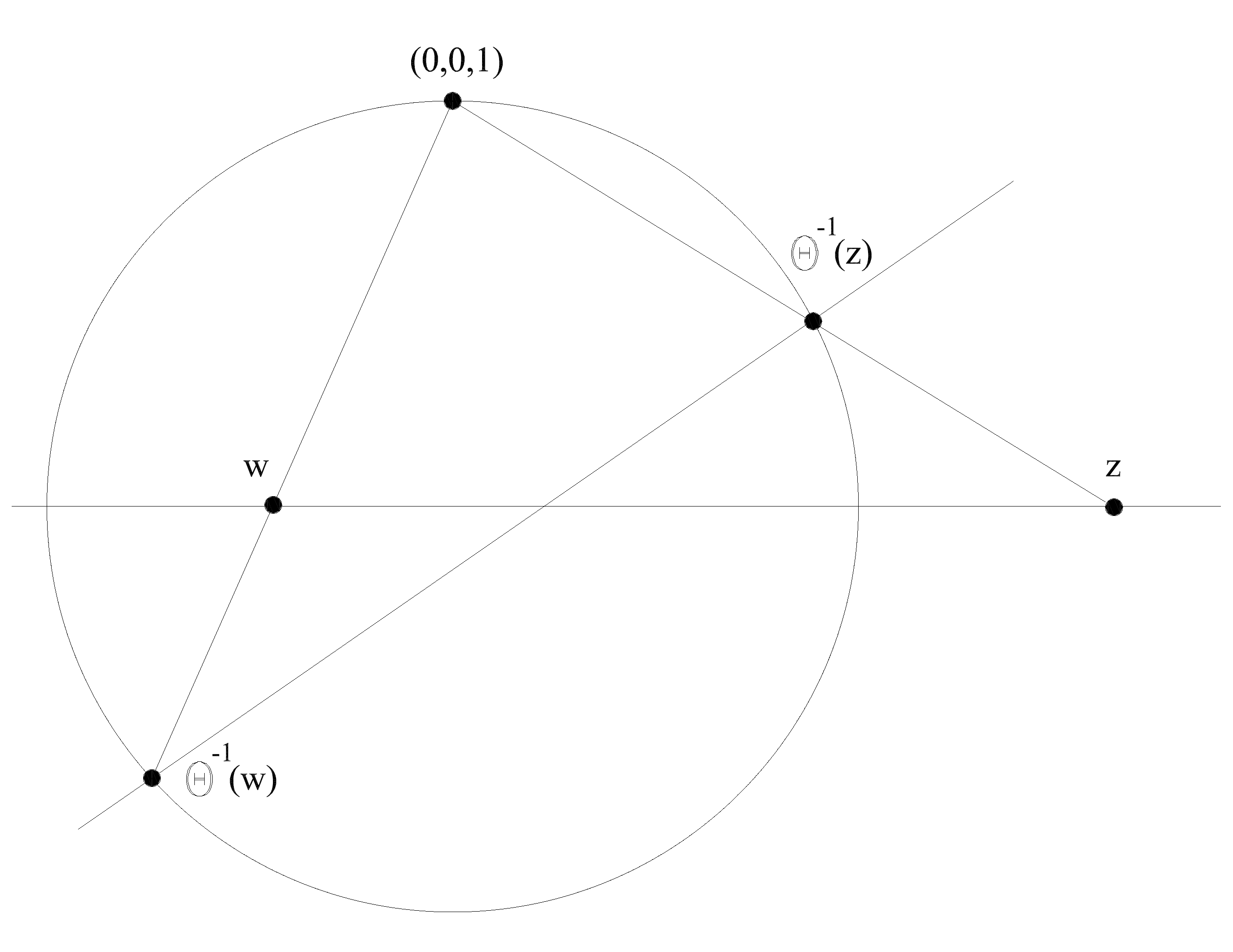}

Show that the triangles with vertices $\{(0,0,1), z, w\}$ and $\{(0,0,1), \Theta^{-1}(z), \Theta^{-1}(w)\}$ are similar.

\vski

b) Conclude that the elliptical cone passing through $\CCC$ is intersected by $\PP$ at the same angle as by $\CC=\{x_3 = 0\}$. Conclude that the image under $\Theta$ of $\CCC$ is a circle in $\CC$.

\vski

\ccases{stereoflip} a) For nonzero $z \in \CC$, show that $\Arg z = \Arg \frac{1}{\bar z}$.

\vski

b) Show that the map $f(z) = \frac{1}{\bar z}$ with $f(0) = \ff$ and $f(\ff) = 0$ on $\CC$ corresponds to a reflection over the plane $\{x_3 = 0\}$ on the Riemann sphere.

%\section{Notes}

\theendnotes

\setcounter{cccases}{0}
\chapter{Dynkin's formula and L\'evy's Theorem} \lll{itochap}

In this chapter we begin with the formula which governs the manner at which the expectation of $u(B_t)$ evolves over time, where $u$ is any $C^2$ function on a domain in $\CC$. We will see that this formula, which we will refer to as {\it Dynkin's formula}, leads naturally to the definitions of harmonic and analytic functions and the power series expansion of each. It may therefore be considered one of the fundamental formulas in complex analysis, in addition to its obvious importance in probability theory.
%
%\endnote{What we refer to as Dynkin's formula is in fact the expectation taken of both sides of

%\begin{equation} \label{itonormal}
%u(B_{T_\ga})
%\end{equation}
%
%\rrr{itonormal} is what is commonly referred to as Dynkin's formula in most texts. Here we make no use of the stochastic integral terms, and defining them would take us to far astray, so we have simplified by taking expectations.

\section{Dynkin's formula} \label{itoform}

The {\it Laplacian differential operator} $\Del$ is defined by $\Del = \frac{d^2}{dx^2} + \frac{d^2}{dy^2}$. This can be applied to a real-valued function on a domain in $\CC$, and if $u$ is complex-valued, so that $u = u_1+iu_2$ with $u_1, u_2$ real-valued, then we simply have $\Del u = (\Del u_1) + i (\Del u_2)$. The following theorem indicates the prominent position of the Laplacian in the study of Brownian motion.

\begin{theorem} \label{itothm}
Suppose that $u(z) = u(x,y)$ is a $C^2$ function on a bounded domain $\Om \subseteq \CC$, and let $a \in \Om$. Suppose $\tau < T_\Om$ with $E[\tau]<\ff$ is a stopping time such that $u$ and its derivatives up to order 2 are a.s. bounded on the set $\{B_t: 0 \leq t \leq \tau, B_0 = a \mbox{ a.s.}\}$. Then

\begin{equation} \label{cleo}
E_a[u(B_{\tau})] - u(a) = \frac{1}{2} E_a \int_{0}^{\tau} \Del u(B_s) ds.
\end{equation}
If $u$ extends continuously to $\dd \Om$, then the same conclusion holds even if we relax the requirement $\tau < T_\Om$ to $\tau \leq T_\Om$.
\end{theorem}

{\bf Proof:} By considering real and imaginary parts separately, we may assume $u$ is real valued. We begin by partitioning $[0,\ff)$ as $0=s_0 < s_1 < s_2 < \ldots$ with $s_n \nearrow \ff$, and write using independence

\begin{equation} \label{donas}
\begin{split}
E_a[u(B_{\tau})] - u(a) & = \sum_{n=0}^{\ff} E_a[u(B_{s_{n+1} \wedge \tau}) - u(B_{s_n \wedge \tau})] \\
& \approx \sum_{n=0}^{\ff} E_a[(u(B_{s_{n+1}}) - u(B_{s_n}))1_{\{\tau \geq s_n\}}]
\end{split}
\end{equation}

Taylor's Theorem allows us to write

%\begin{equation} \label{}
%\begin{split}
%u(B_{s_{n+1}}) - u(B_{s_n}) & = \\
%u_x(B_{s_n})(R_{s_{n+1}} - R_{s_n}) + u_y(B_{s_n})(I_{s_{n+1}} - I_{s_n})& + u_{xx}(B_{s_n})(R_{s_{n+1}} - R_{s_n})^2 \\
%+ 2 u_{xy}(B_{s_n})(R_{s_{n+1}} - R_{s_n})(I_{s_{n+1}} - I_{s_n}) + u_{yy}(B_{s_n})(& I_{s_{n+1}} - I_{s_n})^2 + o(|B_{s_{n+1}} - B_{s_n}|^2)
%\end{split}
%\end{equation}

\bea \label{manti}
&& \hspace{1cm} u(B_{s_{n+1}}) - u(B_{s_n}) =
\\ \nn && u_x(B_{s_n})(R_{s_{n+1}} - R_{s_n}) + u_y(B_{s_n})(I_{s_{n+1}} - I_{s_n}) + \frac{1}{2} u_{xx}(B_{s_n})(R_{s_{n+1}} - R_{s_n})^2
\\ \nn && + u_{xy}(B_{s_n})(R_{s_{n+1}} - R_{s_n})(I_{s_{n+1}} - I_{s_n}) + \frac{1}{2} u_{yy}(B_{s_n})(I_{s_{n+1}} - I_{s_n})^2 + o(|B_{s_{n+1}} - B_{s_n}|^2),
\eea
where as before $R$ and $I$ are the real and imaginary parts of $B$, respectively. Note that $(R_{s_{n+1}} - R_{s_n})$ and $(I_{s_{n+1}} - I_{s_n})$ are independent of $1_{\{\tau \geq s_n\}}$, since $\tau$ is a stopping time, and also independent of the terms of the form $u_x(B_{s_n}), u_{yy}(B_{s_n})$, etc. Thus, when we plug \rrr{manti} into \rrr{donas} we may factor the expectation and use the facts that $E_a[R_{s_{n+1}} - R_{s_n}] = E_a[I_{s_{n+1}} - I_{s_n}] = E_a[(R_{s_{n+1}} - R_{s_n})(I_{s_{n+1}} - I_{s_n})] = 0$ and $E_a[(R_{s_{n+1}} - R_{s_n})^2] = E_a[(I_{s_{n+1}} - I_{s_n})^2] = \frac{1}{2} E_a[|B_{s_{n+1}} - B_{s_n}|^2] = (s_{n+1}-s_n)$. We arrive at

%\begin{equation} \label{}
%E_a[u(B_{s_{n+1}}) - u(B_{s_n})] = \frac{1}{2} E_a[u_{xx}(B_{s_n}) + u_{yy}(B_{s_n})](s_{n+1} - s_n) + o(s_{n+1} - s_n).
%\end{equation}

%Returning to \rrr{donas}, we have

\begin{equation} \label{}
E_a[u(B_{\tau})] - u(a) = \frac{1}{2}\sum_{n=0}^{\ff} E_a[(u_{xx}(B_{s_n}) + u_{yy}(B_{s_n}))1_{\{\tau \geq s_n\}}](s_{n+1} - s_n) + \sum_{n=0}^{\ff} o(s_{n+1} - s_n) E_a[1_{\tau \geq s_n}].
\end{equation}

We now let the partition become finer, so that the mesh size goes to $0$. Since $E_a[\tau] < \ff$, the last term approaches 0, and the first term approaches $\frac{1}{2} E_a \int_{0}^{\tau} \Del u(B_s) ds$, giving \rrr{cleo}. We made use of a number of approximations along the way, but since $u$ and its derivatives are bounded in the interior of $\ga$ these approximations behave correctly in the limit, and the result follows. \qed

{\bf Remark:} The boundedness assumptions on $u$ and $\Om$ are important; see Exercise \ref{ex_unbddom}. The requirement that $\tau$ must be a stopping time is also important; see Exercise \ref{ex_nostoppingtime}.

\vski

The conditions on the stopping time in the statement of the theorem may appear to be difficult to check in practice. However, we will use primarily one type of stopping time in order to apply this theorem: the hitting time of a bounded, closed curve $\ga$ which lies in $\Om$ and separates $a$ from $\dd \Om$.

\hspace{1.2in} \includegraphics[width=80mm,height=65mm]{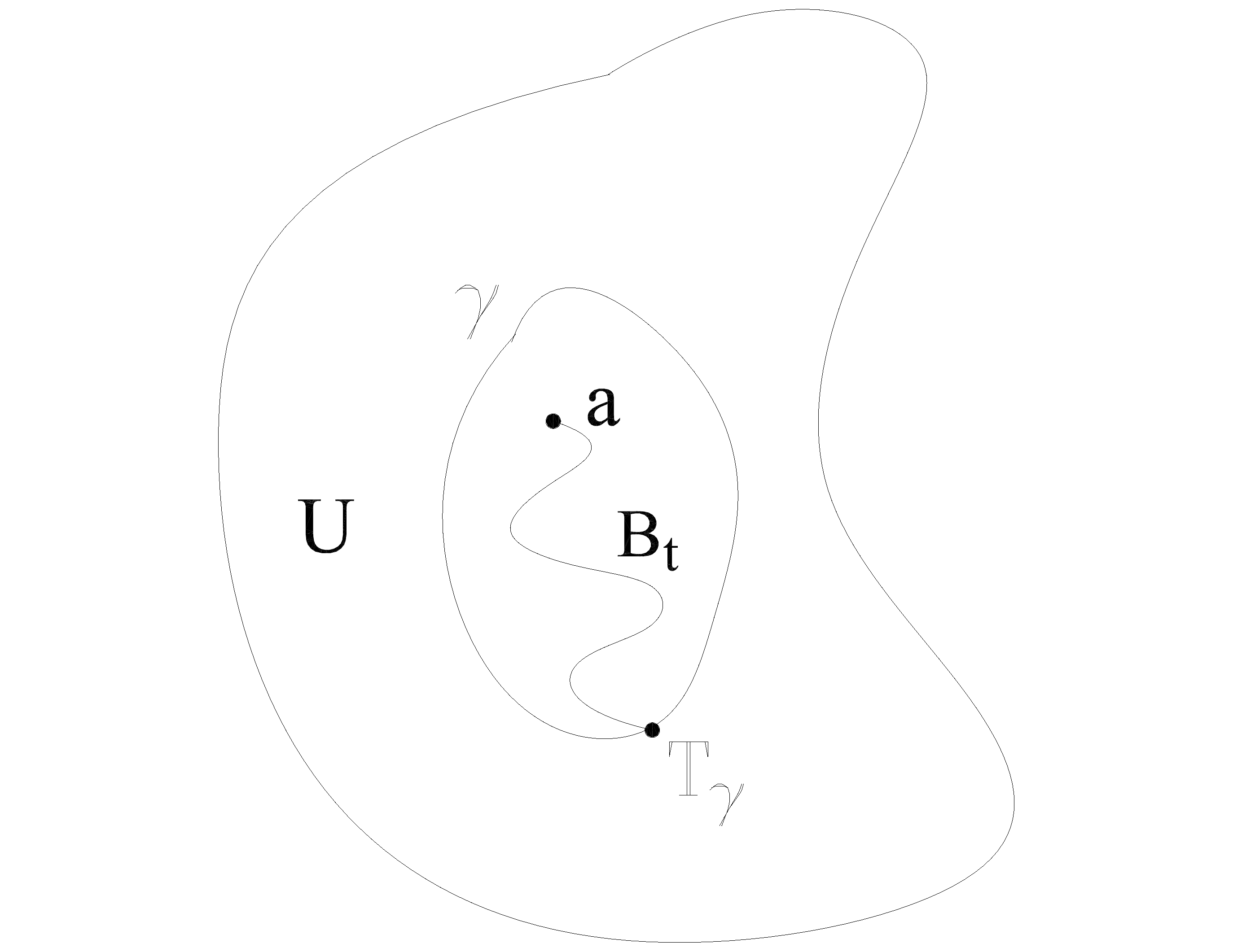}

The fact that $u$ is $C^2$ implies easily that $T_\ga$ automatically satisfies the required boundedness condition on $u$ and its derivatives, so that it must only be shown that $E[T_\ga] < \ff$. Since any bounded curve is contained in a disk centered at $a$, this is an immediate consequence of the following proposition.

\begin{proposition} \label{exittimedisc}
$E_0[T_{r\DD}] = \frac{r^2}{2}$.
\end{proposition}

{\bf Proof:} Fix large $N>0$ and let $\tau = \tau(N) = T_{r\DD} \wedge N$. Note that $E[\tau] < N < \ff$, so we may apply Dynkin's formula with $u(x + yi) = x^2+y^2$. We have $u(0) = 0$ and $\Del u \equiv 4$; we therefore obtain

\begin{equation} \label{}
E_0[u(B_\tau)] = \frac{1}{2} E_0 \int_{0}^{\tau} 4 ds = 2 E_0[\tau].
\end{equation}
Since $u \leq r^2$ on $\{B_t: 0 \leq t \leq \tau\}$, we see that $E[\tau] \leq r^2/2$ independently of $N$, and we may let $N \nearrow \ff$ and apply monotone convergence to conclude that $E[T_{r\DD}] < \ff$. This implies that $T_{r\DD} < \ff$ a.s. so that $E_0[u(B_{T_{r\DD}})] = r^2$, and also allows us to apply Dynkin's formula directly to $T_{r\DD}$ in order to conclude that $r^2 = 2 E_0[T_{r\DD}]$. \qed

We now have a large class of stopping times for which we can apply Dynkin's formula. However, for arbitrary $u$ it is of course difficult to say much about the integral term $\int_{0}^{\tau} \Del u(B_s) ds$; nevertheless, there is one instance in which this expression is easy to determine, and this leads to an immensely valuable definition. We will say that the function $u$ is {\it harmonic} if $\Del u(z) = 0$ for all $z$ in the domain of $u$. It is then immediate from Theorem \ref{itothm} that $E_a[u(B_{\tau})] = u(a)$. Taking $u(x+yi) = x + yi$, which is clearly harmonic, gives

\begin{corollary} \label{}
Suppose $u$ is harmonic and $\tau$ is a stopping time satisfying the conditions of Theorem \ref{itothm}. Then $E_a[u(B_\tau)] = u(a)$.
\end{corollary}
In particular, if we apply this corollary to the identity function $u(z)=z$, which is trivially harmonic, we see that if $E[\tau]<\ff$ is a stopping time such that $\sup\{|B_t|: 0 \leq t \leq \tau, B_0 = a \mbox{ a.s.}\}< \ff$ then $E_a[B_{\tau}] = a$. This fact admits a fruitful interpretation. If $f$ is a function such that $f(B_t)$ is a Brownian motion, then we will have $E_a[f(B_\tau)] = f(a)$. This will hold for a Brownian motion run at a different speed, since all that matters is the values of $B_\tau$, and therefore by Proposition \ref{trublu} will hold for any linear function $f(z) = cz+b$; but since $\tau$ is a variable time it will also hold for a Brownian motion run at a variable speed. In order to make this precise, we define a {\it time change} to be a family of stopping times $\{C_t: t \geq 0\}$ with respect to a given filtration which is a.s. non-decreasing and right continuous in $s$. We will then say that a process $X_t$ is a time changed Brownian motion if there is a time change $C_t$ of $X_t$ such that $X_{C_t}$ is a Brownian motion; or, equivalently, there is a Brownian motion $B_t$ and a time change $C_t$ of $B_t$ such that $B_{C_t}=X_t$ a.s. The prior argument shows that if we have a function $f$ such that $f(B_t)$ is a time changed Brownian motion, then $E_a[f(B_{\tau})] = f(a)$ for all suitable stopping times $\tau$. It should be no great stretch of the imagination to suppose that this equation is extremely useful, though just how useful it will turn out to be may surprise the reader. Suffice it to say that a large percentage of what follows in this book can be reduced at its core to an application of this result. In any case, the next logical task is to identify a large class of functions which map Brownian motion to a time change of Brownian motion. This is the subject of the next section.

\section{Analytic invariance of Brownian motion} \label{confinv}

%It is natural to look for a class of functions which, in some sense, preserve Brownian motion. The following observation will lead us to a more general result.

%\begin{proposition} \label{rotinvBM}
%Let $B_t$ be a Brownian motion starting at $a$. Then, for $c, d \in \CC$, if we set $Y_t = cB_t + d$ then $Y_{t/|c|^2}$ is a Brownian motion starting at $ca + d$.
%\end{proposition}

%{\bf Proof:} Let $c = \aa + \bb i, d=\rho+ \eta i$. Writing real and imaginary parts $Y_t = \tilde R_t + i\tilde I_t$ we have

%\begin{equation} \label{}
%\tilde R_t = \rho + \aa R_t - \bb I_t , \qquad \tilde I_t = \eta + \aa I_t + \bb R_t.
%\end{equation}

%Proposition \ref{trublu} from Chapter \ref{chzero} shows that $\tilde R_{t/|c|^2}$ and $\tilde I_{t/|c|^2}$ are independent Brownian motions, and the result follows.
%\qed

We have seen in Proposition \ref{trublu} that Brownian motion is invariant under dilations and translations, and therefore under linear maps $f(z) = cz+b$, although the time must be adjusted by the square of the expanding factor $|c|$. It turns out that this phenomenon holds for the much more general class of functions which are "locally dilating", in the following sense. We will call a function $f$ on a domain $\Om$ {\it analytic} if it can be expressed as

\begin{equation} \label{}
f(z) = f(z_0) + c(z-z_0) + o(|z-z_0|)
\end{equation}

around any point $z_0 \in \Om$, where $c \in \CC$ is a constant (depending upon $z_0$). That is, $f$ is a continuous function whose behavior close to $z_0$ is dominated by a dilation term. Equivalently, we may say that $f(z)$ is analytic if

\begin{equation} \label{deriv}
\lim_{\Del z \lar 0} \frac{f(z_0 + \Del z) - f(z_0)}{\Del z}
\end{equation}

exists for all $z_0 \in \Om$. In other words, analytic functions are simply complex differentiable functions, and for that reason we will denote the limit in \rrr{deriv} as $f'(z_0)$; however, the reader is forewarned that the behavior of analytic functions is far more restricted than that of differentiable functions of a real variable, and not to attempt to analogize too much between these two types of differentiable functions. This will be revealed as we progress through the subject.

\vski

The primary connection between analytic functions and Brownian motion lies in the following beautiful result, which is known as L\'evy's Theorem.

%\begin{theorem} \label{holinv}
%Let $f$ be analytic and nonconstant on a domain $\Om$, and let $a \in \Om$. Let $B_t$ be a Brownian motion starting at $a$, and set

%\begin{equation} \label{}
%\si_t = \int_{0}^{t \wedge T_\Om} |f'(B_s)|^2 ds.
%\end{equation}

%Set $C_{t} = \inf\{s: \si_s \geq t\}$. Then $f(B_{C_t})$ is a Brownian motion stopped at the stopping time $C_{\si_{T_\Om}}$.
%\end{theorem}

\begin{theorem} \label{holinv}
Let $f$ be analytic and nonconstant on a domain $U$, and let $a \in U$. Let $B_t$ be a planar Brownian motion which starts at $a$ and is stopped at a stopping time $\tau$ such that  $\tau \leq T(U)$ a.s. Set

\begin{equation*} \label{}
\si(t) = \int_{0}^{t \wedge \tau} |f'(B_s)|^2 ds.
\end{equation*}

for $t \in [0,\tau]$. $\si(t)$ is a.s. strictly increasing and continuous, and we let $C_{t} = \si^{-1}(t)$ on $[0,\si(\tau)]$. Then $\hat B_t = f(B_{C_t})$ is a Brownian motion stopped at the stopping time $\si(\tau)$.
\end{theorem}

We will often refer to the new Brownian motion $\hat B_t$ and stopping time $\si(\tau)$ as the {\it projections} of $B_t$ and $\tau$ under $f$. An intuitive proof of the theorem is simple. We know that dilations preserve Brownian motion, but that the time scale must be adjusted by the square of the expanding factor. Therefore, locally dilating maps preserve Brownian motion as well, except that in the more general setting the expansion factor is different at different points in the domain. We must therefore scale the time of the transformed Brownian motion according to the accumulated square of the expansion factors at all of the points that the Brownian motion has visited, which is done via $\si_t$ and its inverse $C_t$. Q.E.D.

\vski

A rigorous proof of Theorem \ref{holinv} would require quite a bit of work, including the introduction of stochastic integration. Although this is a fascinating and important subject in its own right, we choose not to include it in our development. A reader unsatisfied with the intuitive explanation given here is advised to consult \cite{durBM}, \cite{davis}, or \cite{revyor} for complete proofs of the result.

\vski

Note also that, for a analytic function $f$, we have $f(a) = E_{f(a)}[B_{C_{\si_{T}}}] = E_a[f(B_T)]$; this is the same formula which harmonic functions satisfy, though at this stage we do not know that analytic functions are harmonic, or even that they are $C^2$. As we will see in the last section of this chapter, however, analytic functions are indeed harmonic.

\vski

Now that we have the concept of analytic functions, a natural question arises: what functions are analytic? We begin by noting that the standard proofs from calculus that sums, differences, products, quotients (with nonzero denominators), and compositions of differentiable functions are again differentiable carry over directly into the complex setting, so that we may build analytic functions by combining others, with the formulas for the derivatives given by the ordinary sum, product, quotient, and chain rules. We may therefore commence with the observations that the constant function $f(z) = c$ and identity function $f(z) = z$ are trivially analytic, which immediately gives us all polynomials, with the derivative of $\sum_{n=0}^{d} a_n z^n$ given by $\sum_{n=1}^{d} na_n z^{n-1}$. It is natural to ask whether we may carry this over to "infinite degree polynomials"; that is, if $f(z) = \sum_{n=0}^{\ff} a_n (z-z_0)^n$ is a power series, is $f(z)$ analytic? The following proposition answers this question.

\begin{proposition} \label{redrum}
Consider the expression $\sum_{n=0}^{\ff} a_n (z-z_0)^n$. Let $R \in [0,\ff]$ be given by $\frac{1}{R} = \limsup_{n \lar \ff} |a_n|^{1/n}$.

\begin{itemize} \label{}

\item[(i)] The partial sums $\sum_{n=0}^{M} a_n (z-z_0)^n$ converge on $D(z_0,R)$ and do not converge on $\seg{36}{D(z_0,R)}^c = \{|z-z_0|>R\}$ (if $R=0$, this signifies convergence only at $z=z_0$, and if $R = \ff$ this signifies convergence for all $z$).

\item[(ii)] This convergence is uniform on $D(z_0,r)$ for any $r<R$.

\item[(iii)] If we let $f(z) = \sum_{n=0}^{\ff} a_n (z-z_0)^n = \lim_{M \lar \ff} \sum_{n=0}^{M} a_n (z-z_0)^n$, then $f$ is analytic on $D(z_0,R)$, and $f'(z) = \sum_{n=1}^{\ff} n a_n (z-z_0)^{n-1}$.
\end{itemize}
\end{proposition}

{\bf Proof:} If $|z-z_0|< r < R$, then

\begin{equation} \label{}
\Big| \sum_{n=N}^{M} a_n (z-z_0)^n \Big| \leq \sum_{n=N}^{M} (|a_n|^{1/n} |z-z_0|)^n \leq \sum_{n=N}^{\ff} ((\frac{1}{R}+\eps) r)^n,
\end{equation}
where $\eps$ can be made arbitrarily small by requiring $N$ and $M$ to be sufficiently large. We may therefore take $(\frac{1}{R}+\eps) r < 1$, and it follows that $\sum_{n=N}^{M} a_n (z-z_0)^n \lar 0$ uniformly as $N,M \lar \ff$. This implies that $\sum_{n=0}^{M} a_n (z-z_0)^n$ is a Cauchy sequence, uniformly on $D(z_0,r)$, proving $(ii)$ and the first part of $(i)$. The second part of $(i)$ follows easily in a similar manner. To prove $(iii)$, we use the identity $a^n-b^n = (a-b)(a^{n-1} + a^{n-2}b + \ldots + ab^{n-2} + b^{n-1}$ and the uniform convergence to calculate

\begin{equation} \label{}
\begin{split}
f'(z) & =\lim_{h \lar 0} \frac{\sum_{n=0}^{\ff} a_n (z+h-z_0)^n - \sum_{n=0}^{\ff} a_n (z-z_0)^n}{h} \\
& = \lim_{h \lar 0}((z+h-z_0)-(z-z_0))\frac{\sum_{n=0}^{\ff} a_n ((z+h-z_0)^{n-1} + (z+h-z_0)^{n-2}(z-z_0) + \ldots + (z-z_0)^{n-1}) }{h} \\
& = \lim_{h \lar 0} \sum_{n=0}^{\ff} a_n ((z+h-z_0)^{n-1} + (z+h-z_0)^{n-2}(z-z_0) + \ldots + (z-z_0)^{n-1} ) \\
& = \sum_{n=1}^{\ff} n a_n (z-z_0)^{n-1}.
\end{split}
\end{equation}
\qed

We see that convergent power series give us analytic functions, so by taking the Taylor series of well-known transcendental functions, such as the exponential and trigonometric ones, we obtain a wealth of analytic functions. But what about analytic functions not given by power series? It turns out that there are none, at least in a local sense: every analytic function can be expressed as a convergent power series over any disc in its domain. Our next order of business will be to prove this important fact, which we will do in Section \ref{pifseries}. First, however, we must study a prominent class of analytic functions known as M\"obius transformations.

\section{M\"obius transformations} \lll{mobtrans}

This section is devoted to the class of analytic functions known as {\it M\"obius transformations}. We will see that these functions allow us to transform any three points on the plane to any other three, leading to useful simplifications in many cases. This will lead naturally in the next section to an expression for the exit distribution of Brownian motion from a disc, which in turn leads to the Poisson integral formula and the expansion into power series for harmonic functions, two of of the cornerstones of complex analysis. Furthermore, the class of analytic self-maps of the Reimann sphere is coincident with the class of M\"obius transformations. Fluency in these maps is therefore a must for any complex analyst.

\vski

A M\"obius transformation is a function of the form $\phi(z) = \frac{az+b}{cz+d}$ for $a,b,c,d \in \CC, ad-bc \neq 0$. The determinant condition $ad-bc \neq 0$ is simply to ensure that $\phi(z)$ is not a constant. We have seen two types of M\"obius transformations already: the translation $\phi(z) = z+ b$ and the dilation $\phi(z)=az$. Another simple example is the {\it inversion} $\phi(z) = 1/z$. It is a standard and useful fact that these simple forms can be used to build all other M\"obius transformations.

\begin{lemma} \label{mobbb}
Any M\"obius transformation can be expressed as the composition of translations, dilations, and the inversion.
\end{lemma}

{\bf Proof:} Let $\phi(z) = \frac{az+b}{cz+d}$ be an arbitrary M\"obius transformation. If $c=0$ then $\phi$ is merely a linear map, which is clearly the composition of a translation and dilation. If $c \neq 0$, let $\phi_1(z) = z+d/c$, $\phi_2(z) = 1/z$, $\phi_3(z) = \frac{(bc-ad)}{c^2}z$, and $\phi_4(z) = z+a/c$. Then the reader may check that $\phi = \phi_4 \circ \phi_3 \circ \phi_2 \circ \phi_1$. \qed

If $c \neq 0$ then it may seem that $\phi$ is undefined at $z = -d/c$; however, it is far more profitable to consider $\phi$ as a map from the Riemann sphere to itself with the definitions $\phi(-d/c) = \ff$ and $\phi(\ff) = a/c$. Defined as such, M\"obius transformations are bijections from the Riemann sphere to itself (this is left to the reader to check, in Exercise \ref{mobbi} below), and as quotients of polynomials are analytic at all points as well (analyticity at $-d/c$ means that $1/\phi$ is analytic there, and analyticity at $\ff$ means that $\phi(1/z)$ is analytic at 0). If $c=0$, so that $\phi$ is a linear map, we define $\phi(\ff)=\ff$, and likewise obtain a analytic bijection of the sphere to itself. M\"obius transformations are most naturally understood on the Riemann sphere rather than in the plane, as the following important proposition indicates.

\begin{proposition} \label{mobcirc}
M\"obius transformations map circles to circles.
\end{proposition}

{\bf Remark:} This and ensuing statements about circles should be taken in the context of the Riemann sphere, and recall from Section \ref{zerocompnum} that a circle on the sphere is a circle or straight line in the plane. Thus, the image of a circle or straight line in the plane is either a circle or straight line in the plane, but a circle may be mapped to a line and vice versa.

\vski

{\bf Proof:} Lemma \ref{mobbb} shows us that we need only prove the proposition for translations, dilations, and the inversion. We will take it as evident for translations and dilations, and for the inversion we note that a circle on the sphere can be realized as all points $z$ which satisfy

\begin{equation} \label{}
a|z|^2 + bz + c\bar z + d = 0,
\end{equation}
for some $a,b,c,d \in \CC$ (a straight line in the plane occurs precisely when $a=0$). Replacing $z$ by $\frac{1}{z}$ and multiplying through by $|z|^2$ gives

\begin{equation} \label{}
a + b \bar z + cz + d|z|^2 = 0,
\end{equation}
again the equation of a circle on the sphere (see also Exercise \ref{mobcircex} below). \qed

\begin{proposition} \label{mobthreept}
\begin{itemize} \label{}

\item[(i)] The composition of two M\"obius transformations is again a M\"obius transformation. Furthermore, if $\phi(z) = \frac{az + b}{cz + d}$ is a M\"obius transformation, then $\phi^{-1}$ is as well, and is given by $\phi^{-1}(w) = \frac{dw-b}{-cw + a}$.

\item[(ii)] If $(a,b,c), (a',b',c')$ are a pair of sets of three distinct points on the sphere, then there is a M\"obius transformation taking $a \lar a', b \lar b'$, and $c \lar c'$.

\item[(iii)] If $\phi_1$ and $\phi_2$ are M\"obius transformations which agree at three distinct points then they coincide for all $z$.

\end{itemize}
\end{proposition}

{\bf Remark:} The formula for the inverse in $(i)$ may look familiar, and this is no coincidence; see Exercise \ref{mobmat} below.

\vski

{\bf Proof:} $(i)$ is a simple calculation. For $(ii)$, let $\phi_{a,b,c} = \frac{(b-c)(z-a)}{(b-a)(z-c)}$; here $\phi_{a,b,c}(a) = 0, \phi_{a,b,c}(b) = 1$, and $\phi_{a,b,c}(c) = \ff$. Define $\phi_{a',b',c'}$ accordingly. Then $\phi_{a',b',c'}^{-1} \circ \phi_{a,b,c}$ takes $a \lar a', b \lar b'$, and $c \lar c'$. $(ii)$ allows us to pre- and postcompose by appropriate M\"obius transformations in order to reduce $(iii)$ to showing that a M\"obius transformation which fixes $0,1,$ and $\ff$ must be the identity; however, a M\"obius transformation fixing $\ff$ is simply a linear map, and the only degree 1 polynomial with two fixed points is the identity, so the result follows. \qed

Proposition \ref{mobcirc} shows that circles are important when studying M\"obius transformations, and Proposition \ref{mobthreept} shows that triples of points are as well. This suggests that triples of points are closely connected with circles, and this is indeed so.

\begin{proposition} \label{}
If $a,b,c$ are distinct points on $\hat \CC$, then there is precisely one circle on the sphere which passes through all three.
\end{proposition}

{\bf Proof:} If one of the points is $\ff$, then it is easy to see that the circle must be the line in the plane through the other two. Otherwise, construct the perpendicular bisectors of the line segments $\vseg{ab}$ and $\vseg{bc}$. If the two bisectors are parallel, then $a,b,c$ are colinear, and the circle on the sphere in question is the straight line through the three of them. If, on the other hand, the two bisectors intersect, then this point of intersection is the unique point in the plane equidistant from $a,b,c$, and is therefore the center of the unique circle passing through the three of them. \qed

This and the previous propositions lead to an easy test as to whether a circle can be found which passes through four given points: simply find the M\"obius transformation mapping three of the points to $0,1,\ff$, and see where it sends the fourth point. If the fourth point is sent to the real axis, then the four points lie on a circle, but if it is not then they do not. With this end in mind, define the {\it cross ratio} of the points $z_1,z_2,z_3,z_4$ as

\begin{equation} \label{}
(z_1,z_2;z_3,z_4) = \frac{(z_1-z_3)(z_2-z_4)}{(z_2-z_3)(z_1-z_4)}.
\end{equation}

The reader may compare with the formula given in Proposition \ref{mobthreept} to see that this is the image of $z_4$ under the M\"obius transformation taking $z_1 \lar \ff$, $z_2 \lar 0$, and $z_3 \lar 1$. We see that four points lie on a circle if their cross ratio is real. Proposition \ref{mobcirc} then implies that, for any M\"obius transformation $\phi$, $(z_1,z_2;z_3,z_4)$ is real if and only if $(\phi(z_1),\phi(z_2);\phi(z_3),\phi(z_4))$ is real. However, much more is true.

\begin{proposition} \label{}
For any M\"obius transformation $\phi$ and any four points $z_1,z_2,z_3,z_4$,

\begin{equation} \label{}
(\phi(z_1),\phi(z_2);\phi(z_3),\phi(z_4)) = (z_1,z_2;z_3,z_4).
\end{equation}
\end{proposition}

{\bf Proof:} This is a simple calculation for translations, dilations, and the inversion. Proposition \ref{mobbb} therefore gives the result. \qed

%???Put here BM doesn't see points???

For $a \in \DD$ consider the M\"obius transformation

\begin{equation} \label{diskmap}
\psi_a(z) = \frac{z-a}{1-\bar a z}.
\end{equation}

It can be shown that for $|z| \leq 1$ we have $|\psi_a(z)| \leq 1$, with equality precisely when $|z| = 1$ (see Exercise \ref{diskauto} below). Thus, $\psi_a$ is a analytic self-map of $\DD$ sending $a$ to 0. This map will allow us to give an easy solution to an important problem. Suppose we distinguish an interval $\II$ on the unit circle; for $a \in \DD$, what is $P_a(B_{T_\DD} \in \II)$? That is, what is the probability that a Brownian motion starting at $a$ will exit $\DD$ through $\II$? We first note that the rotation invariance of Brownian motion shows that, if $a=0$, then the quantity $P_0(B_{T_\DD} \in \II)$ must be a rotation invariant measure on $\{|z|=1\}$; and this implies that it is equal to $\frac{m(\II)}{2\pi}$, where $m(\II)$ is the Lebesgue measure, or length, of $\II$. For $a \neq 0$, we may simply use the map $\psi_a$ to map $a$ to $0$ and use the fact that $\psi_a(B_t)$ is a time-changed Brownian motion; note that our question is concerned with the exit probability of the Brownian motion but not the time, and therefore we need not even bother with keeping track of the time change. Clearly, the Brownian motion $B_t$ leaves $\DD$ through $\II$ if, and only if, the Brownian motion $\psi_a(B_t)$ leaves $\DD$ through $\psi_a(\II)$:

\hspace{.2in} \includegraphics[width=160mm,height=110mm]{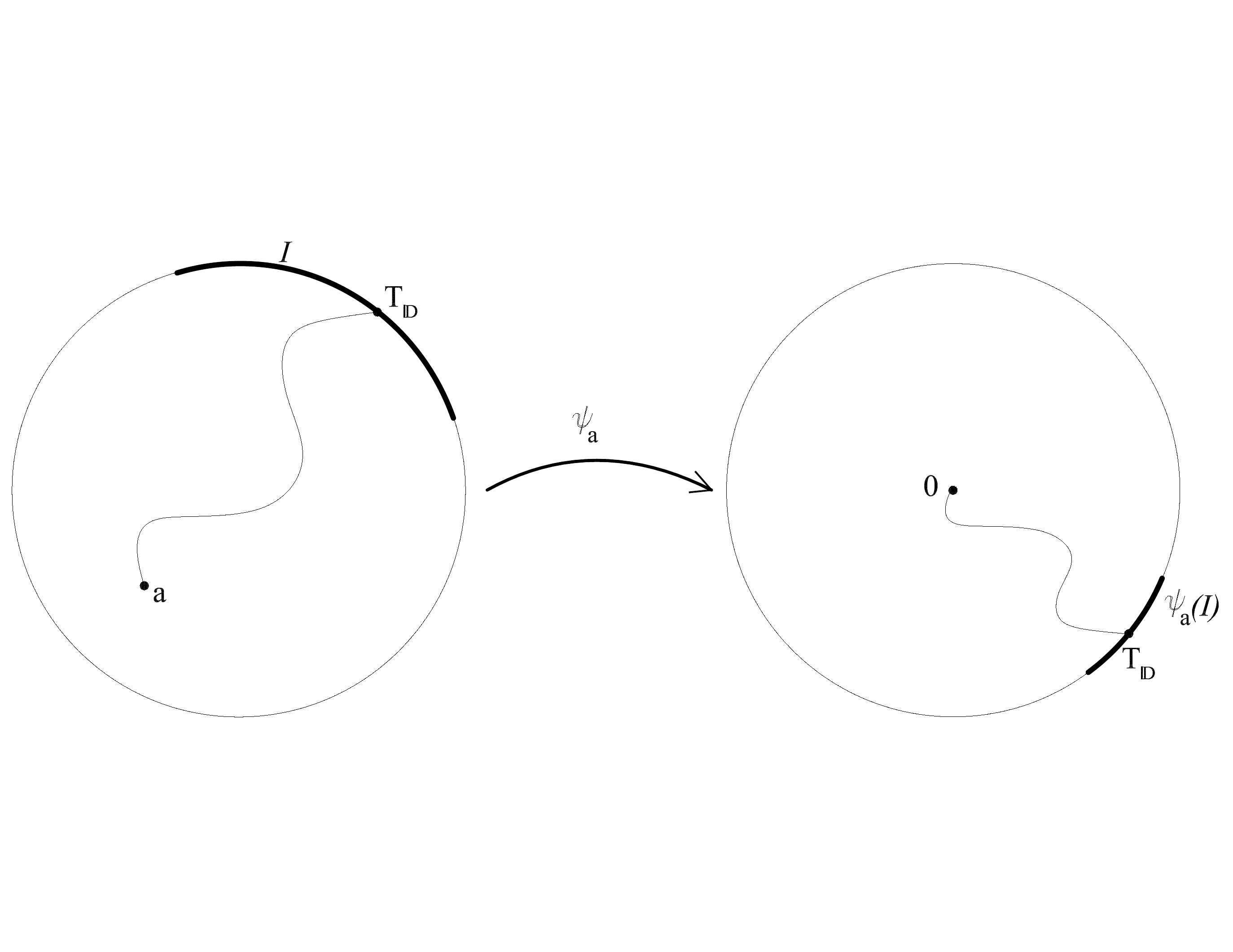}

We obtain the following identity:

\begin{equation} \label{}
P_a (B_{T(\DD)} \in \II) = P_0 (B_{T(\DD)} \in \psi_a(\II)) = \frac{m(\psi_a(\II))}{2\pi}.
\end{equation}
To calculate $m(\psi_a(\II))$, note simply that the map $\psi_a$ expands an infinitesimal portion of $\II$ centered at a point $e^{it} \in \II$  with expansion factor $|\psi_a'(e^{it})|$. It follows that

\begin{equation} \label{}
P_a (B_{T(\DD)} \in \II) = \frac{m(\psi_a(\II))}{2\pi} = \frac{1}{2\pi}\int_{\II} |\psi_a'(e^{it})| dt = \frac{1}{2\pi}\int_{\II} \frac{1-|a|^2}{|1-\bar a e^{it}|^2} dt.
\end{equation}
In the next section we will explore applications of this formula for harmonic and analytic functions.

\section{Poisson's integral formula and power series} \lll{pifseries}

Let us suppose that we have a harmonic function $h$ on $\DD$ which is continuous on $\seg{10}{\DD}$. We then have by Dynkin's formula, for any $a \in \DD$,

\begin{equation} \label{}
h(a) = E_a[B_{T_\DD}] = \int_{0}^{2 \pi} h(e^{i t}) P_a (B_{T(\DD)} \in dt).
\end{equation}

Recalling the formula for exit probabilities from the previous section, we let $\II$ shrink to an infinitesimal interval to obtain the identity

\begin{equation} \label{avgn}
P_a (B_{T(\DD)} \in dt) = \frac{1}{2 \pi}|\psi_a'(e^{it})| dt = \frac{1}{2 \pi} \frac{1-|a|^2}{|1-\bar a e^{it}|^2} dt.
\end{equation}

Let $\PP_r(\th) = \frac{1-r^2}{|1-re^{-i\th}|^2}$; this is known as the {\it Poisson kernel}. If we set $a=re^{i\th}$, then \rrr{avgn} gives us

\begin{theorem}[Poisson integral formula] \label{PIF}

Let $h$ be a harmonic function on $\DD$ which is continuous on $\seg{10}{\DD}$. Then

\begin{equation} \label{PIFform}
h(re^{i\th}) = \frac{1}{2\pi} \int_{0}^{2\pi}  h(e^{it}) \PP_r(\th - t)  dt .
\end{equation}
\end{theorem}

%The quantity $\PP_r(\th) = \frac{1-r^2}{1-2 r \cos \th +r^2}$ is known as the  The Poisson integral formula, which we will abbreviate to PIF, can therefore be expressed as

%\begin{equation} \label{}
%h(re^{i\th}) = \frac{1}{2\pi} \int_{0}^{2\pi}  h(e^{it}) \PP_r(\th - t) dt.
%\end{equation}
%
There are several equivalent forms for the Poisson kernel (see Exercise \ref{Poisskerequiv} below), one of which is $\PP_r(\th) = \sum_{n=-\ff}^{\ff} r^{|n|} e^{in\th}$; this is an extremely useful fact, since it allows us to write

\begin{equation} \label{jimbean}
\begin{split}
h(a) & = h(re^{i\th}) =  \frac{1}{2\pi} \int_{0}^{2\pi} h(e^{it}) \Big(\sum_{n=-\ff}^{\ff} r^{|n|} e^{in(\th-t)}\Big) dt = \\
& = \sum_{n=-\ff}^{\ff} r^{|n|} e^{in \th}  \frac{1}{2\pi} \int_{0}^{2\pi} h(e^{it}) e^{-int} dt \\
& = \bb_0 + \sum_{n=1}^{\ff} \bb_{-n} \bar a^n + \sum_{n=1}^{\ff} \bb_{n} a^n,
\end{split}
\end{equation}
where

\begin{equation} \label{sercoeffform}
\bb_n = \frac{1}{2\pi} \int_{0}^{2\pi} h(e^{it}) e^{-int} dt.
\end{equation}
Note that $h$ is continuous on $\bar \DD$ and thus bounded there, so that interchanging summation and integration in \rrr{jimbean} is justified. We see that $h$ can be expressed as the sum of power series in $a$ and $\bar a$. There is nothing special about the unit disk for the purposes of this argument; it applies equally well to disks with arbitrary centers and radii, and we therefore obtain

\begin{theorem} \label{harmseries}
Let $h(z)$ be harmonic (or analytic) on $D(z_0, R)$. Then we may write

\begin{equation} \label{harmfunexp}
h(z) = \bb_0 + \sum_{n=1}^{\ff} \bb_{-n} \seg{30}{(z-z_0)}^n + \sum_{n=1}^{\ff} \bb_{n} (z-z_0)^n,
\end{equation}
for all $z \in D(z_0, R)$, with $\bb_n$ given by

\begin{equation} \label{sercoeffform}
\bb_n = \frac{1}{2\pi r^n} \int_{0}^{2\pi} h(z_0 + r e^{it}) e^{-int} dt,
\end{equation}
for any $0<r<R$. These series converge uniformly on $D(z_0,R')$ for any $R' < R$. Furthermore, the coefficients in the series are unique, and $h$ is real valued if and only if $\bb_{-n} = \bar \bb_n$ for all $n$.
\end{theorem}

{\bf Proof:} The first statement has already been proved, and the uniform convergence follows as in Proposition \ref{redrum}. The final statement follows from a standard technique in Fourier analysis; see Exercise \ref{uniq} below. \qed

%???Need something about radius of convergence, $\limsup |a_n|^{1/n}$, etc. ???

\vski

We will somewhat informally refer to the first sum in \rrr{harmfunexp} as the {\it $\bar z$-part}, and the second as {\it $z$-part}, of the power series for $h(z)$. Since $z$ and $\bar z$ are simply linear combinations of $x$ and $y$, we obtain

\begin{corollary} \label{}
If $h$ is harmonic or analytic, then $h$ is $C^\ff$.
\end{corollary}

Define the differential operators

\begin{equation} \label{}
\frac{d}{dz} = \frac{1}{2} (\frac{d}{dx} - i \frac{d}{dy}), \qquad \frac{d}{d\bar z} = \frac{1}{2} (\frac{d}{dx} + i \frac{d}{dy}).
\end{equation}
These are known as the {\it Wirtinger derivatives}, and satisfy a number of nice properties, such as

\begin{lemma} \label{wirtprop} Suppose $f,g$ are differentiable functions, and $a,b \in \CC$. Then
\begin{itemize} \label{}

\item[(i)] (Sum/difference rule)

\begin{equation} \label{}
\frac{d}{dz} (a f + b g) = a \frac{df}{dz} + b \frac{dg}{dz}, \qquad \frac{d}{d\bar z} (a f + b g) = a \frac{df}{d\bar z} + b \frac{dg}{d\bar z}.
\end{equation}

\item[(ii)] (Product rule)

\begin{equation} \label{}
\frac{d}{dz} (f g) = f \frac{dg}{dz} + g \frac{df}{dz}, \qquad \frac{d}{d\bar z} (f g) = f \frac{dg}{d\bar z} + g \frac{df}{d\bar z}.
\end{equation}

\item[(iii)] (Chain rule)

\begin{equation} \label{}
\frac{d}{dz} (f \circ g) =  \Big(\frac{df}{dz} \circ g \Big) \frac{dg}{dz} + \Big(\frac{df}{d\bar z} \circ g \Big) \frac{d \bar g}{dz}, \qquad \frac{d}{d\bar z} (f \circ g) =  \Big(\frac{df}{dz} \circ g \Big) \frac{dg}{d \bar z} + \Big(\frac{df}{d\bar z} \circ g \Big) \frac{d \bar g}{d\bar z}.
\end{equation}

\item[(iv)] For integer $n \geq 1$, we have

\begin{equation} \label{}
\frac{d}{dz} z^n = nz^{n-1}, \qquad \frac{d}{d\bar z} \bar z^n = n \bar z^{n-1}, \qquad \frac{d}{dz} \bar z^n = \frac{d}{d \bar z} z^n = 0.
\end{equation}

\item[(v)] If $h$ is $C^2$, then

\begin{equation} \label{}
\Del h = 4\frac{d^2 h}{dz d\bar z} = 4\frac{d^2 h}{d \bar z dz}.
\end{equation}

\end{itemize}
\end{lemma}

The properties given in this lemma are all consequences of standard rules in multivariate calculus, and their proofs are left to the reader (Exercise \ref{exwirtprop}). An immediate consequence of the lemma is the following.

\begin{proposition} \label{wertderseries}
Suppose $h(z)$ is harmonic (or analytic) on $D(z_0, R)$, and given by \rrr{harmfunexp}. Then

\begin{equation} \label{treeroll}
\frac{dh}{dz} = \sum_{n=1}^{\ff} n \bb_{n} (z-z_0)^{n-1}, \qquad \frac{dh}{d \bar z} = \sum_{n=1}^{\ff} n \bb_{-n} \seg{30}{(z-z_0)}^{n-1}.
\end{equation}
\end{proposition}

On the other hand, we are now in a position to prove

\begin{proposition} \label{itsharm}
If $h(a) = E_a[h(B_\tau)]$ for any appropriate stopping time $\tau$, then $h$ is harmonic.
\end{proposition}

{\bf Proof:} The fact that $h(a) = E_a[h(B_\tau)]$ means that the proof of Poisson's Integral formula \rrr{PIFform} applies, and in particular $h$ can be expressed in a series expansion of the form \rrr{harmfunexp}. It is now clear that $\Del h = 4\frac{d^2 h}{dz d\bar z} = 0$ \qed

The reader will note that the condition that $h$ is analytic is relegated to parentheses in the statements of Theorem \ref{harmseries} and Proposition \ref{wertderseries}; this is to signify that, while these results do hold for analytic functions, in fact much more is true. If $f$ is analytic at $z_0$, then the limit in \rrr{deriv} must exist, regardless of the manner in which $\Del z \lar 0$. In particular, we may let $\Del z = x$ be real, separate $f$ into real and imaginary parts $u + iv$, and obtain

\begin{equation} \label{cr1}
\begin{split}
f'(z_0) & = \lim_{x \lar 0} \frac{f(z_0+x)-f(z_0)}{x} = \lim_{x \lar 0} \frac{u(z_0+x)-u(z_0)}{x} + \lim_{x \lar 0} \frac{v(z_0+x)-v(z_0)}{x} \\
& = \frac{du}{dx} (z_0) + i \frac{dv}{dx} (z_0).
\end{split}
\end{equation}
We may then conduct the similar operation where $\Del z = iy$ is purely imaginary, to obtain

\begin{equation} \label{cr2}
\begin{split}
f'(z_0) & = \lim_{y \lar 0} \frac{f(z_0+yi)-f(z_0)}{yi} = \lim_{y \lar 0} \frac{u(z_0+yi)-u(z_0)}{yi} + \lim_{y \lar 0} \frac{v(z_0+yi)-v(z_0)}{yi} \\
& = -i \frac{du}{dy} (z_0) + \frac{dv}{dy} (z_0).
\end{split}
\end{equation}

Equating real and imaginary parts of \rrr{cr1} and \rrr{cr2} gives the {\it Cauchy-Riemann equations}:

\begin{equation} \label{}
\frac{du}{dx} = \frac{dv}{dy}, \qquad \frac{du}{dy} = - \frac{dv}{dx}.
\end{equation}
Analytic functions satisfy these equations at all points in their domains. Note that we may write

\begin{equation} \label{}
\frac{df}{d\bar z} = \frac{du}{d\bar z} + i \frac{dv}{d\bar z} = \frac{1}{2}\Big(\frac{du}{dx} - \frac{dv}{dy}\Big) + \frac{i}{2} \Big(\frac{du}{dy} + \frac{dv}{dx}\Big) = 0.
\end{equation}
Thus, the second series in \rrr{treeroll} is identically 0 for analytic functions, which means that all coefficients $\bb_{-n}$ must be zero as well. Theorem \ref{harmseries} and Proposition \ref{wertderseries} are therefore eclipsed for analytic functions by

\begin{theorem} \label{holseries}
Let $f(z)$ be analytic on $D(z_0, R)$. Then we may write

\begin{equation} \label{holfunexp}
f(z) = \sum_{n=0}^{\ff} \bb_{n} (z-z_0)^n.
\end{equation}
This series converges uniformly on $D(z_0,R')$ for any $R' < R$. Furthermore, $f'$ is analytic as well, and we have

\begin{equation} \label{}
f'(z) = \frac{df}{dz} = \sum_{n=1}^{\ff} n \bb_{n} (z-z_0)^{n-1}.
\end{equation}
We also have $\frac{df}{d \bar z} = 0$.
\end{theorem}

{\bf Remark:} The relationship between harmonic and analytic functions can now be clarified. If $f$ is analytic then the series expansion \rrr{holfunexp} holds, and thus $Re(f) = \frac{f+\bar f}{2}$ and $Im(f) = \frac{f-\bar f}{2i}$ are both power series in $z$ and $\bar z$ and are therefore harmonic. On the other hand, if$f=u+iv$ is $C^2$ and satisfies the Cauchy-Riemann equations, then

\begin{equation} \label{}
\frac{d^2u}{dx^2} + \frac{d^2u}{dy^2} = \frac{d^2v}{dydx} - \frac{d^2v}{dxdy} = 0, \qquad \frac{d^2v}{dx^2} + \frac{d^2v}{dy^2} = -\frac{d^2u}{dydx} + \frac{d^2u}{dxdy} = 0,
\end{equation}
so that $u$ and $v$ are each individually harmonic. $f$ is therefore harmonic and admits the expansion given in \rrr{harmfunexp}, but with all coefficients $\bb_{-n}$ equal to 0 by the Cauchy-Riemann equations; it is therefore analytic.

%On the other hand, analytic functions are clearly harmonic, but if $f=u+iv$ is analytic then in fact we have by the Cauchy-Riemann equations

%as well. We mention finally that , then it can be shown that $f$ must necessarily be analytic; we will not make use of this fact, and therefore do not prove it, but it can be found in \cite[Thm. 2.29]{conway}.

\section{Exercises}

\ccases{ex_unbddom} Show that Theorem \ref{itothm} fails for the function $u(z)=u(x + yi)=y$ for the stopping time $T_K$ where $K=\{y=-1\}$.

\vski

\ccases{ex_nostoppingtime} Let $T_{k} = \inf\{t \geq 0: Im(B_t) = k\}$ and let $\hat T_{-1} = \sup \{ 0 \leq t \leq T_{-2} : Im(B_t) = -1\}$. Let $\tau = T_1 \wedge \hat T_{-1}$. Note that $E[\tau] < \ff$ but $\tau$ is not a stopping time. Show that Dynkin's formula fails for $\tau$ with the function $u(x+yi) = y$.

\vski

%\ccases{fintimehit} Let $\ga_r = \{|z|=r\}$.
%a) Show that $E_0[T_{\ga_r}] < \ff$. (Need hint here)
%b) Show that $E_0[T_{\ga_r}] = \frac{r^2}{2}$. (Hint: Apply Dynkin's formula with the non-harmonic function $u(x+iy) = x^2 + y^2$.

\ccases{mobbi} Show that M\"obius transformations are bijections from the Riemann sphere to itself. Check also that they are analytic at $\ff$, and at $-d/c$ when $c \neq 0$.

\vski

\ccases{mobmat} $SL_2(\CC)$ is the space of all $2 \times 2$ matrices with complex entries and determinant 1. This is a group under matrix multiplication. The {\it M\"obius group} is the set of all M\"obius transformations $\frac{az + b}{cz + d}$ with $ad-bc=1$. This is a group under composition. Show that the map $\Theta$ defined by $\Theta(\frac{az + b}{cz + d}) = \left( \begin{array}{ccc}
a & b  \\
c & d \end{array} \right)$ is an isomorphism.

\vski

\ccases{diskauto} Let $\psi_a(z) = \frac{z-a}{1-\bar a z}$ with $|a|<1$. Show that $|\psi_a(z)| \leq 1$, with equality precisely when $|z| = 1$.

\vski

%\ccases{upperhpPIF} $a)$ Show that the M\"obius transformation $\phi(z) = \frac{z-i}{z+i}$ maps the upper half-plane conformally to  Suppose that PIF for upper half plane.

%\vski

\ccases{Poisskerequiv} Prove the following equivalent forms for the Poisson kernel:

\begin{equation} \label{}
\PP_r(\th) = \frac{1-r^2}{1-2 r \cos \th +r^2} = \sum_{n=-\ff}^{\ff} r^{|n|} e^{in\th} = Re \Big( \frac{1+re^{i\th}}{1-re^{i\th}} \Big)
\end{equation}

\vski

\ccases{exwirtprop} Verify $(ii)-(v)$ of Lemma \ref{wirtprop}.

\vski

\ccases{mobcircex} Prove Proposition \ref{mobcirc} by employing the result of Exercise \ref{stereoflip} from Chapter \ref{chzero}.

\vski

\ccases{uniq}

$a)$ Show that, if $r$ is an integer, then

\begin{equation} \label{}
\int_{0}^{2\pi} e^{irt} dt = \left \{ \begin{array}{ll}
2 \pi & \qquad  \mbox{if } r=0  \\
0 & \qquad \mbox{if } r \neq 0 \;
\end{array} \right.
\end{equation}

$b)$ Suppose

\begin{equation} \label{harmfunexp}
h(z) = \bb_0 + \sum_{n=1}^{\ff} \bb_{-n} \seg{30}{(z-z_0)}^n + \sum_{n=1}^{\ff} \bb_{n} (z-z_0)^n,
\end{equation}
for all $z \in D(z_0, R)$. Use part (a) to show that

\begin{equation} \label{sercoeffform}
\bb_n = \frac{1}{2\pi r^n} \int_{0}^{2\pi} h(z_0 + r e^{it}) e^{-int} dt.
\end{equation}

Show also that $h$ is real-valued if, and only if, $\bb_{-n} = \bar \bb_n$ for all $n$.

\vski

\ccases{antianal} A function $f$ is {\it anti-analytic} if $\bar f$ is analytic. Show that if $f$ is a nonconstant, anti-analytic function, then $f(B_t)$ is a time changed Brownian motion, and identify the time change.

\vski

%\ccases{exCRconv} Show that a function which satisfies the CR equations is analytic.

%\vski

\ccases{concirc} Suppose that $\CCC_1, \CCC_2$ are two circles on the sphere which do not intersect (so, they may be two circles in the plane or one line and one circle). Show that there is a M\"obius transformation taking them to two concentric circles (i.e. circles with the same center) in the plane. (Solution: Take one circle to the imaginary axis, and then translate and scale so that the other has center on $\RR$ and intersects $\RR$ at $\frac{1}{r}$ and $r$. Then $\phi(z) = \frac{z-1}{z+1}$ does the trick.

\setcounter{cccases}{0}
\chapter{Examples and basic properties of harmonic and analytic functions}
\label{altrings} % So I can \ref{altrings} later.
\section{General properties of harmonic functions} \lll{logargsubs}

We begin by studying $h(z) = \ln |z| = \frac{1}{2} \ln |z|^2 = \frac{1}{2} \ln(x^2+y^2)$. A straightforward calculation reveals $h$ to be harmonic (showing $\frac{d^2h}{dx^2} + \frac{d^2h}{dy^2} = 0$ is not difficult, but the slick way is to calculate $\frac{dh}{dzd\bar z}$ and use Lemma \ref{wirtprop}). Thus, for $r < R$ the function $h_r^R(z) = \frac{\ln R-\ln |z|}{\ln R-\ln r}$ is harmonic on $A= A_{r,R} := \{r < |z| < R\}$, and it is equal to 1 on $\{|z|=r\}$ and 0 on $\{|z|=R\}$. Fix $a$ with $r < |a| < R$; by It\=o's formula we have $h_r^R(a) = E_a[h_r^R(B_{T_A})] = P(|B_{T_A}| = r)$. We therefore obtain the probability that $B_t$ exits $A$ at the inner circle $\{|z|=r\}$ when started at $a$. We can now derive two important facts about planar Brownian motion by this observation. The first to notice is that, upon letting $R \nearrow \ff$, we see that $h_r^R(a) \nearrow 1$; this tells us that $B_t$ is {\it recurrent}: regardless of the starting point, $B_t$ will hit any open set eventually with probability $1$. On the other hand, if we hold $R$ fixed and let $r \searrow 0$, then $h_r^R(a) \searrow 0$ as well; this tells us that $P_a(B_t = 0 \mbox{ for some } t>0) = 0$, and the same holds if 0 is replaced by any $b \neq 0$. This may seem counterintuitive, since after all $P(\cup_{z \in \CC} \{B_t = z \mbox{ for some } t>0\}) = P(B_t \mbox{ hits some point in } \CC) = 1$, but since there are uncountably many points in $\CC$ we do not run into a contradiction by summing probabilities. This curious property of Brownian motion is commonly expressed by a phrase like, "planar Brownian motion is recurrent but doesn't see points". We also can note that this phenomenon sheds new light the holomorphic invariance of Brownian motion, Theorem \ref{holinv}, as the reader may have noticed a potential problem with the intuitive proof given for that result: a holomophic function is not a local dilation at a point at which its derivative is 0. However, we will see that the zero set of a holomorphic function is a countable set, and therefore Brownian motion avoids the zero set of the derivative with probability 1, and again we do not run into trouble.

\vski

Let us calculate the power series for $h(z)$ around $z=1$. We have

\begin{equation} \label{}
\frac{d}{dz} \Big(\frac{1}{2}\ln z \bar z \Big)= \frac{1}{2z} = \frac{1}{2(1+(z-1))} = \frac{1}{2}(1 - (z-1) + (z-1)^2 - \ldots ).
\end{equation}

Proposition \ref{wertderseries} tells us that integrating this series will give us the $z$-part of the series for $h$, while the final statement of Theorem \ref{harmseries} and the fact that $h$ is real valued tell us that the $\bar z$-part is obtained by conjugation. Using $h(1) =0$, we arrive at the following expression:

\begin{equation} \label{}
h(z) = \ln |z| = \frac{1}{2} \sum_{n=1}^{\ff}\frac{(-1)^{n-1}}{n}\seg{30}{(z-1)}^n + \frac{1}{2} \sum_{n=1}^{\ff}\frac{(-1)^{n-1}}{n}(z-1)^n.
\end{equation}
The reader will note that this reduces to the familiar series $\ln x = \sum_{n=1}^{\ff}\frac{(-1)^{n-1}}{n}(x-1)^n$ when $z=x$ is real.

\vski

Let us now consider the function $\tilde h(z) = \Arg z = \tan^{-1} \frac{y}{x} = \tan^{-1} -i (\frac{z-\bar z}{z+\bar z})$, defined on $\CC \backslash (-\ff,0]$. Another straightforward calculation shows that $\Del \tilde h = 4 \frac{d^2 \tilde h}{dzd\bar z} = 0$, and $\tilde h$ is harmonic. Suppose that $-\pi < \th_1 < \th_2 < \pi$, and suppose $\th_1 < \Arg a < \th_2$. If we set $\tilde h_{\th_1}^{\th_2}(z) = \frac{\th_2 - \Arg z}{\th_2 - \th_1}$, we see that $\tilde h_{\th_1}^{\th_2}$ is equal to 1 on $\{\Arg z = \th_1\}$ and equal to 0 on $\{\Arg z = \th_2\}$. Thus, if we set $W=W_{\th_1}^{\th_2} = \{\th_1 < \Arg z < \th_2\}$, then $\tilde h_{\th_1}^{\th_2}(a) = E_a[h_{\th_1}^{\th_2}(B_{T_W})] = P_a(\Arg(B_{T_W}) = \th_1)$ (there is a subtlety here: since the domain $W$ is unbounded, there are many harmonic functions with the given boundary values; however, in Chapter \ref{exitime} we will see that there is only one such function which remains bounded, and this must then be equal to $P_a(\Arg(B_{T_W}) = \th_1)$). That is, we obtain the probability that $B_t$ hits the set $\{\Arg z = \th_1\}$ before hitting $\{\Arg z = \th_2\}$ when started at $a$, and this probability is linear in $\Arg a$.

\vski

Let us compute the power series for $\tilde h$ around $z=1$. We have

\begin{equation} \label{}
\frac{d}{dz} \Big( \tan^{-1} i (\frac{z-\bar z}{z+\bar z}) \Big)= \frac{-i}{2z} = \frac{-i}{2(1+(z-1))} = \frac{-i}{2}(1 - (z-1) + (z-1)^2 - \ldots ).
\end{equation}
As before, Proposition \ref{wertderseries} tells us that the $z$-part of the series for $\tilde g$ comes by integrating this series, while the $\bar z$-part is again obtained by conjugation. Again we have $\tilde h(1) =0$, so we obtain

\begin{equation} \label{}
\tilde h(z) = \Arg z = \frac{i}{2} \sum_{n=1}^{\ff}\frac{(-1)^{n-1}}{n}\seg{30}{(z-1)}^n + \frac{-i}{2} \sum_{n=1}^{\ff}\frac{(-1)^{n-1}}{n}(z-1)^n.
\end{equation}

This looks suspiciously like the series for $\ln |z|$, and we arrive at an important idea. We note that we can write

\begin{equation} \label{}
h(z) + i\tilde h(z) = \ln |z| + i\Arg z =  \sum_{n=1}^{\ff}\frac{(-1)^{n-1}}{n}(z-1)^n,
\end{equation}
and Theorem \ref{holseries} tells us that we have a holomorphic function. Thus, the seemingly unrelated harmonic functions $h$ and $\tilde h$ have an intimate relationship with each other, namely that they can be taken to be the real and imaginary parts of a holomorphic function. We will call any such pair {\it harmonic conjugates}. The reader may note a seeming incongruity with this example, however, which is that the functions $h$ and $\tilde h$ do not have the same domains of definition: $h$ is harmonic on $\CC \bsh \{0\}$, while $\tilde h$ can only be defined on $\CC \bsh (-\ff. 0]$, as the function $\Arg z$ cannot be made to be continuous on $(-\ff. 0]$.

\vski

We are led immediately to a fundamental question: suppose $h(z)$ is any real-valued harmonic function; can we find another real-valued harmonic function $\tilde h(z)$ such that $h+i\tilde h$ is analytic? We note that the previous example shows that this is not always the case for all domains, however we can show easily using power series that it is the case when the domain of definition is a disk.

\begin{theorem} \label{}
Suppose $h(z)$ is a harmonic function on $D(z_0,R)$. Then there is another harmonic function $\tilde h(z)$ such that $h+i\tilde h$ is holomorphic on $D(z_0,R)$. This function is unique up to adding a constant.
\end{theorem}

{\bf Proof:} Although the argument is implicit in the example of the logarithm, we prove it in generality here. We know that we may write

\begin{equation} \label{harmfunexp91}
h(z) = \bb_0 + \sum_{n=1}^{\ff} \bb_{-n} \seg{30}{(z-z_0)}^n + \sum_{n=1}^{\ff} \bb_{n} (z-z_0)^n.
\end{equation}
$\tilde h$ must now be chosen such that the power series for $h+i\tilde h$ contains no $\bar z$-part. If we write

\begin{equation} \label{harmfunexp92}
\tilde h(z) = \aa_0 + \sum_{n=1}^{\ff} \aa_{-n} \seg{30}{(z-z_0)}^n + \sum_{n=1}^{\ff} \aa_{n} (z-z_0)^n,
\end{equation}
we see that we must have

\begin{equation} \label{}
\sum_{n=1}^{\ff} \bb_{-n} \seg{30}{(z-z_0)}^n + i \sum_{n=1}^{\ff} \aa_{-n} \seg{30}{(z-z_0)}^n = 0,
\end{equation}
and we must therefore have $\aa_{-n} = i \bb_{-n}$. Having obtained this, the requirement that $\tilde h$ be real-valued forces $\aa_n = \bar \aa_{-n} = -i \bar \bb_{-n} = -i \bb_n$. We therefore obtain the unique (up to choice of $\aa_0$) harmonic conjugate as

\begin{equation} \label{harmfunexp93}
\tilde h(z) = \aa_0 + i \sum_{n=1}^{\ff} \bb_{-n} \seg{30}{(z-z_0)}^n - i \sum_{n=1}^{\ff} \bb_{n} (z-z_0)^n,
\end{equation}

and we have

\begin{equation} \label{}
h+i\tilde h = \bb_0 + \aa_0 + 2 \sum_{n=1}^{\ff} \bb_{n} (z-z_0)^n.
\end{equation}
\qed

We will see later that a disk may be replaced by any simply connected domain.

\vski

Let us now establish some other properties of harmonic functions. 

\begin{theorem}{\rm [Mean value property]}  \label{}
Suppose $h$ is harmonic on a domain $W$, and the set $\{|z-a| \leq r\} \subseteq W$. Then

\begin{equation} \label{}
h(a) = \frac{1}{2\pi} \int_{0}^{2\pi} h(a + re^{i\th}) d \th.
\end{equation}

\end{theorem}

{\bf Proof:} By Dynkin's formula, $h(a) = E_a[h(B_{T(D(a,r))}]$, and by the rotation invariance of Brownian motion the distribution of $B_{T(D(a,r))}$ on $\{|z-a| = r\}$ is uniform. The result follows. \qed

In other words, $h(a)$ is equal to the average of its values on a circle in $W$ centered at $a$. Note that this result can also be obtained from the Poisson Integral Formula in Chapter 1.
 
\begin{theorem}[Maximum/minimum principle] \label{harmmaxthm}
Suppose $W$ is a bounded domain. If $h(z)$ is harmonic on $W$ and continuous on $cl(W)$ then
$$
\inf_{w \in \dd W} h(w) \leq h(z) \leq \sup_{w\in \dd W} h(w)
$$
for all $z\in W$. Furthermore, if $h$ is nonconstant, then

$$
\inf_{w \in \dd W} h(w) < h(z) < \sup_{w\in \dd W} h(w)
$$

for all $z\in W$.
\end{theorem}

Note that if $W$ is the unit disk $\DD$, then this result follows from the Poisson Integral Formula \rrr{PIFform} by virtue of the fact that $\frac{1}{2\pi} \int_{0}^{2\pi}  \PP_r(t)  dt = 1$. The general result will be shown in the next section. 

\begin{theorem}[Uniqueness principle] \label{harmmaxthm}
Suppose $W$ is a bounded domain. If $h_1(z), h_2(z)$ are harmonic on $W$ and continuous on $cl(W)$ with $h_1(z) = h_2(z)$ for all $z \in \dd W$, then $h_1(z) = h_2(z)$ for all $z \in W$.
\end{theorem}

{\bf Proof:} Apply the maximum/minimum principle to the harmonic function $h_1 - h_2$. \qed

\begin{theorem} \label{harmsing}
If $h(z)$ is harmonic and bounded on $\DD^\times = \{0<|z|<1\}$, then we can define $h(0)$ such that $h(z)$ so extended is harmonic on all of $\DD$.
\end{theorem}

{\bf Proof:} Dynkin's formula implies that $h(a) = E_a[h(B_{T(\DD^\times)})]$. As was shown at the beginning of this chapter, Brownian motion does not see points, and thus $P_a(B_{T(\DD^\times)} =0) =0$. The distribution of $B_{T(\DD^\times)}$ is therefore the same as that of $B_{T(\DD)}$, and thus $h(a) = E_a[h(B_{T(\DD^\times)})]$. As was shown in Chapter 1, $h(z)$ can therefore be expressed as a power series in $z$ and $\bar z$, and this series provides the required extension to all of $\DD$. \qed

\section{The Dirichlet Problem}

The Dirichlet problem in $\CC$ is the following.

\vski

{\bf Dirichlet problem} {\it If $U$ is a bounded domain, and $h$ a continuous real-valued function on $\dd U$, to extend $h$ to a harmonic function on $U$ which is continuous on the closure of $U$.}

\vski

As we will see, in many cases the Dirichlet problem may be solved by running a Brownian motion. To be precise, we will show that $h(a) = E_a[h(B_{T(U)})]$ provides the desired solution whenever $U$ is a sufficiently nice domain. However, there are domains, such as $\DD^\times = \{0<|z|<1\}$, for which the Dirichlet problem is not solvable. To see this, suppose that $h$ is harmonic on $\DD^\times$, continuous on its closure, with $h(z) = 1$ for $|z|=1$ and $h(0) = 0$. By Theorem \ref{harmsing}, $h$ is in fact harmonic on all of $\DD$, and thus by Theorem \ref{harmmaxthm} $0=h(0) \geq \inf \{|z|=1\} h(z) = 1$, a contradiction. Nonetheless, we have the following classical result.

\begin{theorem} \label{bigguy}
The Dirichlet problem is solvable for any bounded simply connected domain $U$ in $\CC$.
\end{theorem}

{\bf Proof:} Since $U$ is bounded, $T(U) < \ff$ a.s. We can therefore extend $h$ from $\dd U$ to $cl(U)$ by setting $h(a) = E_a[h(B_{T(U)})]$, and we must show that $h$ so extended on $U$ and continuous on $cl(U)$. Harmonicity follows from the Strong Markov property: 

$$h(a) = E_a[h(B_{T(U)})] = E_a[E_{B_{T(D(a,r))}(\omega')}[h(B_{T(U)}(\omega))]] = \frac{1}{2\pi} \int_{0}^{2\pi} h(a+re^{it})dt;$$

note that in the final expectation in this string of equalities, the inner expectation is taken with respect to $\omega\in \Om$ while the outer one is taken with respect to $\omega' \in \Om'$, where the probability spaces $\Om$ and $\Om'$ are independent. To show continuity we must prove that a Brownian motion started near a boundary point of $U$ is likely to exit near that point, in particular

\begin{lemma} \label{key}
If $y$ is a boundary point of a simply connected domain $U$, then given any $\eps >0$ we can find $\dd>0$ such that $P_a(|B_{T(U)}-y|<\eps) > 1 - \eps$ whenever $|a-y| < \dd$.
\end{lemma}

The continuity of $h$ on $cl(U)$ now follows from from this by writing

$$E_a [h(B_{T(U)})] = E_a [h(B_{T(U)})1_{\{|B_{T(U)}-y|<\eps\}}] + E_a [h(B_{T(U)})1_{\{|B_{T(U)}-y| \geq \eps\}}]$$

and noting that for $a$ sufficiently close to $y$ the first term on the right is as close as we like to $h(y)$ and the second can be made as small as we like (using the boundedness of $h$).

\vski

{\bf Proof of Lemma \ref{key}:} Consider the domain $V_M = \{Re(z) < M\} \bsh L$, where $L$ is the ray in $\CC$ emanating from $2 \pi i$ horizontally to the left; see Figure \ref{figdir}. We note that if $M$ is sufficiently large then $P_0(B_{T(V_M)} \in L) > 1-\eps$; this is evident, since as $M \nearrow \ff$ we have $P_0(B_{T(V_M)} \in L) \nearrow P_0(T(V_\ff) < \ff)$, where $V_\ff = \CC \bsh L$, and $P_0(T(V_\ff) < \ff) = 1$ (proof: by L\'evy's Theorem, $arg(B_t-2\pi i)$ is a time-changed one-dimensional Brownian motion, is therefore recurrent, and thus is equal to $-\pi$ at some finite time a.s.). We may therefore choose $\dd$ sufficiently small so that $P_0(B_{T(V_M)} \in L) > 1-\eps$ when $M = \ln \frac{\eps}{2\dd}$, and ensure also that $\dd < \frac{\eps}{4}$. We claim now that $P_a(B_{T(U)} \in B(y,\eps)) > 1 - \eps$ whenever $|a-y| < \dd$.

\begin{center}
\begin{figure} \lll{figdir}
\vspace{-2in}
\includegraphics[width=5 in,height=7 in]{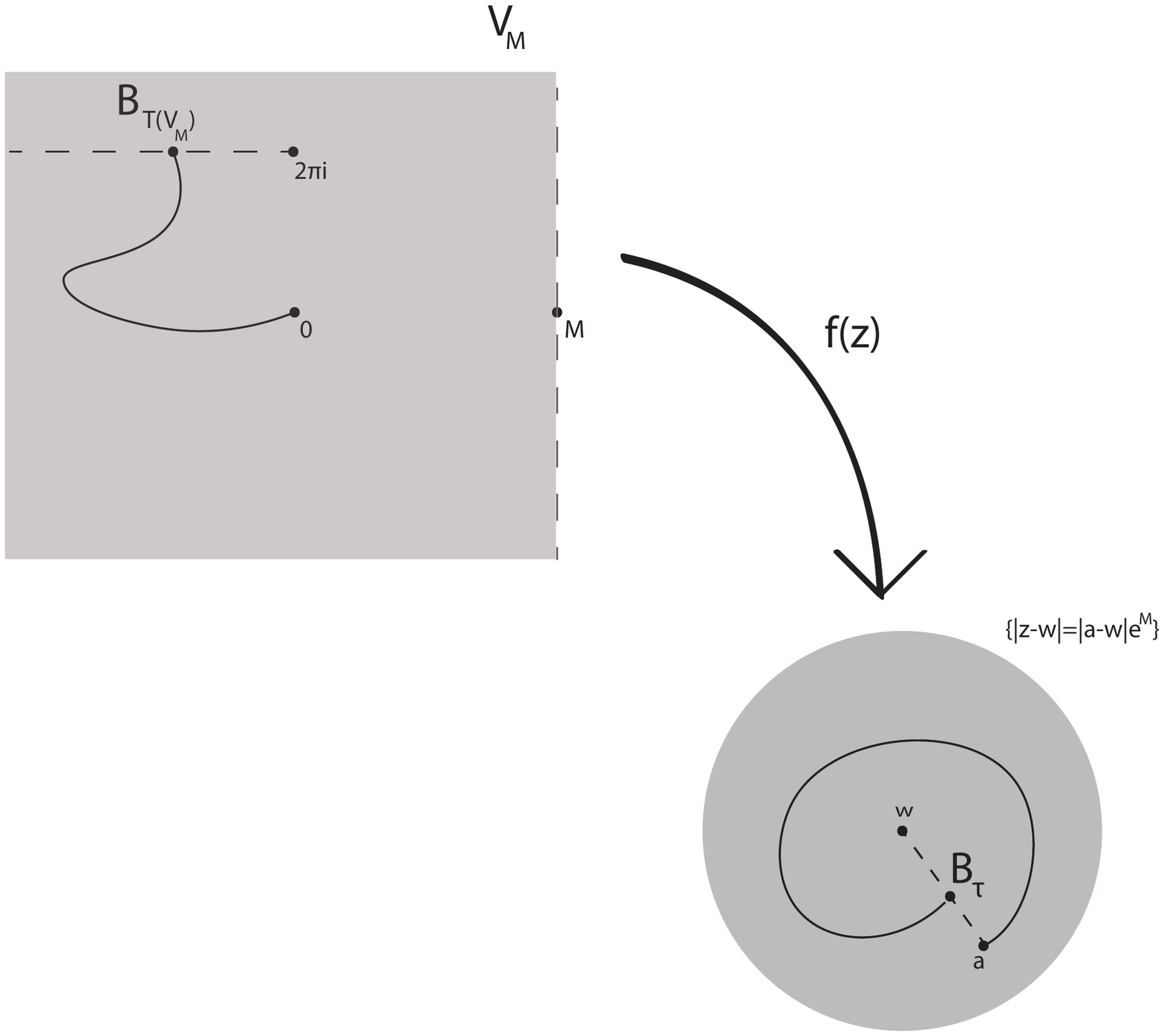}
\vspace{-1in}
\caption{The mapping by $f$ from $V_M$ to $\{|z-y|< e^M\}$.}
\end{figure}
\end{center}

\vspace{-.4in}

To prove this, suppose that $|a-y| < \dd$, and let $w$ be the point in $\dd U$ closest to $a$; note that $w$ is not necessarily the same as $y$, but $|a-w| < \dd$ which implies by the triangle inequality that $|y-w| < 2 \dd < \eps$. Consider now the analytic function $f(z) = (a-w)e^z + w$: it takes $0$ to $a$, takes $L$ to the line segment $[a,w)$ connecting $a$ and $w$, and wraps the vertical line $\{Re(z) = M\}$ around the circle $\{|z-w| = |a-w|e^M = (\frac{|a-w|}{\dd})\frac{\eps}{2}\}$. We will realize a Brownian motion starting at $a$ via L\'evy's Theorem as the image under $f$ of a Brownian motion starting at $0$. The projection $\tau = \si(T(V_M))$ under $f$ of the stopping time $T(V_M)$ is again a stopping time; we can describe it as "the first time when either $B_t$ lies on the line segment $(a,w)$ having wound once anticlockwise about $w$ or $|B_t-w| = (\frac{|a-w|}{\dd})\frac{\eps}{2}$". Note that $P_a(B_{\tau} \in (a,w)) = P_0(B_{T(V_M)} \in L) > 1-\eps$. Furthermore, on the event $\{B_{\tau} \in (a,w)\}$ we may connect $B_{\tau}$ to $a$ by a straight line segment and form a closed path which winds about the boundary point $w$; if this path remained entirely within $U$ this would contradict simple connectivity, and since $(a,w)$ lies within $U$ we conclude that $B_t$ intersects $\dd U$ before $\tau$. Thus, on the event $\{B_{\tau} \in (a,w)\}$ we have $T(U) < \tau \leq \inf\{t\geq 0:|B_t - w| = (\frac{|a-w|}{\dd}) \frac{\eps}{2}\}$. The triangle inequality implies easily that $\{|z-w| < (\frac{|a-w|}{\dd}) \frac{\eps}{2} \} \subseteq \{|z-w| < \frac{\eps}{2} \} \subseteq \{|z-y| < \eps \}$, and we conclude that $P_a(T(U) < \inf\{t\geq 0:|B_t - y| = \eps\})>1-\eps$. Since $\{T(U) < \inf\{t\geq 0:|B_t - y| = \eps\}\} \subseteq \{|B_{T(U)}-y|<\eps\}$, the lemma is proved. \qed

\vski

{\bf Remarks:} \begin{itemize} \label{}

\item The ray $L$ can be replaced in the proof by the union of the rays emanating to the left from all points of the form $n 2\pi i$ for integer $n \neq 0$. This would form an even smaller stopping time $\tau$, but would not change the proof in any meaningful way.

\item It may be interesting to note that the choice of $\dd$ depends only on $\eps$ and not on $U$.

\item The proof given here can be easily extended to other situations, such as when the complement of $U$ is formed of finitely many connected components, none of which is a single point.

\item The boundedness requirement on $U$ can in many cases be relaxed, but the Dirichlet problem and even the maximum/minimum and uniqueness principles for harmonic functions are rather subtle for unbounded domains. See Chapter \ref{exittimes} for more details on this.

\end{itemize}

\section{The Identity Theorem, Open Mapping Theorem, and Maximum Principle}

We now discuss a number of results on general analytic and harmonic function. We begin with a result, commonly referred to as the Identity Theorem, which states that if two analytic functions coincide on a set which contains an accumulation point than they agree everywhere. Note that this is not true for real-valued differentiable functions on the line, as for instance the function $f(x) = e^{-1/x}$ for $x>0$ and $f(x) = 0$ for $x \leq 0$ is $C^\ff$ on the real line, agrees with the constant function 0 on the negative real axis, and yet is clearly not equal to 0 for all $x$.

\begin{theorem} \label{identT}
Let $f$ be a nonconstant analytic function on a domain $W$, and let $\ZZ$ be the zero set of $f$; that is, let $\ZZ = \{z\in W: f(z) = 0\}$. Then $\ZZ$ is a discrete set in $W$.
\end{theorem}

{\bf Proof:} Let $z_0 \in \ZZ$. We know that in some disk centered at $z_0$ we can write $f(z) = \sum_{n=1}^{\ff} \bb_n (z-z_0)^n$. Let $N \geq 1$ be the smallest index so that $\bb_N \neq 0$; we can then write $f(z) = (z-z_0)^N \sum_{n=N}^{\ff} \bb_n (z-z_0)^{n-N}$. As $z \lar z_0$, $\sum_{n=N}^{\ff} \bb_n (z-z_0)^{n-N} \lar \bb_N \neq 0$, and it follows by continuity that this expression is nonzero on some (possibly smaller) disk centered at $z_0$. Since $(z-z_0)^N$ is only 0 at $z_0$, we see that $f$ is not zero in this disk except at $z_0$, and $z_0$ is therefore not an accumulation point of $\ZZ$. \qed

\begin{corollary} \label{identcor}
If $f, g$ are analytic on a domain $W$, and the set $\{z:f(z)=g(z)\}$ has an accumulation point in $W$, then $f(z) = g(z)$ for all $z$ in $W$.
\end{corollary}

{\bf Proof:} Apply Theorem \ref{identT} to the analytic function $f-g$. \qed

We now will prove that the image under a nonconstant analytic function of an open set is again open. This is commonly known as the {\it Open Mapping Theorem}. Note that, again, nothing like this is true for real-valued differentiable functions on the line, as for instance the image of the real line under the sine function is the closed interval $[-1,1]$. First, we need a lemma about Brownian motion.

\begin{lemma} \label{hitexit}
Let $a \in \CC$, and for any $\dd > 0$ let $\tau_\dd = \inf \Big\{t \geq 0 : B_t \in \{|z - a| = \dd\}\Big\}$. Fix $r>0$, and let $V$ be an open set contained in $D(a,r)$, where $D(a,r) = \{|z-a| < r\}$. Let $B_t$ be a Brownian motion starting at $a$. Then $P(B_t \in V \mbox{ for some } t \leq \tau_r) > 0$.
\end{lemma}

{\bf Proof:} We have

\begin{equation} \label{}
P_a(B_t \in V \mbox{ for some } t \leq \tau_r) \geq P_a(B_t \in \CCC \mbox{ for some } t \leq \tau_r) \geq P_a(B_{\tau_{r'}} \in \CCC) = \frac{\th_2 - \th_1}{2\pi}.
\end{equation}
\qed

\begin{theorem} \label{OMT}
Let $f$ be a nonconstant analytic function on a domain $W$. Then $f(W)$ is a domain.
\end{theorem}

{\bf Proof:} We must prove that $f(W)$ is open and connected. Analyticity is not needed for connectedness: any continuous image of a connected set is connected, since if $A$ and $B$ disconnect $f(W)$ then $f^{-1}(A)$ and $f^{-1}(B)$ disconnect $W$. Openness is a bit more tricky. Let $v \in f(W)$, and choose $a \in W$ such that $f(a) = v$. By Corollary \ref{identcor} we can choose a small $r>0$ such that $\bar D(a,r) \subseteq W$ and $f(w) \neq v$ on $\{|w-a|=r\}$. If we start a Brownian motion at $a$ and stop it at the hitting time $\tau$ of $\{|w-a|=r\}$, then $\hat B_t = f(B_{C_t})$ is a Brownian motion starting at $v$ and stopped at the stopping time $C^{-1}_\tau$, which satisfies $f(\hat B_{C^{-1}_\tau}) = f(B_{\tau}) \in f(\{|w-a|=r\})$. Let $\ga$ denote the curve $f(\{|w-a|=r\})$. In traveling from $v$ to $\ga$, $\hat B_t$ must cross first the circle $\{|w-v|=m\}$, where $m = \inf \{|w-v|: w \in \ga \} > 0$; so if we let $\hat \tau$ be the first time $\hat B_t$ hits $\{|w-v|=m\}$ then $C^{-1}_\tau \geq \hat \tau$ a.s.

\vspace{-1.5in}

\includegraphics[width=180mm,height=145mm]{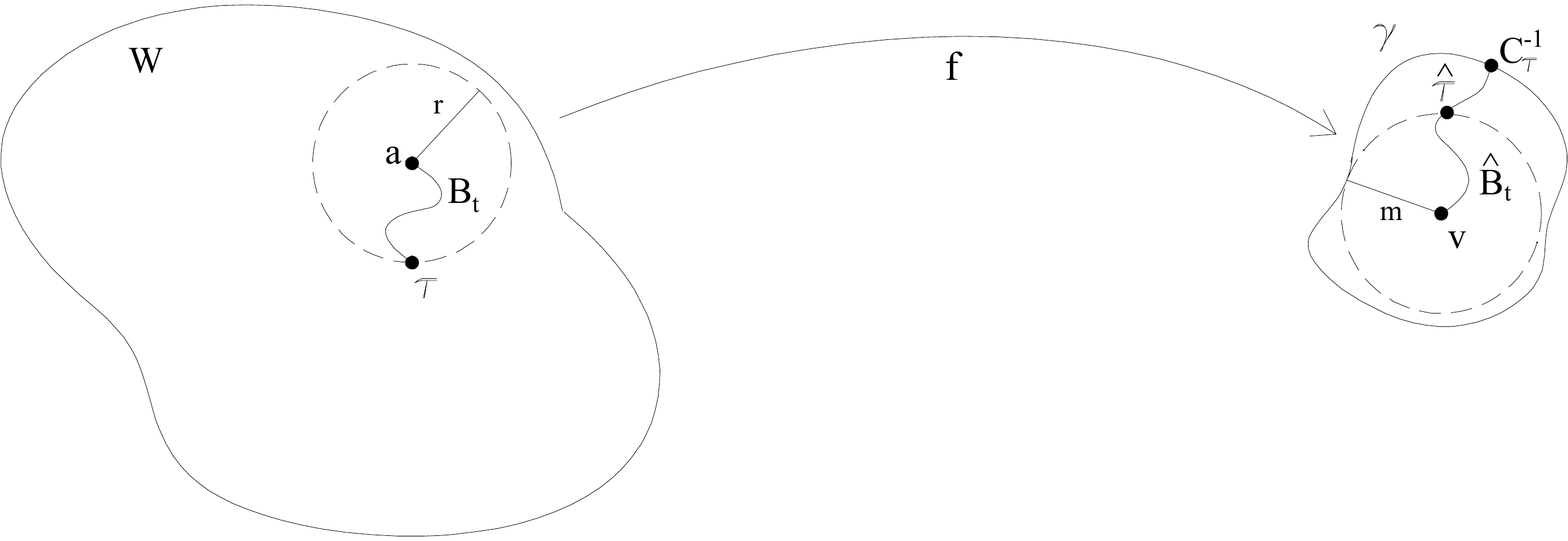}

\vspace{-1.55in}

Lemma \ref{hitexit} now shows that $\hat B_t = f(B_{C_t})$ hits any open set in $\{|w-v|<m\}$ with positive probability before time $\hat \tau$, and this shows in particular that $\{|w-v|<m\}$ is contained in the closure of $f(\bar D(z,r))$. However, $\bar D(z,r)$ is a compact set, so its image under the continuous map $f$ is compact and thus closed, and therefore $\{|w-v|<m\} \subseteq f(\bar D(z,r))\subseteq f(W)$. This proves that $f(W)$ is open. \qed

\begin{theorem} \label{maxprinanal}
Let $f$ be a nonconstant analytic function on a domain $W$. Then $f$ does not attain its maximum value on $W$; that is, there is no point $z$ such that $|f(z)| = \sup_{w \in W} |f(w)|$ for all $w \in W$.
\end{theorem}

{\bf Proof:} Let $\MM = \{z \in W: |f(z)| = \sup_{w \in W} |f(w)|\}$. $\MM$ is closed, being the preimage of a point under the continous map $f$; but we may argue that it is open as well, as follows. For $a \in \MM$, choose a small $r>0$ such that $\bar D(a,r) \subseteq W$, and let $\tau$ be the first time Brownian motion hits $\{|w-a|=r\}$. Then $|f(a)| = |E_{a}[f(B_\tau)]| \leq E_{a}[|f(B_\tau)|] \leq \sup_{|w-a|=r} |f(w)| \leq |f(a)|$. All inequalities in this string must then be equalities, which can only happen if $|f(w)| = |f(a)|$ almost surely on $\{|w-a|=r\}$. By continuity, $|f(w)| = |f(a)|$ for all $w$ on this circle, and by repeating this argument for any radii less than $r$ we see that $|f(w)| = |f(a)|$ for all $w \in D(a,r)$. This proves $\MM$ is an open set, and since it is both open and closed the connectedness of $W$ now proves that $\MM = \emptyset$ or $\MM = W$. If $\MM = W$ then $f(W) \subseteq \{|w| = const.\}$, but this implies by Theorem \ref{OMT} that $f$ is constant, a contradiction. Thus, $\MM = \emptyset$, proving the theorem. \qed

\begin{corollary} \label{}
Let $f$ be a nonconstant analytic function on a domain $W$ with $f(z) \neq 0$ on $W$. Then $f$ does not attain its minimum value on $W$; that is, there is no point $z$ such that $|f(z)| = \inf_{w \in W} |f(w)|$ for all $w \in W$.
\end{corollary}

{\bf Proof:} Apply Theorem \ref{maxprinanal} to the analytic function $1/f$. \qed

We should mention that the Maximum Principle applies to real-valued harmonic functions as well. The proof given here carries over directly, with the only difference being that dealing with the modulus is no longer necessary, and this removes the requirement that the function never be zero to deduce the Minimum Principle.

\begin{theorem} \label{}
Let $h$ be a nonconstant harmonic function on a domain $W$. Then $f$ does not attain its maximum (resp. minimum) value on $W$; that is, there is no point $z$ such that $h(z) = \sup_{w \in W} h(w)$ (resp. $h(z) = \inf_{w \in W} h(w)$) for all $w \in W$.
\end{theorem}

\section{The exponential and trigonometric functions}

Let us examine the exponential function $e^z = e^{x+yi} = e^x(\cos y + i \sin y)$ more closely. This is an analytic function mapping the plane to itself, with $|e^z| = e^x$ and $\Arg(e^z) = y$. We see that this function and the logarithm defined earlier in the chapter are inverses wherever $\Log$ can be applied, in the sense that $e^{\Log z} = e^{\ln |z| + i \Arg z} = |z| e^{i \Arg z} = z$, and $\Log e^z = \ln e^x + i \Arg (\cos y + i \sin y) = x + \tilde y i$, where $\tilde y$ is the value in $(-\pi, \pi]$ such that $y = \tilde y (\mbox{mod } 2\pi)$. Since the argument of $e^z$ can be any value, and the modulus of $e^z$ can be any positive value, the image of the plane under this map is $\CC \bsh \{0\}$. The vertical line $\{y = y_0\}$ in the plane gets mapped to a curve of constant argument; that is, to the infinite ray connecting 0 to $\ff$, which can be parameterized as $\{te^{i y_0}: 0 < t < \ff\}$. The horizontal line $\{x = x_0\}$ in the plane gets mapped to a curve of constant modulus; that is, to the circle centered at 0 of radius $e^{x_0}$. This geometric understanding of the function gives rise to several useful consequences for Brownian motion, as we now describe.

\vski

To begin with, it is a standard fact that one-dimensional Brownian motion is recurrent; in fact we showed in Section \ref{logargsubs} that planar Brownian motion is recurrent, and it is immediate from this that the one-dimensional variety is as well by considering the real or imaginary part of the planar case. However, there are proofs of recurrence in one dimension which do not use work done in the plane. We may therefore deduce planar recurrence by noting that recurrence on the line shows that the real part of a Brownian motion $B_t$ starting at any point $x_0 + i y_0$ will eventually hit any other point $x_1$, which means $B_t$ hits $\{x=x_1\}$ with probability 1. Taking the image under the exponential map, $B_t$ is mapped to a time changed Brownian motion $\tilde B_t$, and $\{x=x_1\}$ is mapped to $\{|z| = e^{x_1}\}$. We conclude that

\begin{equation} \label{}
P_{e^{x_0 + y_0 i}}(\tilde B_t \in \{|z| = e^{x_1}\} \mbox{ for some } t > 0) = 1,
\end{equation}
which is merely a restatement of recurrence. However, with a bit more work we can obtain an extension of this result which will be important to us. Let $r<R$ be two positive values, and define sequences of stopping times $\tau^{(r)}_n$ and $\tau^{(R)}_n$ as follows. Let $\tau^{(r)}_1 = \inf\{t>0: B_t \in \{|z| = r\}\}$, then $\tau^{(R)}_1 = \inf\{t>\tau^{(r)}_1: B_t \in \{|z| = R\}\}$. For larger $n$ define these similarly as $\tau^{(r)}_n = \inf\{t>\tau^{(R)}_{n-1}: B_t \in \{|z| = r\}\}$ and $\tau^{(R)}_n = \inf\{t>\tau^{(r)}_n: B_t \in \{|z| = R\}\}$. In words, $\tau^{(r)}_1$ is the first time $B_t$ hits $\{|z| = r\}$, then $\tau^{(R)}_1$ is the first time after $\tau^{(r)}_1$ that $B_t$ hits $\{|z| = R\}$, then $\tau^{(r)}_1$ is the first time after $\tau^{(R)}_1$ that $B_t$ hits $\{|z| = r\}$, etc. Recall Theorem \ref{StrongMP}, the Strong Markov Property, from Chapter \ref{chzero}, which states that a Brownian motion behaves identically following a stopping time to a Brownian motion starting at time 0; together with the recurrence already established this guarantees that the random variables $\tau^{(r)}_1, \tau^{(R)}_n - \tau^{(r)}_n,$ and $\tau^{(r)}_n - \tau^{(R)}_{n-1}$ are all finite almost surely, regardless of the starting point of $B_t$. This in turn guarantees the almost sure finiteness of $\tau^{(r)}_n$ and $\tau^{(R)}_n$ for all n; that is, $B_t$ will make infinitely many trips back and forth between $\{|z| = r\}$ and $\{|z| = R\}$. Let a small $\dd>0$ be given, and apply our knowledge with $r = \dd/2, R= \dd$. Note that for the duration of a trip between $\{|z| = \dd/2\}$ and $\{|z| = \dd\}$ the Brownian motion $B_t$ must reside entirely within the disk $D(0,\dd)$. Thus, we have

\begin{equation} \label{}
\int_{0}^{\ff} 1_{D(0,\dd)}(B_s) ds \geq \sum_{n=1}^{\ff} (\tau^{(\dd)}_n - \tau^{(\dd/2)}_n).
\end{equation}
Set $Y_n = (\tau^{(\dd)}_n - \tau^{(\dd/2)}_n)$. It is not hard to see by the Strong Markov Property and rotation invariance that the $Y_n$'s form a sequence of independent, identically distributed, positive random variables. We therefore have, for any $a \in \CC$,

\begin{equation} \label{}
E_a[e^{-\int_{0}^{\ff} 1_{D(0,\dd)}(B_s) ds}] \leq E_a[ e^{-\sum_{n=1}^{\ff} Y_n}] = \lim_{n \lar \ff} E_a[e^{-Y_1}]^n = 0.
\end{equation}
This can only occur if $\int_{0}^{\ff} 1_{D(0,\dd)}(B_s) ds = \ff$ a.s. By translating if necessary any open set $A$ so that it contains $0$, we have proved the following.

\begin{proposition} \label{opencamp}
For Brownian motion starting at any point and for any open set $A$ in $\CC$, we have

\begin{equation} \label{}
\int_{0}^{\ff} 1_{A}(B_s) ds = \ff \mbox{ a.s.}
\end{equation}
\end{proposition}

That is, Brownian motion spends an infinite amount of time in any open set with probability 1. This will be important in the next chapter, when we study general functions analytic on all of $\CC$.

Now do winding and Spitzer.

\section{Entire and almost entire functions}

An {\it entire} function is a function analytic on all of $\CC$. We will call a function {\it almost entire} if it is analytic on a domain $W$ such that $W^c$ is a countable, closed set. This is a natural definition in relation to Brownian motion, since a Brownian motion started in $W$ will never leave $W$, just as a Brownian motion can never leave the domain of an entire function. If $f$ is entire or almost entire, then $f(B_t)$ is a time-changed Brownian motion, which means that its paths are precisely the paths of a Brownian motion, provided that we can show that nothing strange happens with the time change $\si_t = \int_{0}^{t} |f'(B_s)|^2 ds$; that is, we need to know that $\si_t$ doesn't blow up for some $t$, and also doesn't approach some finite limit as $t \lar \ff$ (which would correspond to the image Brownian motion $f(B_t)$ stopping). That is, we need the following result.

\begin{theorem} \label{entirego}
Let $f$ be almost entire and nonconstant. Then $\si_t = \int_{0}^{t} |f'(B_s)|^2 ds$ is almost surely finite for any $t$, and $\si_t \nearrow \ff$ almost surely as $t \lar \ff$.
\end{theorem}

{\bf Proof:} Since $f$ is nonconstant, $f'$ is not identically 0, which means we can choose $\eps>0$ so that the open set $A=\{z: |f'(z)|>\eps\}$ is nonempty. Then

\begin{equation} \label{}
 \int_{0}^{\ff} |f'(B_s)|^2 ds \geq \int_{0}^{\ff} \eps^2 1_A(B_s)ds = \eps^2 \int_{0}^{\ff} 1_A(B_s)ds = \ff \mbox{ a.s.}
\end{equation}
by Proposition \ref{opencamp}. This proves the second part of the statement, and the first is immediate, since $f'$ is continuous and thus bounded on the a.s. compact sets $\{B_s:0 \leq s \leq t\}$. \qed

We can use this fact to obtain some results on entire and almost entire functions. For instance, since Brownian motion is recurrent and the paths of $f(B_t)$ are the paths of a Brownian motion, $f(B_t)$ must be recurrent as well. This gives

\begin{theorem} \label{entiredense}
Let $f$ be nonconstant and entire or almost entire. Then the image of $f$ is dense in $\CC$.
\end{theorem}

This gives several immediate corollaries which are quite important.

\begin{corollary}[Liouville's Theorem] \label{}
A bounded entire or almost entire function is a constant.
\end{corollary}

\begin{corollary}[Fundamental Theorem of Algebra] \label{}
Suppose $p(z) = \bb_k z^k + \ldots + \bb_0$ is a degree $k$ polynomial, with $\bb_0, \ldots, \bb_k \in \CC$ and $\bb_k \neq 0$. Then $p$ has a root; that is, there is an $a \in \CC$ such that $f(a) = 0$.
\end{corollary}

{\bf Proof:} Since $p(\CC)$ is dense in $\CC$, we can find a sequence of points $\{z_n\}$ such that $p(z_n) \lar 0$. The compactness of $\hat \CC$ allows us to extract a convergent subsequence, and passing to this subsequence if necessary allows us to assume that $z_n \lar a \in \hat C$. We must show that $a \neq \ff$. Note that we can write

\begin{equation} \label{}
|p(z)| = |z|^k \Big|\bb_k + \frac{\bb_{k-1}}{z} + \ldots + \frac{\bb_0}{z^k}\Big|.
\end{equation}
As $z \lar \ff$, we have $|\bb_k + \frac{\bb_{k-1}}{z} + \ldots + \frac{\bb_0}{z^k}| \lar |\bb_k| \neq 0$, so that $|p(z)| \lar \ff$; thus, $a \neq \ff$. We then have $a \in \CC$, and by continuity we obtain $f(a) = \lim_{n \lar \ff} f(z_n) = 0$. \qed

In fact, $p(z)$ will have $k$ roots; see Exercise \ref{rootspoly} below.

\vski

Theorem \ref{entiredense} is just the beginning for entire functions. We will show later using Brownian motion that, if $f$ is entire, then $f(\CC)$ is all of $\CC$ with at most one point removed. This is known as Picard's Little Theorem.

\section{Exercises}

\ccases{rootspoly}  Suppose $p(z) = \bb_k z^k + \ldots + \bb_0$ is a degree $k$ polynomial, with $\bb_0, \ldots, \bb_k \in \CC$ and $\bb_k \neq 0$. Show that $p$ has $k$ roots, counting multiplicities. That is, show that there are values $a_1, \ldots, a_k, c \in \CC$ such that $p(z) = c(z-a_1)(z-a_2)\ldots (z-a_k)$. (Hint: If $p(a) = 0$, consider $p(z)-p(a)$).

\setcounter{cccases}{0}
\chapter{Simply connected domains} \lll{chscdom}

\section{Contour integration and Cauchy's Theorem} \lll{CTsec}

One of the most important theorems in complex analysis is Cauchy's Theorem, which is the subject of this section. Ordinarily, Cauchy's Theorem is formulated using contour integration, but there is an equivalent form which will be more useful for our purposes. First, we need a definition.

\vski

A domain $W$ is {\it simply connected} if $W^c$ is connected, where the complement is taken in $\hat \CC$. It is important for this definition that $\ff$ be considered a point in $W^c$, for the first domain in the picture below may appear to have a disconnected complement, but in fact it is connected in the sphere. Armed with this definition, we may state the theorem.

\begin{theorem} [Cauchy's Theorem] \label{CTdiff}
Let $f$ be an analytic function on a simply connected domain $W$. Then there is an analytic antiderivative of $f$ on $W$; that is, there is an analytic function $F$ on $W$ such that $F'(z) = f(z)$.
\end{theorem}

In order to prove this theorem, we introduce contour integration. Let $\ga(t)$ be a continuous, piecewise smooth function from an interval $[a,b] \subseteq \RR$ into $\CC$ which is only non-smooth at finitely many points in $[a,b]$; we will refer to such a function and its image as a {\it rectifiable curve}. The image $\ga ([a,b])$ will be a connected curve in $\CC$, and we will abuse notation somewhat by simply referring to this curve as $\ga$. The {\it contour integral along $\ga$ }(also known as the {\it complex line integral}), denoted $\int_\ga f dz$,
is defined by the formula

\be
\int_{\ga} f(z) dz = \int_{a}^b f(\ga(t)) \ga'(t) dt.
\ee

It can be shown easily that this quantity does not depend on the parameterization chosen (see Exercise \ref{intprop}). A curve $\ga$ is {\it closed} if $\ga(a) = \ga(b)$. The following is an equivalent form of Cauchy's Theorem.

\begin{theorem} [Cauchy's Theorem - contour integral form] \lll{CTint}
If $f(z)$ is analytic on a simply connected domain $W$, and $\ga$ is a closed path in $W$, then

\be \lll{outfield}
\int_{\ga} f(z)\, dz = 0.
\ee
\end{theorem}

The two forms of Cauchy's Theorem are equivalent by the following proposition.

\begin{proposition} \label{CTequiv}
A function $f$ on a domain $W$ has an analytic antiderivative on $W$ if, and only if, $\int_{\ga} f(z)\, dz = 0$ for any closed curve $\ga \subseteq W$.
\end{proposition}

{\bf Remark:} Note that $W$ is not assumed to be simply connected here.

\vski

{\bf Proof of Proposition \ref{CTequiv}:} Suppose that $f$ has an analytic antiderivative $F$ on $W$. Then, for any closed curve $\ga$, we have

\begin{equation} \label{}
\int_{\ga} f(z)\, dz = \int_{a}^b f(\ga(t)) \ga'(t) dt = \int_{a}^b \Big(\frac{d}{dt} F(\ga(t))\Big)dt = F(\ga(b))-F(\ga(a)) = 0.
\end{equation}
On the other hand, if \rrr{outfield} holds, then we may fix $z_0 \in W$ and, for any $z \in W$, define

\begin{equation} \label{}
F(z) = \int_{z_0}^z f(w) dw,
\end{equation}
where the notation $\int_{z_0}^z$ indicates that the integral is taken over any curve traveling from $z_0$ to $z$. This is well defined, for if $\ga_1, \ga_2$ are two such such curves, then we can form a new curve $\ga$ which follows $\ga_1$ from $z_0$ to $z$ and then follows $\ga_2$ backwards from $z_0$ to $z$; $\ga$ is a closed curve, which means that $\int_{\ga} f(z)dz =0$, and the integral over $\ga_2$ traversed backwards is equal to the integral over $\ga_2$ traversed forward (see Exercise \ref{intprop}); it follows that $\int_{\ga_1} f(z)dz = \int_{\ga_2} f(z)dz$. With this notation, it is clear from \rrr{outfield} and Exercise \ref{intprop} that

\begin{equation} \label{}
\int_{z_0}^{z_1} f(w) dw + \int_{z_1}^{z_2} f(w) dw + \int_{z_2}^{z_0} f(w) dw = 0,
\end{equation}
since a curve traveling from $z_0$ to $z_1$ to $z_2$ and back to $z_0$ is closed, and thus

\begin{equation} \label{}
\begin{split}
\lim_{h \lar 0} \frac{F(z+h)-F(z)}{h} &= \lim_{h \lar 0} \frac{\int_{z_0}^{z+h}f(w) dw - \int_{z_0}^{z}f(w) dw}{h} \\
& = \lim_{h \lar 0} \frac{1}{h}\int_{z}^{z+h}f(w) dw.
\end{split}
\end{equation}
For $h$ sufficiently small we can take this interval over the line segment connecting $z$ and $z+h$, which we can parameterize as $\ga(t) = \{z + th: 0 \leq t \leq 1\}$. We then have

\begin{equation} \label{}
\lim_{h \lar 0} \frac{F(z+h)-F(z)}{h} = \lim_{h \lar 0} \frac{1}{h} \int_{0}^{1}f(z+ht) h dt = \lim_{h \lar 0} \int_{0}^{1}f(z+ht) dt = f(z).
\end{equation}
Thus, $F$ is analytic in $W$, and $F'(z) = f(z)$. \qed

\begin{corollary} [Morera's Theorem]\label{}
If a continuous function $f$ on a domain $W$ satisfies $\int_{\ga} f(z)\, dz = 0$ for any closed curve $\ga \subseteq W$ then $f$ is analytic.
\end{corollary}

{\bf Proof:} Proposition \ref{CTequiv} shows that $f = F'$ for an analytic function $F$, and Theorem \ref{holseries} then shows that $f$ is itself analytic. \qed

We will now prove Cauchy's Theorem; in particular, we will prove the integral form, Theorem \ref{CTint}.

\vski

{\bf Proof of Theorem \ref{CTint}:} We will progress through a sequence of cases, ending with the general case.

\vski

{\bf Case 1:} $\ga$ is contained in a disk $D(z_0,r)$ which is contained in $W$. Here we expand $f$ in a power series around $z_0$, $f(z) = \sum_{n=0}^{\ff} \bb_n (z-z_0)^n$, and we are assured that the radius of convergence of the series is at least $r$. If we let $F(z) = \sum_{n=0}^{\ff} \frac{\bb_n}{n+1} (z-z_0)^{n+1}$, then the radius of convergence of the series for $F$ is the same as the series for $f$ (see Exercise \ref{radconvsame}), and we can differentiate termwise to get $F'(z) = f(z)$. By Proposition \ref{CTequiv}, $\int_{\ga} f(z) dz = 0$.

\vski

{\bf Case 2:} $\ga$ is a triangle. Exercise \ref{intprop} shows that we can assume that $\ga$ is traversed counterclockwise. Here Case 1 does not necessarily directly apply, but the idea is to express $\ga$ as the sum of a finite number of smaller triangles, each of which is contained in a disk contained in $W$; note that it will be here that we make use of simple connectivity, since it guarantees us that the interior of $\ga$ is contained in $W$. Let $\ga_1^1 = \ga$, then let $\ga_1^2, \ldots, \ga_4^2$ be the four triangles obtained by connecting the midpoints of the sides of $\ga_1^1$, and then $\ga_1^3, \ldots, \ga_{16}^3$ be the sixteen triangles obtained by connecting the midpoints of the sides of each $\ga_j^2$, and so forth.

\vspace{-1.6in}

\includegraphics[width=180mm,height=145mm]{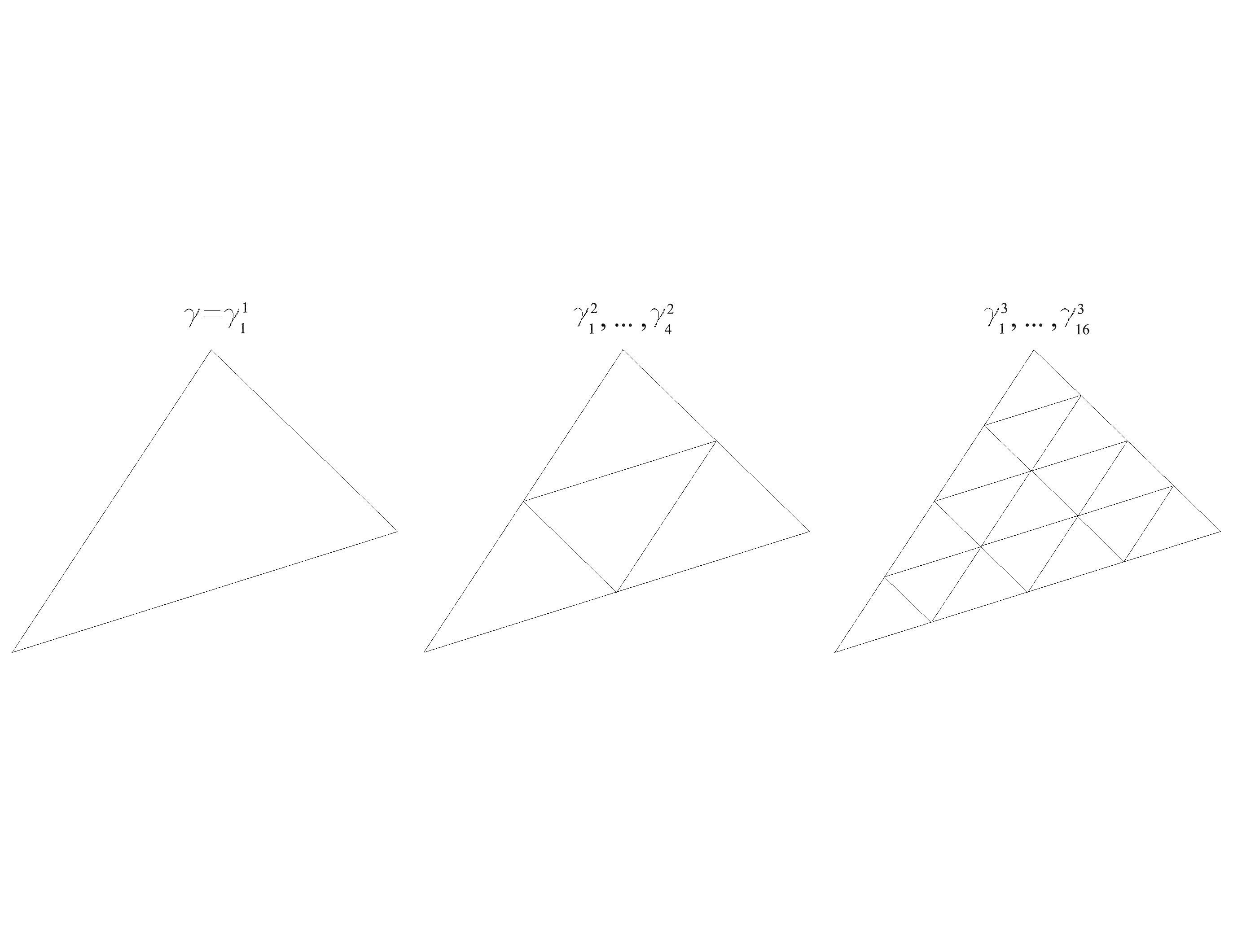}

\vspace{-1.55in}

If we sum the integrals of $f$ over each of the triangles at a given level traversed counterclockwise, we see that all but the outside line segments are integrated over in both directions. Exercise \ref{intprop} (ii) then shows that the contribution of each of these interior segments vanishes, and we obtain

\begin{equation} \label{}
\int_{\ga} f(z) dz = \sum_{j=1}^{4^{n-1}} \int_{\ga_j^n} f(z) dz.
\end{equation}
Let $d(z,\dd W) = \inf \{|z-w|: w \in \dd W\}$. By compactness (see Exercise \ref{compdist}), there is an $\eps>0$ such that $D(z,\eps) \subseteq W$ for all $z$ on $\ga$ and in the interior of $\ga$. We now need only choose $n$ large enough so that the diameter of $\ga_j^n$ is less than $\eps$, and we see that each of the triangles $\ga_j^n$ is contained in a disk contained in $W$. It follows from Case 1 that the integral over each $\ga_j^n$ is $0$, and thus $\int_{\ga} f(z) dz = 0$.

\vski

{\bf Case 3:} $\ga$ is a simple polygon; that is, $\ga$ is a piecewise linear curve with no self-intersections. Again we can assume that $\ga$ is traversed counterclockwise. It is not hard to show that any such polygon admits a triangulation $\TT$ (see Appendix Jordan Curve), and the integral over each triangle is 0 by Case 2. If we sum the integrals over the triangles traversed counterclockwise, then as in Case 2 each interior line segment will be traversed once in each direction, and by canceling these quantities we see that $\int_{\ga} f(z) dz$ is equal to the sum of the integrals over the triangles in $\TT$.  We conclude that $\int_{\ga} f(z) dz = 0$.

\vspace{.2in}

\hspace{.6in} \includegraphics[width=150mm,height=80mm]{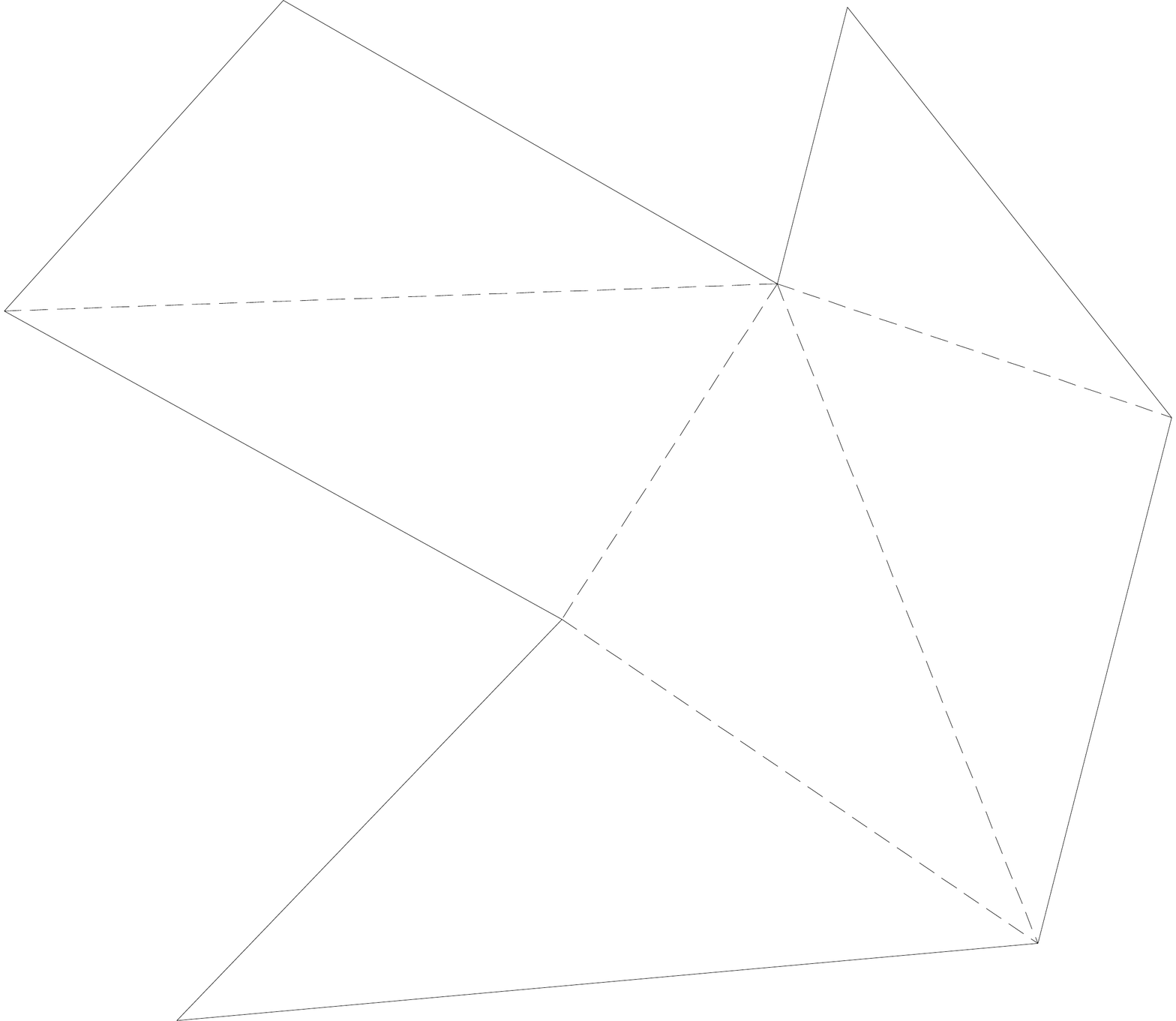}

%\vspace{-1.55in}

\vski

{\bf Case 4:} $\ga$ is a polygonal path. Here $\ga$ has a finite set of intersections, and we can break $\ga$ into the union of a finite number of simple polygons. The integral over each of these components of $\ga$ is 0, and thus $\int_{\ga} f(z) dz = 0$.

\vski

{\bf Case 5:} $\ga$ is arbitrary. Again by compactness there is an $\eps>0$ such that $D(z,\eps) \subseteq W$ for all $z$ on $\ga$, and we can choose $a=t_0 < t_1 < \ldots < t_N = b$ such that $\ga \subseteq \cup_{k=1}^N D(\ga(t_k),\eps)$. We can then choose $z_k \in D(\ga(t_k),\eps) \cap D(\ga(t_{k+1}),\eps)$ for $k = 1, \ldots ,N$. We now let $\tilde \ga$ be the polygonal path formed by taking the union of the line segments from $z_k$ to $z_{k+1}$.

\vspace{.2in}

\hspace{.6in} \includegraphics[width=140mm,height=100mm]{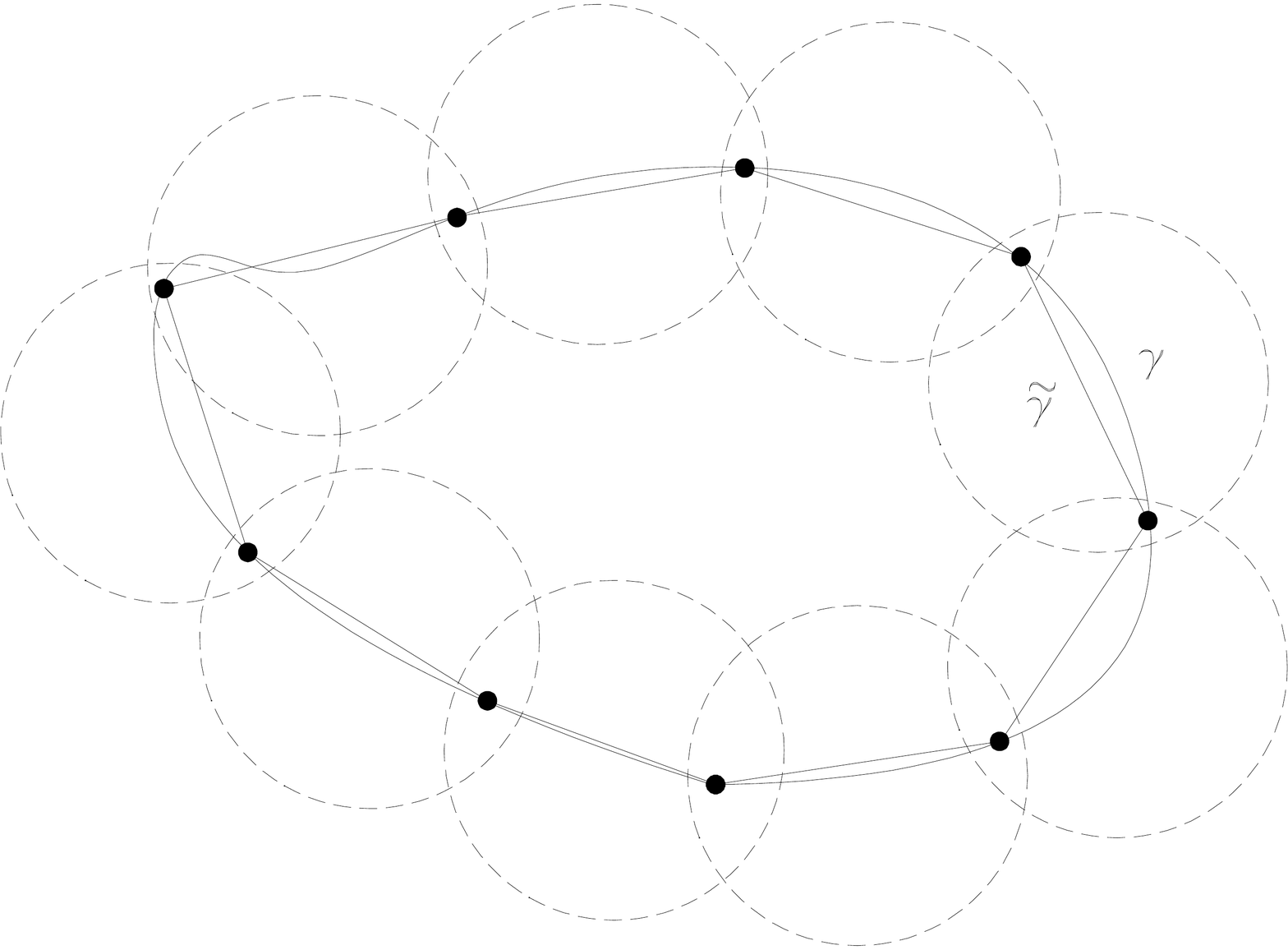}

%\vspace{-1.55in}

Since the line segment connecting $z_k$ and $z_{k+1}$ lies in the disk $D(\ga(t_{k+1}),\eps)$, the integral over this segment is the same as the integral over the portion of $\ga$ connecting $z_k$ to $z_{k+1}$ by Case 1. It follows that $\int_{\ga} f(z) dz =\int_{\tilde \ga} f(z) dz$, and since $\tilde \ga$ is a polygonal path we have $\int_{\ga} f(z) dz = 0$ by Case 4. Thus, $\int_{\ga} f(z) dz = 0$. This completes the proof of Cauchy's Theorem. \qed

We now deduce some important consequences.

\begin{corollary} \label{}
Let $f$ be an analytic function on a simply connected domain $W$, with $f(z) \neq 0$ for all $z \in W$. Then there is an analytic logarithm of $f$ on $W$; that is, there is an analytic function $g$ on $W$ such that $e^g = f$.
\end{corollary}

{\bf Proof:} If $f \neq 0$ on $W$, Cauchy's theorem tells us that we may find an antiderivative $\tilde g$ of the analytic function $\frac{f'}{f}$. We then have

\begin{equation} \label{}
\frac{d}{dz} fe^{-\tilde g} = f' e^{-\tilde g} + f e^{-\tilde g} \tilde g' = f' e^{-\tilde g} + f e^{-\tilde g} \frac{f'}{f} = 0.
\end{equation}
It follows that $f = Ce^{\tilde g}$ for some constant $C \neq 0$. Our desired logarithm is then $g = \tilde g + \tilde C$, where $\tilde C$ is chosen so that $e^{\tilde C} = C$. \qed

\begin{corollary} \label{}
Let $f$ be an analytic function on a simply connected domain $W$, with $f(z) \neq 0$ for all $z \in W$. Then there is an analytic $n$-th root of $f$ on $W$; that is, there is an analytic function $g$ on $W$ such that $g^n = f$.
\end{corollary}

{\bf Proof:} Choose $\tilde g$ so that $e^{\tilde g} = f$. The desired $n$-th root is then $g=e^{\tilde g/ n}$. \qed

\begin{corollary} \label{}
Let $h$ be a (real valued) harmonic function on a simply connected domain $W$. Then there is a harmonic conjugate of $h$ on $W$; that is, there is a harmonic function $\tilde h$ on $W$ such that $h+ i \tilde h$ is analytic.
\end{corollary}

{\bf Proof:} Let us first note that if $h$ were the real part of an analytic function $f$, then by the Cauchy-Riemann equations we could write $f'= h_x - i h_y$. We therefore begin by setting $g=h_x - i h_y$ and letting $f$ be an antiderivative of $g$. If $f=u+iv$, then the Cauchy-Riemann equations again show that $f' = u_x - i u_y$, and it follows that $h_x = u_x, h_y=u_y$. From this we conclude that there is a constant $c$ such that $h +c = u$, and we see that we can take $\tilde h = v - ic$ to obtain the conjugate to $h$. \qed

\section{Green's function}

In the field of analysis, Green's function $G(x,y)$ on regions of $\RR^n$ is formally defined to be the solution of $LG(x,y) = \dd(y-x)$, where $L$ is a linear differential operator. In complex analysis, where the Laplacian is the differential operator of most importance, for a given domain $W \subseteq \CC$ and $z \in W$ the Green's function of the Laplacian is generally defined by the following.

\begin{defn} \label{anal}

The Green's function $G_W(z,w)$ on a domain $W$ is a function in $w$ on $W \bsh \{z\}$ satisfying the following properties.

\begin{itemize} \label{analdef}

\item[(i)] $G_W(z,w)$ is harmonic and positive on $W \bsh \{z\}$.

\item[(ii)] $G_W(z,w) \lar 0$ as $w \lar bd(W)$, where $bd(W)$ is the set-theoretic boundary of $W$ (Note that the boundary is to be taken in the Riemann sphere $\hat \CC := \CC \cup \{ \ff \}$, so that if $W$ is unbounded then $\ff \in W$).

\item[(iii)] $G_W(z,w) + \frac{1}{\pi} \ln |w-z|$ extends to be continuous (and therefore harmonic) at $w=z$.

\end{itemize}

\end{defn}

Note that the factor of $\frac{1}{\pi}$ in $(iii)$ is not standard but has been chosen to align with the probabilistic considerations to follow. Not every domain has a Green's function as defined above, as for instance it can be shown that no such function can exist on the punctured disk $\DD^\times = \{0 < |z| < 1\}$, or more generally on a domain with isolated singularities. The Green's function is of tremendous importance in analysis on $\RR^n$, including complex analysis, and the question of what domains possess a Green's function has been keenly studied by analysts over the years. As will be seen below, the existence of the Green's function on simply connected domains leads to a proof of the celebrated Riemann mapping theorem.

\vski

We will now bring Brownian motion into the picture, which gives us an advantage in that we may attach Green's function to a stopping time rather than to a domain, leading to a useful generalization. Let us begin by letting $\rho_t^{T(W)}(z,w)$ be the probability density function at point $w$ and time $t$ of a Brownian motion started at $z$ and stopped at the exit time $T(W)$ of $W$. We then can calculate formally, for any suitable function $f$,

\begin{equation} \label{exploctime}
\begin{split}
E_z \int_{0}^{T(W)} f(B_s) ds & =  \int_{0}^{T(W)} E_z[f(B_s)] ds \\
& = \int_{0}^{\ff} \int_{W} f(w) \rho_s^W(z,w) dA(w) ds \\
& = \int_{W} \Big( \int_{0}^{\ff} \rho_s^W(z,w) ds \Big) f(w) dA(w).
\end{split}
\end{equation}
We see that the function

\begin{equation} \label{greeny}
G_{T(W)}(z,w) =  \int_{0}^{\ff} \rho_s^W(z,w) ds,
\end{equation}
is the density of the {\it occupation measure} $\mu(A) = \int_{0}^{T(W)} 1_A(B_s)ds$. We will refer to \rrr{greeny} as the {\it Green's function of $T(W)$,} and show that it satisfies properties $(i)$ through $(iii)$ above, but first let us note that the definition \rrr{greeny} makes sense if $T(W)$ is replaced by an arbitrary stopping time $\tau$. Naturally, we will denote the resulting function as $G_{\tau}(z,w)$. Our first order of business is showing that $G_{\tau}(z,w)$ is conformally invariant, and subsequently providing an extension of this statement to analytic functions which are not necessarily injective.

\begin{proposition} \label{greenconfinv}
Green's function is conformally invariant. That is, if $f$ is a conformal map on $W$, then $G_{T(W)}(z,w) = G_{T(f(W))}(f(z),f(w))$.
\end{proposition}

{\bf Proof:} As was described above, $G_{T(W)}(z,w)$ is the density of the occupation measure $\mu_z(A) = \int_{0}^{T(W)} 1_A(B_s)ds$, where the subscript of $z$ signifies that $B_0 = z$ a.s. If $\la$ denotes Lebesgue measure in the plane, then we have

\begin{equation} \label{}
G_{T(W)}(z,w) = \lim_{\dd \searrow 0} \frac{\mu_z(D(w,\dd))}{\la(D(w,\dd))}; \qquad G_{f(\Om)}(f(z),f(w)) = \lim_{\dd \searrow 0} \frac{\mu_{f(z)}(f(D(w,\dd)))}{\la(f(D(w,\dd)))}.
\end{equation}

Let $\eps > 0$ be given. Conformality implies that $f'(w) \neq 0$, and that for sufficiently small $\dd$ we have $D(f(w),|f'(w)|\dd(1-\eps)) \subseteq f(D(w,\dd))) \subseteq D(f(w),|f'(w)|\dd(1+\eps))$, so that $\frac{\la(f(D(w,\dd)))}{\la(D(w,\dd))} \in (|f'(w)|^2(1-\eps)^2,|f'(w)|^2(1+\eps)^2)$. On the other hand, in L\'evy's Theorem the scaling factor for time is $|f'|^2$ as well; that is, we can choose $\dd$ sufficiently small so that $\frac{|f'(w')|}{|f'(w)|} \in (1-\eps,1+\eps)$ for all $w' \in D(w,\dd)$, and then we will have $\frac{\mu_\Om(f(D(w,\dd)))}{\mu_\Om(D(w,\dd))} \in (|f'(w)|^2(1-\eps),|f'(w)|^2(1+\eps))$. Combining these two estimates and letting $\eps \lar 0$ shows that $\lim_{\dd \searrow 0} \frac{\mu_\Om(D(w,\dd))}{\la(D(w,\dd))} = \lim_{\dd \searrow 0} \frac{\mu_{f(\Om)}(f(D(w,\dd)))}{\la(f(D(w,\dd)))}$, and the result follows. \qed

It is important to note, however, that L\'evy's theorem in fact does not require maps to be injective, permitting general nonconstant analytic functions as well. This allows us to extend the previous proposition (with essentially the same proof) as follows.

\begin{theorem} \label{chimass}
Let $W$ be a domain, and suppose $f$ is a function analytic on $W$. Let $B_t$ be a Brownian motion starting at $a$, and $\tau$ a stopping time such that the set of Brownian paths $\{B_t: 0 \leq t \leq \tau\}$ lie within $W$ a.s. Let $\hat \tau = \si(\tau)$, where $\si$ is defined by L\;evy's Theorem. Then

\begin{equation} \label{cherry2}
G_{\hat \tau} (f(z),w) = \sum_{w' \in f^{-1}(\{w\})} n(f,w') G_\tau (z,w'),
\end{equation}

where $n(f,w')$ is the order of the zero of $f-w$ at $w'$.
\end{theorem}

The proof of the Proposition \ref{greenconfinv} applies to this more general statement at all points which are not images under $f$ of a point $w'$ at which $f'(w')=0$ (the {\it critical values} of $f$), with the only difference being that a point $w$ with multiple preimages will accumulate mass in the projected Green's function corresponding to mass accumulated at each of the preimages by the initial Brownian motion. The $n(f,w')$ term is not necessary in order to calculate a density of the occupation measure, as in fact $n(f,w') = 1$ except on the zero set of $f'$, which is a discrete set and therefore of Lebesque measure 0; the term is required merely to make the right side of \rrr{cherry2} continuous in $w$ when $G_\tau$ is as well (such as when $\tau$ is an exit time of a domain), since if $w$ is a critical value of $f$ then points near $w$ will have as many preimages as does $w$ only when multiplicities of preimages are counted. For an example of this phenomenon see Exercise ???.

\vski

Our eventual goal is to prove the Riemann mapping theorem in the next section, and we must first calculate several of the simplest Green's functions. We note first the two-dimensional Gaussian density $\rho_t (z,w) = \frac{1}{2\pi t} e^{-|z-w|^2/(2t)}$; it is evident that $\int_{0}^{\ff}\rho_t (z,w)dt = \ff$ for any $z,w$, and thus when $\tau = \ff$ we have $G_{\tau}(z,w) = \ff$.

\vski

As usual, let $\HH = \{Im(z)>0\}$. The density $\rho_t^{T(\HH)} (z,w)$ must capture the probability of Brownian paths near $w$ at time $t$, but only those which have not previously intersected $\RR$. By the reflection principle (ref ???), the processes $B_t$ and

\be \hat B_t = \left \{ \begin{array}{ll}
B_t & \qquad  \mbox{if } t \leq T(\HH)  \\
\overline B_t & \qquad \mbox{if } t > T(\HH) \;,
\end{array} \right. \ee

have the same law; but $B_t$ is near $w$ and $t > T(\HH)$ precisely when $\hat B_t$ is near $\bar w$ with $t > T(\HH)$, and this occurs precisely when $\hat B_t$ is near $\bar w$, since the Brownian motion cannot travel from $z$ to $\bar w$ without first crossing $\RR$. We conclude that $\rho_t^{T(\HH)} (z,w) = \rho_t (z,w) - \rho_t (z,\bar w) = \frac{1}{2\pi t} (e^{-|z-w|^2/(2t)} - e^{-|z-\bar w|^2/(2t)})$. We therefore have

\begin{equation} \label{hp}
\begin{split}
G_{T(\HH)}(z,w) & = \int_{0}^{\ff} \frac{1}{2\pi t} (e^{\frac{-|z-w|^2}{2t}} - e^{\frac{-|z-\bar w|^2}{2t}}) dt \\
& = \frac{1}{2\pi} \int_{0}^{\ff} \frac{1}{t} \int_{|z-w|}^{|z-\bar w|} \frac{a}{t}e^{\frac{-a^2}{2t}} da dt \\
& = \frac{1}{2\pi} \int_{|z-w|}^{|z-\bar w|} \int_{0}^{\ff} \frac{a}{t^2} e^{\frac{-a^2}{2t}} dt da \\
& = \frac{1}{\pi} \int_{|z-w|}^{|z-\bar w|} \frac{1}{a} \int_{0}^{\ff} e^{-u} du da \\
& = \frac{1}{\pi} \ln \Big( \frac{|z-\bar w|}{|z-w|}\Big).
\end{split}
\end{equation}

Now, $\phi(z) = -i \Big(\frac{z-1}{z+1}\Big)$ maps $\DD$ conformally onto $\HH$, sending $0$ to $i$. Thus,

\begin{equation} \label{}
\begin{split}
G_{T(\DD)}(0,z) & =G_{T(\HH)}(i,\phi(z)) = \frac{1}{\pi} \ln \Big( \frac{|i-i(\frac{\bar z -1}{\bar z + 1})|}{|i+i(\frac{z -1}{z + 1})|}\Big) \\
& = \frac{1}{\pi} \ln \Big| \frac{2z+2}{2|z|^2 + 2z}\Big| = \frac{1}{\pi} \ln \frac{1}{|z|}.
\end{split}
\end{equation}

Let us now calculate the Green's function of a stopping time defined by the winding of Brownian motion. If our Brownian motion $B_t$ starts at 1 then it will never hit 0 a.s., and we can then define $arg(B_t)$ as a continuous process (with $arg(B_0) = 0$). For any posi-

\begin{wrapfigure}{r}{8cm}
\vspace{-.7in}
\includegraphics[width=8cm]{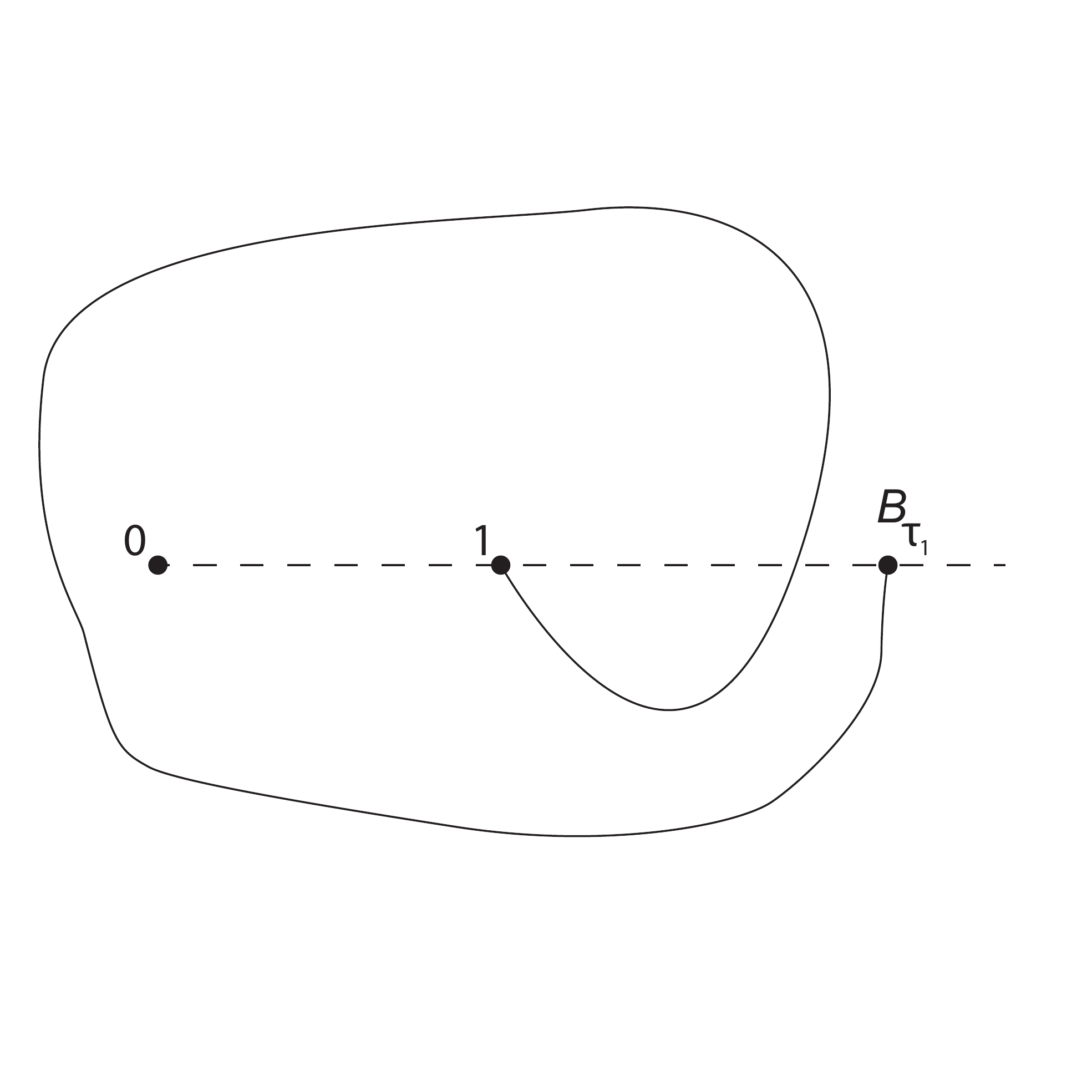}
\vspace{-.7in}
\caption{The stopping time $\tau_1$}
\end{wrapfigure}

tive integer $n$ we then define $\tau_n = \inf\{t: arg(B_t) = \pm 2 \pi n\}$. In other words, $\tau_n$ is the first time at which the Brownian motion has wound around the origin in either direction $n$ times. The image to the right gives an intuitive illustration of a Brownian path up until time $\tau_1$. If we let $U = \{Re(z)>0\}$ denote the right half-plane, we see that the function $z \lar z^{4n}$ transforms a Brownian motion starting at $1$ and stopped at $T(U)$ into one starting at 1 and stopped at $\tau_n$. We may therefore calculate $G_{\tau_n}$ by using Proposition \ref{chimass}. For any point $re^{i\th}$, with $r >0$ and $\th \in (-\pi, \pi]$, we have

\begin{equation} \label{emily}
\nonumber G_{\tau_n}(1,re^{i\th}) = \sum_{\substack{k \in \ZZ \\ (\frac{\th}{4n}+ \frac{2\pi k}{4n}) \in [-\frac{\pi}{2},\frac{\pi}{2}]}} G_{T(U)}(1,r^{\frac{1}{4n}} e^{i(\frac{\th}{4n}+ \frac{2\pi k}{4n})}).
\end{equation}

The explicit formula is now easy to produce if desired, using the fact that reflection across the $y$-axis is given by $z \lar -\bar z$, so that \rrr{hp} gives $G_{T(U)}(z,w) = \frac{1}{\pi} \ln \frac{|z+\bar w|}{|z-w|}$; however, our purpose in examining this stopping time is the proof of the Riemann mapping theorem in the next section, and for that only a few basic properties of $G_{\tau_n}$ are required. To begin with, since $G_{\tau_n}$ is expressed as the sum of a finite number of finite terms, it is finite for all $re^{i\th} \neq 1$. Furthermore, $G_{\tau_n}(1,re^{i\th}) + \frac{1}{\pi} \ln |1-re^{i\th}|$ can be extended continuously at $re^{i\th} = 1$, since the only term in the sum which has a singularity corresponds to $k=0$ and has singular part $-\frac{1}{\pi} \ln |1-r^{\frac{1}{4n}} e^{i\frac{\th}{4n}}|$, and we have

\begin{equation} \label{}
%\begin{split}
\lim_{re^{i\th} \lar 1} \frac{1}{\pi} \ln |1-re^{i\th}|-\frac{1}{\pi} \ln |1-r^{\frac{1}{4n}} e^{i\frac{\th}{4n}}| = \lim_{\zeta \lar 1} \ln \frac{1-\zeta^{4n}}{1-\zeta} = \ln 4n.
%\end{split}
\end{equation}

Finally, if $r \lar 0$ then all of the preimages of $re^{i\th}$ approach 0 as well, and we conclude that $G_{\tau_n}(1,re^{i\th}) \lar 0$; an analogous argument shows that $G_{\tau_n}(1,re^{i\th}) \lar 0$ whenever $r \lar \ff$.
%For certain domains $G_\Om$ is infinite everywhere, and for others it is not. This function is of fundamental importance, and satisfies the following properties.

%\begin{theorem} \label{}
%\begin{itemize} \label{}

%\item[(i)] If $G_{\Om}(z,w)< \ff$ for some $z,w \in \Om$, then it is finite for all pairs $z,w \in G$ with $z \neq w$. In this case we will say Green's function {\it exists} on $\Om$.

%\item[(ii)] If $G_\Om$ exists, then as a function of $w$ with fixed $z$, $G_{\Om}(z,w)$ is harmonic on $\Om$, approaches $0$ as $w$ approaches $\dd \Om$, and $G_{\Om}(z,w) - \frac{1}{\pi} \ln |w-z|$ is harmonic on $\Om$.

%\item[(iii)] If $G_\Om$ exists, then $G_{\Om}(z,w) = G_{\Om}(w,z)$ for all $z,w \in \Om$.

%\item[(iv)] If $f: \Om_1 \lar \Om_2$ is conformal, then, for any $z,w \in \Om_1$ we have $G_{\Om_1}(z,w) = G_{\Om_2}(f(z),f(w))$.

%\end{itemize}

%\end{theorem}

%\begin{proposition} \label{}
%\be
%G_\DD(0,w) = \frac{1}{\pi} \ln |w|
%\ee
%\end{proposition}

%{\bf Proof:} The given function clearly satisfies $(ii)$, so we need only show that it is the only such function. So suppose $h(w)$ satisfies the conditions in $(ii)$. Then $h(w) - \ln |w|$ is harmonic on $\DD$ and equal to 0 on $\dd \DD$. However, we know that $h(w) - \ln |w|$ can be expanded into a power series on $\DD$ with coefficients given by \rrr{sercoeffform}. These coefficients must all be zero, which implies that $h(w) - \ln |w| = 0$ on $\DD$, and thus $h(w) = \ln |w|$. \qed

\section{Conformal mapping and the Reimann mapping theorem}

The Riemann mapping theorem is as follows.

\vski

{\bf Riemann mapping theorem} If $W \subset \CC$ is a simply connected domain with $W \neq \CC$, and $a \in \Om$, then there is a conformal map $f$ from $W$ onto $\DD$, with $f(a)=0$.

\vski

Note that if such a map $f$ existed, then the function $\frac{-1}{\pi}\ln | f(z) |$ would be satisfy the conditions $(i)-(iii)$ given for the Green's function in the previous section, and thus would be equal to $G_W(a,z)$. We have already seen that the Green's function exists for $W$, and therefore the trick is to use $G_U$ to construct $f$. %For this we will use a clever argument due to Greg Lawler (personal communication)...

\section{Schwarz-Christoffel transformations}

Another important class of conformal maps are the {\it Schwarz-Christoffel transformations}. In the simplest form, this maps the upper half-plane $\HH = \{Im(z) > 0 \}$ conformally onto a polygon. It is easiest to understand this type of map by looking first at the map $f(z) = z^m$, for $0<m\leq 2$, which takes $\HH$ conformally to the domain $\{0< Arg(z) < m \pi\}$.  We see that this function turns an angle of $\pi$ (a straight angle) into an angle of $m \pi$.  Note also that the derivative, up to a constant, is $z^{m-1}$.

\vski

Suppose now that we want to find a conformal map $f$ from $\HH$ onto a polygon with vertices $w_1, \ldots, w_n$ with corresponding angles (to the left) of $\aaa_1\pi, \ldots, \aaa_n\pi$.  Such a conformal map will extend to a bijection between the boundaries of the domains, and there will therefore be points $z_1, \ldots , z_{n-1}$ on $\RR$ which are mapped to $w_1, \ldots, w_{n-1}$, and we will arrange it so that $\ff$ is mapped to $w_n$.  Between any two of these points, $z_j$ and $z_{j+1}$ say, the argument of $f'(z)$ must remain constant, but then it must have a jump of $\aaa_j \pi$ at $z_j$.  This suggests that

$$
f'(z) = C (z-z_1)^{\aaa_1-1}(z-z_2)^{\aaa_2-1} \ldots (z-z_{n-1})^{\aaa_2-1}.
$$

This turns out to be correct, and we have the following theorem.

\begin{theorem} \label{}
The conformal map $f$ from $\HH$ onto a polygon with vertices $w_1, \ldots, w_n$ with corresponding angles (to the left) of $\aaa_1\pi, \ldots, \aaa_n\pi$  is given by

$$
f(z) =A+ C \int^{z}(\zeta-z_1)^{\aaa_1-1}(\zeta-z_2)^{\aaa_2-1} \ldots (\zeta-z_{n-1})^{\aaa_2-1} d\zeta,
$$

for some complex constants $A,C$, where $w_j = f(z_j)$ for $j=1, \ldots, n-1$.
\end{theorem}

Most explicit conformal maps that are known can be expressed as compositions of M\"obius transformations and Schwarz-Christoffel transformations. On the other hand, the class of conformal maps is in fact very large (as is shown by the Riemann Mapping Theorem below), we just don't have formulas for most of them.

\section{Exercises}

\ccases{idententire} a) Find all entire functions $f(z)$ such that $f(1/n) = 1/n^2$ for all $n \in \NN$.

b) Find all entire functions $f(z)$ such that $f(1/n) = 1/\sqrt{n}$ for all $n \in \NN$.

\vski

\ccases{openaxes} Suppose $f(z)$ is an analytic function on a domain $W$ with the property that $f(z)^n \in \RR$ for all $z \in W$. Show that $f$ is a constant function.

\vski

\ccases{lindprep} Suppose that $\ga$ is a {\it Jordan curve}: that is, $\ga$ can be realized as the image of a continuous and one-to-one function $\ga(t)$ from $[0,1]$ into $\CC$. $\ga$ then determines a bounded domain $W$ in its interior. Suppose that $f(z)$ is an analytic function on $W$, and suppose that there are points $w_1, \ldots , w_n \in \ga$ such that $\limsup_{\{z \lar \dd W \backslash \{w_1, \ldots ,w_n\}\}} |f(z)| \leq m < \ff$. Suppose further that there is some finite $M > m$ such that $|f(z)| \leq M$ for all $z \in W$. Show that in fact $|f(z)| \leq m$ for all $z \in W$.

\hspace{1.2in} \includegraphics[width=90mm,height=70mm]{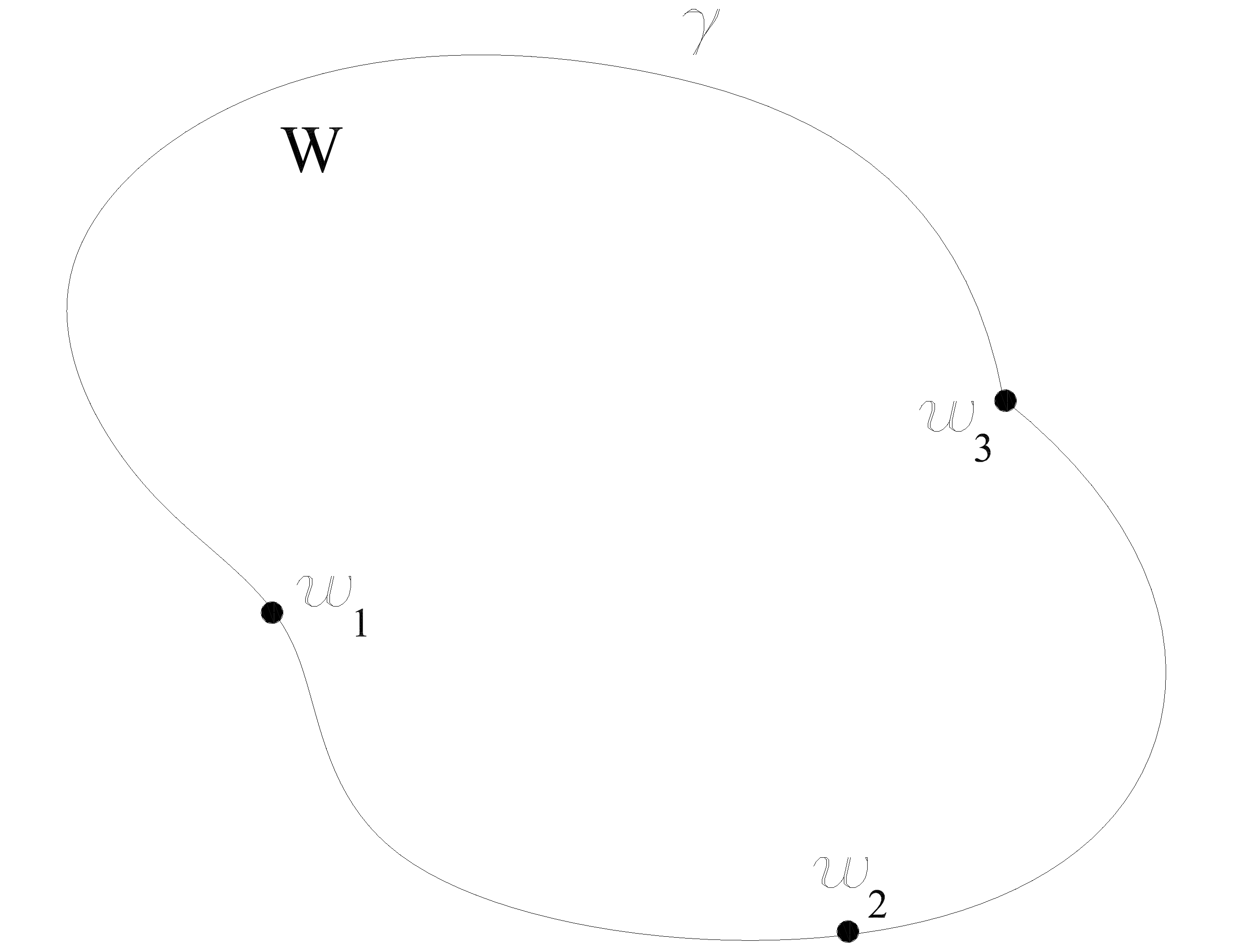}

\vski

\ccases{rotatethis} Let $W$ be the square domain $\{-1<x<1,-1<y<1\}$. Suppose $f(z)$ is a bounded analytic function on $W$ such that $f(z) \lar 0$ as $z$ approaches the bottom edge of the square; that is, $\lim_{\{z \lar \{y=-1\}\}} |f(z)| = 0$. Prove that $f(z) = 0$ for all $z \in W$.

\vski

\ccases{k-roots} Suppose $p(z) = \bb_k z^k + \ldots + \bb_0$ is a degree $k$ polynomial, with $\bb_0, \ldots, \bb_k \in \CC$ and $\bb_k \neq 0$. Show that $p$ has $k$ roots, counting multiplicities. That is, show that there are values $a_1, \ldots, a_k, c \in \CC$ such that $p(z) = c(z-a_1)(z-a_2)\ldots (z-a_k)$. (Hint: If $p(a) = 0$, consider $p(z)-p(a)$).

\vski

\ccases{intprop} i) (Independence of parameterizations) Suppose $\ga_1: [a_1,b_1] \lar \CC$ and $\ga_2: [a_2,b_2] \lar \CC$ are rectifiable curves, and suppose that there is a smooth increasing function $\phi : [a_1,b_1] \lar [a_2,b_2]$ with $\phi(a_1) = a_2, \phi(b_1) = b_2$ such that $\ga_1(t) = \ga_2(\phi(t))$ for all $t \in [a_1,b_1]$. Show that

\begin{equation} \label{}
\int_{\ga_1} f(z)dz = \int_{\ga_2} f(z)dz
\end{equation}
for any analytic function $f$ defined in a neighborhood of $\ga_1=\ga_2$.

\vski

ii) (Path reversal) Let $\ga: [a_1,b_1] \lar \CC$ be a rectifiable curve. Define $\tilde \ga : [-b_1,-a_1] \lar \CC$ by $\tilde \ga (t) = \ga (-t)$. Then $\tilde \ga$ is simply the curve $\ga$ traversed in the opposite direction. Show that

\begin{equation} \label{}
\int_{\tilde \ga} f(z)dz = - \int_{\ga} f(z)dz
\end{equation}
for any analytic function $f$ defined in a neighborhood of $\ga$.

\vski

iii) (Concatenating paths) Suppose $\ga_1: [a,b] \lar \CC$ and $\ga_2: [c,d] \lar \CC$ are rectifiable curves such that $\ga_1(b) = \ga_2(c)$. Define a new curve $\tilde \ga: [a, b+d-c] \lar \CC$ by

\begin{equation} \label{17}
\tilde \ga(t) := \begin{cases} \ga_t(t) & t \in [a,b],\\
\ga_2(t-b+c) & t \in [b,b+d-c].
\end{cases}
\end{equation}
$\tilde \ga$ is the concatenation of $\ga_1$ and $\ga_2$. Show that

\begin{equation} \label{}
\int_{\tilde \ga} f(z)dz = \int_{\ga_1} f(z)dz + \int_{\ga_2} f(z)dz
\end{equation}
for any analytic function $f$ defined in a neighborhood of $\tilde \ga$.

\vski

\ccases{compdist} Let $d(z,\dd W)$ be as in Section \ref{CTsec}. Show that the function $z \lar d(z,\dd W)$ is continuous and positive on $W$. Conclude that if $K$ is a compact set in $W$ then there is some $\eps > 0$ such that $d(z,\dd W) > 0$ for all $z \in K$.

\vski

\ccases{radconvsame} Show that the radius of convergence is the same for the three series $\sum_{n=0}^{\ff} \bb_n (z-z_0)^n, \sum_{n=0}^{\ff} \frac{\bb_n}{n+1} (z-z_0)^{n+1},$ and $\sum_{n=1}^{\ff} n \bb_n (z-z_0)^{n-1}$.

\vski

\ccases{polytriang} Polygon admits triangulation.

{\bf \Large 1.} Let $\ga = \{e^{it}: 0 \leq t \leq 2\pi\}$. Evaluate the following two integrals.

\begin{equation} \label{}
\int_{\ga} \frac{1}{z}dz, \qquad \int_{\ga} \frac{1}{z^2}dz.
\end{equation}

{\bf \Large 2.} Let $\ga(t): [0, 5\pi]$ be defined as follows.

\begin{equation} \label{17}
\ga(t) = \begin{cases} (t+\frac{3\pi}{2})e^{it} & 0 \leq t \leq 2 \pi,\\
\frac{3\pi}{2} + 2\pi e^{-it} & 2\pi \leq t \leq 3 \pi,\\
(t-\frac{7 \pi}{2}) + i \cos (t-\frac{7 \pi}{2}) & 3 \pi \leq t \leq 5\pi.
\end{cases}
\end{equation}

Sketch $\ga$, and evaluate

\begin{equation} \label{}
\int_{\ga} e^{1/z} dz .
\end{equation}

{\bf \Large 3.}  Suppose $p(z) = \bb_k z^k + \ldots + \bb_0$ is a degree $k$ polynomial, with $\bb_0, \ldots, \bb_k \in \CC$ and $\bb_k \neq 0$. Show that $p$ has $k$ roots, counting multiplicities. That is, show that there are values $a_1, \ldots, a_k, c \in \CC$ such that $p(z) = c(z-a_1)(z-a_2)\ldots (z-a_k)$. (Hint: If $p(a) = 0$, consider $p(z)-p(a)$).

\ccases{exlineremove} Suppose $W$ is a domain and $C$ is a line intersecting $W$. Suppose further that $f$ is a continuous function on $W$ which is analytic on $W \bsh C$. Show that in fact $f$ is analytic on all of $W$. (Hint: Morera's Theorem)

\setcounter{cccases}{0}
\chapter{The reflection principle and consequences} \label{chrefl}

In this chapter we study a number of important consequences of the reflection principle for Brownian motion. Let $T_\RR$ be the hitting time of $\RR$, and set

\begin{equation} \label{}
\hat B_t = \begin{cases} B_t & t \leq T_\RR,\\
\seg{10}{B_t} & t > T_\RR.
\end{cases}
\end{equation}

The reflection principle for Brownian motion says that $\hat B_t$ is itself a Brownian motion.

\section{The reflection principle for harmonic functions} \label{secreflharm}

The reflection principle for Brownian motion leads immediately to the reflection principles for harmonic functions.

\begin{theorem} \label{reflharmthm}
Suppose that $h^+(z)$ is harmonic on a domain $W^+ \subseteq \HH$. Suppose further that $h^+(z) \lar 0$ as $Im(z) \lar 0$ in $W^+$. Let $W^- = \{ \bar z: z \in W^+ \}$ be the reflection of $W^+$ over $\RR$, and let $W = W^+ \cup W^- \cup (\dd W^+ \cap \RR)^o$, where $(\dd W^+ \cap \RR)^o$ is the interior of the set $(\dd W^+ \cap \RR)$, taken as a subset of $\RR$. For $z \in W$, define

\begin{equation} \label{17}
h(z) = \begin{cases} h^+(z) & z \in W^+,\\
-h^+(\bar z) & z \in W^-,\\
0 & z \in (\dd W^+ \cap \RR)^o.
\end{cases}
\end{equation}
Then $h(z)$ is harmonic on $W$.
\end{theorem}
{\bf Remark:} Note that it is not being claimed that $W$ is a domain, since it will not be connected if $(\dd W^+ \cap \RR)^o = \emptyset$. However, $W$ is a domain if $(\dd W^+ \cap \RR)^o \neq \emptyset$, and an open set in either case.

\vski

{\bf Proof:} Harmonicity is a local property, so it suffices to consider $z$ in the various regions of $W$ separately. Clearly, $h$ is harmonic on $W^+$. If $z_0 \in W^-$, we know from Theorem \ref{harmseries} that $h^+$ can be expanded into a power series in a neighborhood of $\seg{8}{z_0}$, so that we may write

\begin{equation} \label{}
\begin{split}
h(z) = -h^+(\bar z) & = -\Big( \bb_0 + \sum_{n=1}^{\ff} \bb_{-n} \seg{30}{(\bar z-\seg{8}{z_0})}^n + \sum_{n=1}^{\ff} \bb_{n} (\bar z-\seg{8}{z_0})^n \Big)\\
& = -\bb_0 + \sum_{n=1}^{\ff} (-\bb_{-n}) (z-z_0)^n + \sum_{n=1}^{\ff} (-\bb_{n}) \seg{30}{(z-z_0)}^n.
\end{split}
\end{equation}
This shows that $h$ is harmonic in a neighborhood of $z_0$. Finally, if $z_0 \in (\dd W^+ \cap \RR)^o$, then we can choose $r>0$ such that $\bar D(z_0,r) \subseteq W$. Let $D = D(z_0,r)$. If we can show that $E_z[h(B_{T_{D}})] = h(z)$ for any $z \in D$, then by Proposition \ref{itsharm} it will follow that $h$ is harmonic on $D$, and in particular at $z_0$. If $z \in (D \cap \RR)$, then by the symmetry of $h$ it is clear that $E_z[h(B_{T_{D}})] = 0 = h(z)$. Suppose now that $z \in (W^+ \cap D)$. Let $D^+ = (W^+ \cap D)$, and write

\begin{equation} \label{}
E_z[h(B_{T_{D}})] = E_z[h(B_{T_{D^+}})1_{\{T_{D^+} < T_{D}\}}] + E_z[h(B_{T_{D^+}})1_{\{T_{D^+} \geq T_{D}\}}].
\end{equation}
Note that, since $\hat B_t$ is itself a Brownian motion, we must have $E_z[h(B_{T_{D^+}})1_{\{T_{D^+} \geq T_{D}\}}] = E_z[h(\hat B_{T_{D^+}})1_{\{T_{D^+} \geq T_{D}\}}]$. On the other hand, the definition of $h$ shows that

\be
E_z[h(\hat B_{T_{D^+}})1_{\{T_{D^+} \geq T_{D}\}}] = E_z[h(\seg{18}{B_{T_{D^+}}})1_{\{T_{D^+} \geq T_{D}\}}] = -E_z[h(B_{T_{D^+}})1_{\{T_{D^+} \geq T_{D}\}}].
\ee
We conclude that $E_z[h(B_{T_{D^+}})1_{\{T_{D^+} \geq T_{D}\}}] = 0$. Note also that $E_z[h(B_{T_{D^+}})1_{\{B_{T_{D^+}} \in \RR\}}] = 0$, since $h = 0$ on $\RR$. We therefore obtain

\begin{equation} \label{}
\begin{split}
E_z[h(B_{T_{D}})] & = E_z[h(B_{T_{D^+}})1_{\{T_{D^+} < T_{D}\}}] = E_z[h(B_{T_{D^+}})1_{\{T_{D^+} < T_{D}\}}] + E_z[h(B_{T_{D^+}})1_{\{B_{T_{D^+}} \in \RR\}}] \\
& = E_z[h(B_{T_{D^+}})] = h(z),
\end{split}
\end{equation}
with the last equality due to the harmonicity of $h$ on $W^+$ and continuity to $\RR$. \qed

A similar result holds for analytic functions.

\begin{theorem} \label{reflanalthm}
Suppose that $f^+(z)$ is analytic on a domain $W^+ \subseteq \HH$. Suppose further that $f^+(z)$ extends continuously to the boundary component $(\dd W^+ \cap \RR)^o$, defined as in Theorem \ref{reflharmthm}. Let $W$ also be defined as in Theorem \ref{reflharmthm}. For $z \in W$, define

\begin{equation} \label{17}
f(z) = \begin{cases} f^+(z) & z \in W^+,\\
\seg{25}{f^+(\bar z)} & z \in W^-,\\
\lim_{w \in W^+, w \lar z}f^+(w) & z \in (\dd W^+ \cap \RR)^o.
\end{cases}
\end{equation}
Then $f(z)$ is analytic on $W$.
\end{theorem}

{\bf Proof:} The proof is very similar to that of Theorem \ref{reflharmthm}. It is clear that $f$ is analytic on $W^+$. For $z_0 \in W^-$, we can write $f(z) = f^+(z) = \sum_{n=0}^{\ff} \bb_n(z-\seg{10}{z_0})^n$ in a neighborhood of $\seg{10}{z_0}$, so that in a neighborhood of $z_0$ we have

\begin{equation} \label{}
f(z) = \seg{25}{f^+(\bar z)} = \sum_{n=0}^{\ff} \seg{50}{\bb_n(\bar z-\seg{10}{z_0})^n} = \sum_{n=0}^{\ff} \seg{10}{\bb_n}(z-z_0)^n.
\end{equation}
$f$ is therefore analytic in $W^-$. Analyticity on $(\dd W^+ \cap \RR)^o$ follows from the same argument as was used for Theorem \ref{reflharmthm}; alternatively, the result in Exercise \ref{exlineremove} of Chapter \ref{chscdom} may be applied. \qed

Both reflection principles as stated are rather restrictive: it may seem that they may be applied only in very specific instances. However, let us recall that the real axis is merely a circle on the sphere, and since any circle can be taken to any other by a M\"obius transformation, the analytic invariance of Brownian motion suggests that these principles hold in a much more general setting. To formulate the more general result, we need to understand what it means to reflect a point over a circle. If $C$ is a circle on the sphere, we will find a M\"obius transformation $\phi$ taking $C$ to $\RR$, then for any $z$ define $z^* = \phi^{-1}(\seg{25}{\phi(z)})$. $z^*$ will be referred to as the {\it reflection of $z$ over $C$}. We remark that there will be many M\"obius transformations taking $C$ to $\RR$; however, it can be shown (see Exercise \ref{exreflform} below) that the value of $z^*$ is independent of the choice of $\phi$. Armed with this definition, we can state the more general reflection principle for analytic functions.

\begin{theorem} \label{genrefthm}
Let $C_1, C_2$ be circles in the sphere, and let $*_1, *_2$ denote the reflection over these two circles, respectively. Suppose $W^+$ is a domain lying in one component of $C_1^c$, and suppose further that $f^+$ is analytic in $W^+$ and extends continuously to $\dd C_1$ with $f^+(C_1) \subseteq C_2$. Let $W^- = \{ z^{*_1}: z \in W^+\}$, and set $W = W^+ \cup W^- \cup (\dd W^+ \cap C_1)^o$, where as before $(\dd W^+ \cap C_1)^o$ denotes the interior of $(\dd W^+ \cap C_1)$ taken as a subset of $C_1$. Define

\begin{equation} \label{17}
f(z) = \begin{cases} f^+(z) & z \in W^+,\\
f^+(z^{*_1})^{*_2} & z \in W^-,\\
\lim_{w \in W^+, w \lar z}f^+(w) & z \in (\dd W^+ \cap C_1)^o,
\end{cases}
\end{equation}

Then $f(z)$ is analytic on $W$.
\end{theorem}

{\bf Proof:} Let $\phi_1, \phi_2$ be M\"obius transformations taking $C_1, C_2$ to $\RR$. Then apply Theorem \ref{reflanalthm} to $\phi_2 \circ f \circ \phi_1^{-1}$. \qed

\section{Picard's Theorem} \label{secpicard}

The reflection principle allows us to construct a remarkable analytic function from $\HH$ onto $\CC \bsh \{-1,1\}$ known as the {\it modular function}. Let $T=\{|z|>1, y>0, -1<x<-1\}$. $T$ is a curvilinear triangle on the sphere; that is, a region with three sides which are each segments of circles on the sphere. Let $f$ be the conformal map from $T$ onto $\HH$ which extends continuously to a map from $\bar T$ to $\bar \HH$ fixing $-1, 1,$ and $\ff$; the existence of such a map $g$ without the final condition is guaranteed by the Riemann Mapping Theorem and Carath\'eodory's Theorem, and we may then let $f=\phi \circ g$, where $\phi$ is the M\"obius transformation taking $g(-1), g(1),$ and $g(\ff)$ to $-1, 1,$ and $\ff$, respectively. The boundary components $\{x=-1, y>0\}, \{|z|=1, y>0\},$ and $\{x=1, y>0\}$ are then mapped under $f$ to $\{x<-1, y = 0\}, \{-1<x<1, y = 0\},$ and $\{x>1, y = 0\}$, respectively. These sets are all segments of circles on the sphere, so the reflection principle allows us to reflect $f$ over each boundary component of $T$, extending $f$ to a function defined on the sets $T_1, T_2,$ and $T_3$ as illustrated in the picture below. The image of each of these sets under $f$ will be the reflection of $\HH$ over $\RR$, which is the lower half-plane $\{y<0\}$. We then reflect again over the boundary arcs of $T_1, T_2, T_3$, extending $f$ into more triangular domains whose image is now the upper half-plane again. Continuing in this manner, we are abe to tessellate $\HH$ with triangles, thereby extending $f$ to a map from $\HH$ onto $\CC \bsh \{-1,1\}$, as desired. In the following picture, the three boundary components are colored red, blue, and green, with the colors corresponding to the boundary components in the image $\HH$. The dotted colored lines correspond to the reflections in the domain of $f$, which are still mapped to the corresponding colors in the image.

\vspace{-1.5in}

\includegraphics[width=180mm,height=145mm]{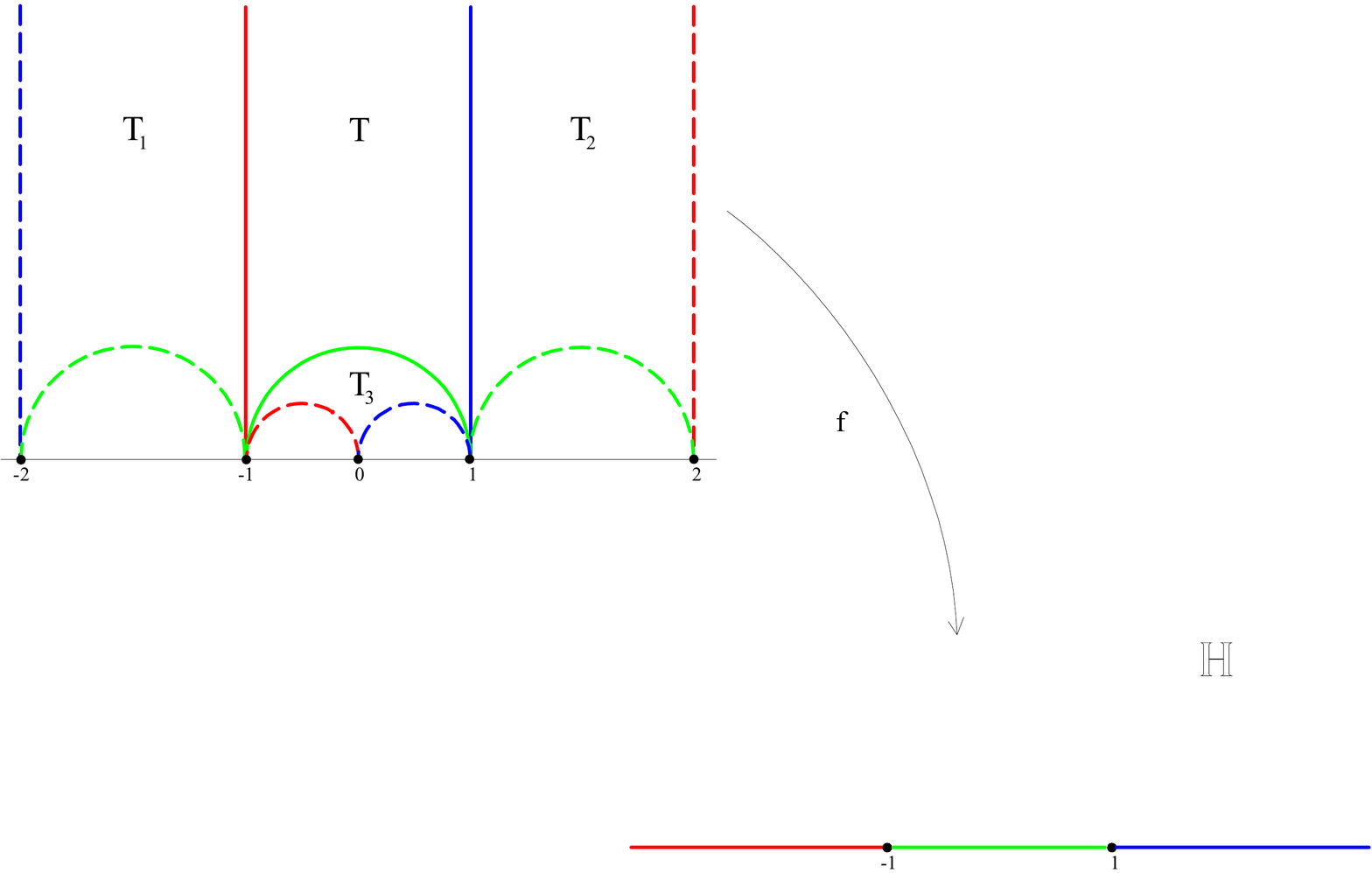}

%\vspace{-1.55in}

Every point $w$ in $\CC \bsh \{-1,1\}$ will have many preimages under $f$ in $\HH$, and if $\ga$ is a curve in $\HH$ connecting any two such preimages $z_1, z_2$, then $f(\ga)$ will be a closed curve beginning and ending at $w$. However, if we take $z_1 \neq z_2$ then the curve $f(\ga)$ will in fact enclose one or both of $\{-1,1\}$. For instance, in the picture below, $f(\ga)$ first crosses $\{x<-1, y = 0\}$ from below and then $\{x>1, y = 0\}$ from above before returning to $w$ and encircling both of $\{-1,1\}$.

\vspace{-1.5in}

\includegraphics[width=180mm,height=145mm]{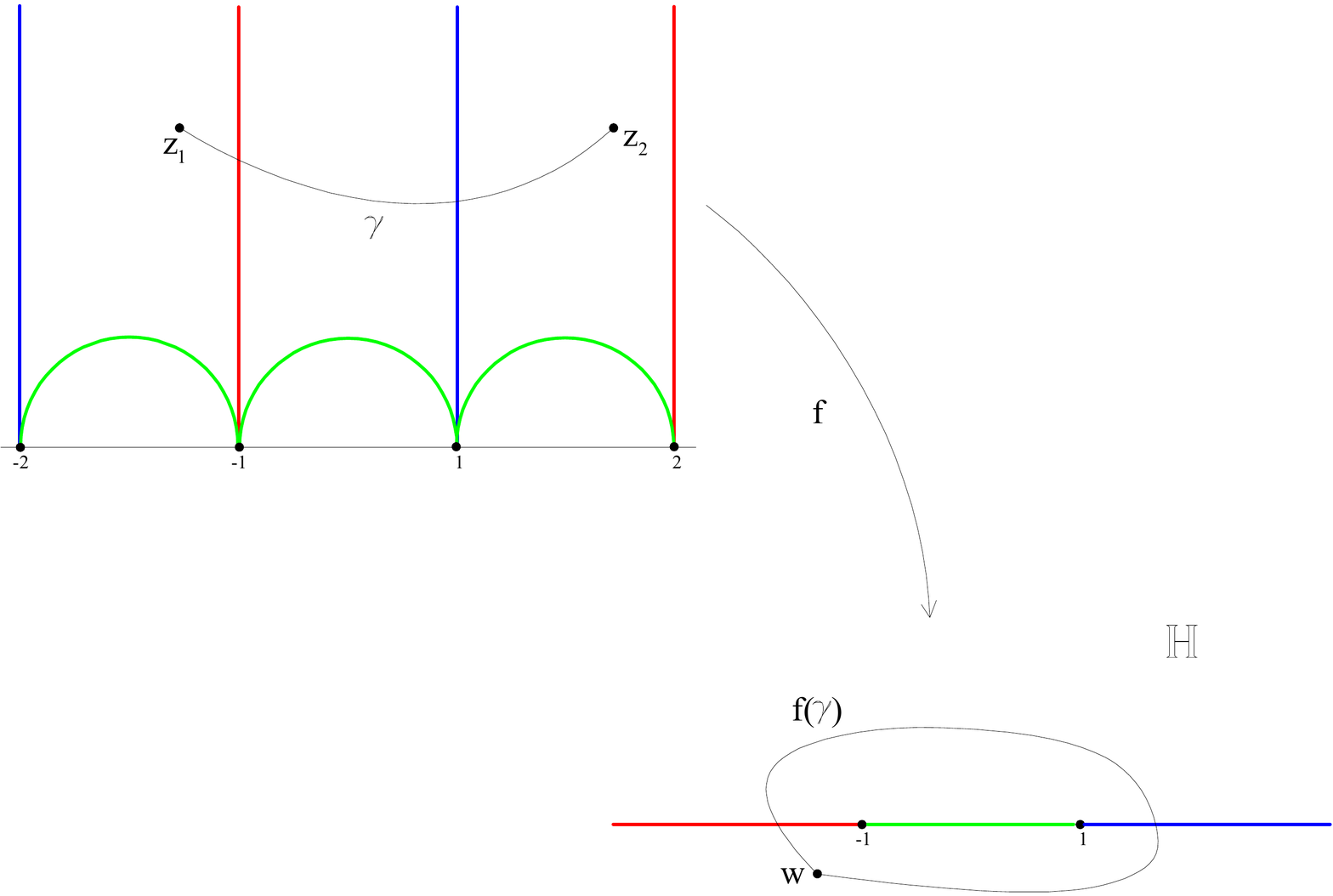}

The reader is invited to sketch more examples to get the idea. Essentially, if $f(\ga)$ is a closed curve starting and finishing at $w$, then $\ga$ connects two preimages of $w$, and the particular preimages which $\ga$ connects determines the exact manner in which $f(\ga)$ winds its way through around $-1$ and $1$. We need a way to make this more precise. \\

Suppose $W$ is a domain and $\ga_1, \ga_2: [0,1] \lar W$ are two curves such that $\ga_1(0) = \ga_2(0), \ga_1(1) = \ga_2(1)$. We will say $\ga_1$ {\it is homotopic to $\ga_2$ in $W$} if there is a a continuous map $\tilde \ga: [0,1] \times [0,1] \lar W$ with the following properties:

\begin{itemize} \label{}

\item[(i)] $\tilde \ga(0,t) = \ga_1(t)$ for $t \in (0,1)$.

\item[(ii)] $\tilde \ga(1,t) = \ga_1(t)$ for $t \in (0,1)$.

\item[(iii)] $\tilde \ga(s,0) = \ga_1(0) = \ga_2(0)$ and $\tilde \ga(s,1) = \ga_1(1) = \ga_2(1)$ for $s \in (0,1)$.

\end{itemize}

The map $\tilde \ga$ is called a {\it homotopy}, and as $s$ goes from $0$ to $1$ the curve $\tilde \ga(s,t)$ represents a continuous deformation of the curve $\ga_1(t)$ into $\ga_2(t)$. In case of the existence of such a $\tilde \ga$, we will write $\ga_1 \sim_{{}_W} \ga_2$, or simply $\ga_1 \sim \ga_2$ when the domain $W$ is understood. Note that in deforming $\ga_1$ into $\ga_2$, $\tilde \ga$ is not allowed to leave $W$. This is the part of the definition which will allow us to make rigorous a statement like "$f(\ga)$ encircles $\{-1,1\}$". A constant function $\ga_2(t) = w \in W$ will be referred to as a {\it point}, and therefore the phrase {\it $\ga_1$ is homotopic to a point} means simply that $\ga_1$ can be continuously contracted to a trivial curve. \\

Let us return to the modular function. If $\ga$ is a curve connecting two preimages $z_1,z_2$ of $w$ under $f$, then as discussed before $f(\ga)$ is a closed curve beginning and ending at $w$, and it can be shown that $f(\ga)$ is homotopic to a point if, and only if, $z_1 = z_2$. With $z_1$ fixed, the many possible choices of $z_2$ therefore correspond to different homotopy classes of closed curves beginning and ending at $w$. Let us see how this relates to Brownian motion. If we start a Brownian motion $B_t$ at $a \in \HH$ and stop it upon hitting $\RR$, then $f(B_t)$ is a time-changed Brownian motion in $\CC \bsh \{-1,1\}$ stopped upon exiting $\CC \bsh \{-1,1\}$. However, we know that with probability 1 a Brownian motion will never hit $\{-1,1\}$, so the process $f(B_t)$ in fact traces out a full Brownian path in $\CC$ as $t \nearrow T_\HH$. Let $\dd>0$ be chosen so that $D=D(f(a),\dd) \subseteq \CC \bsh \{-1,1\}$. We know from the recurrence of Brownian motion that $f(B_t)$ lies with $D$ for values of $t$ arbitrarily close to $T_\HH$, which implies that $B_t$ intersects components of the set $f^{-1}(D)$ for values of $t$ arbitrarily close to $T_\HH$. However, as $t \nearrow T_\HH$, Brownian motion leaves the component $O_a$ of $f^{-1}(D)$ containing $a$. Let us therefore suppose we have $t$ such that $B_t \in f^{-1}(D) \bsh O_a$. We choose $\tilde a$ in the component of $f^{-1}(D)$ containing $B_t$, and let $\tilde \ga$ be a curve lying in $f^{-1}(D)$ connecting $B_t$ and $\tilde a$ and $\ga = \{B_s: 0 \leq s \leq t\} + \tilde \ga$. $\ga$ is a curve connecting two distinct preimages of $f(a)$ in $\HH$, so $f(\ga)$ is a closed curve beginning and ending at $f(a)$ which is not homotopic to a point in $\CC \bsh \{-1,1\}$. We have proved an important fact about planar Brownian motion.

\begin{theorem} \label{}
Let $B_t$ be a Brownian motion starting at $a \in \CC \bsh \{-1,1\}$, and let $D$ be a small disk such that $a \in D \subseteq \CC \bsh \{-1,1\}$. Then there is a stopping time $\tau< \ff$ a.s. such that if $t > \tau$ and $B_t \in D$, and $\ga = \{B_s: 0 \leq s \leq t\} + \tilde \ga$ where $\tilde \ga \subseteq D$ connects $B_t$ to $a$ then $\ga$ is not homotopic to a point.
\end{theorem}

That is, $B_t$ eventually gets tangled about $\{-1,1\}$ and never becomes untangled. Let us contrast this with the case in which only one point is removed from the plane.

\begin{theorem} \label{}
Let $B_t$ be a Brownian motion starting at $a \in \CC \bsh \{0\}$, and let $D$ be a small disk such that $a \in D \subseteq \CC \bsh \{0\}$. Then a.s. there are arbitrarily large values of $t$ such that $B_t \in D$ and $\ga = \{B_s: 0 \leq s \leq t\} + \breve \ga$ ,where $\breve \ga \subseteq D$ connects $B_t$ to $a$, is homotopic to a point.
\end{theorem}

{\bf Proof:} $B_t$ can be realized as a time change of $e^{\hat B_t}$, where $\hat B_t$ is a Brownian motion starting at some primage $\hat a$ of $a$ under the exponential map. $\hat B_t$ returns to the component of $f^{-1}(D)$ containing $\hat a$ for arbitrarily large values of $t$, and for such values we can connect $\hat B_t$ to $\hat a$ by a curve $\hat \ga$ so that $\{\hat B_s: 0 \leq s \leq t\} + \tilde \ga$ is homotopic to a point. This homotopy, call it $\tilde \ga$, can be projected, so that in fact $e^{\tilde \ga}$ gives a homotopy of $\ga$ to a point in $\CC \bsh \{0\}$. \qed

Thus, Brownian motion gets tangled around two points, but doesn't around only one point. This property has a very interesting corollary known as {\it Picard's Theorem}.

\begin{corollary} \label{}
Suppose $f(z)$ is an entire function or an analytic function defined on $\CC \bsh \{b\}$, and suppose there are two points $w_1.w_2 \in \CC$ such that the image of $f$ is contained in $\CC \bsh \{w_1,w_2\}$. Then $f$ is a constant.
\end{corollary}

{\bf Proof:} Suppose $f$ is a counterexample to this statement. If $B_t$ is a Brownian motion starting at a point in the domain of $f$, then we know from the previous results that we can find a $t$ such that $\ga \{\hat B_s: 0 \leq s \leq t\}$ is homotopic to a point but $f(\ga)$ is not. This is a contradiction, since we can project a homotopy through $f$. \qed

\section{Laurent series} \label{seclaurent}

Let $D^\times(w, \dd)$ refer to $D(w,\dd)$ with $w$ removed; that is, $D^\times(w, \dd) = \{0< |z-w| < \dd\}$. If a function $f$ is analytic on a domain $W$ which does not contain $w$ but which does contain $D^\times(w, \dd)$ for some $\dd > 0$, then we will say that $f$ has a {\it singularity at $w$}. We will examine the behavior of $f(z)$ near $w$ closely. Let us distinguish three cases.

\begin{itemize} \label{}

\item[Case 1] $\lim_{z \lar w} f(z)$ exists and equals some $v \in \CC$, and if we define $f(w) = v$ then $f$ so extended is analytic on $W \cup \{w\}$. This is known as a {\it removable singularity}.

\item[Case 2] $\lim_{z \lar w} |f(z)| = \ff$. This is called a {\it pole}.

\item[Case 3] Neither $(i)$ nor $(ii)$ occurs. This is called an {\it essential singularity}.

\end{itemize}

It is a tautological and trivial statement that one of these three must occur, but we will see in fact that functions behave in a very specific way in the neighborhood of an essential singularity. Before discussing this in more detail, let us discuss the terminology a bit. A removable singularity is so named because it is a singularity which is not really a singularity, since it can be removed. A pole is so named because if one were to graph the function $|f(z)|$ in $\RR^3$, the graph would have a vertical asymptote at $z=w$, resembling a pole of infinite height. It will be seen, furthermore, that if $f$ has a pole at $w$ then $w$ is a removable singularity for $\frac{1}{f}$; in that sense, a pole is not really a true singularity, since $f$ can be thought of as analytic at $w$ if we define $f(w) = \ff$. The true singularities, then, are the essential ones, which justifies their name. \\

The following important theorem goes a long way towards clarifying the behavior of analytic functions in a neighborhood of the various types of singularities. It is known as the Casorati-Weierstrass Theorem.

\begin{theorem} \label{casweier}
Let $f(z)$ be analytic on $D^\times(w,\dd)$ for some $\dd>0$. Then the following hold.

\begin{itemize} \label{}

\item[(i)] If $f$ is bounded on $D^\times(w,\dd)$ then $f$ has a removable singularity at $w$.

\item[(ii)] If $f(D^\times(w,\dd))$ is not dense in $\CC$, then $f$ has either a pole or removable singularity at $w$.

\end{itemize}
\end{theorem}

{\bf Proof:} Suppose first that $f$ is bounded on $D^\times(w,\dd)$. Let $B_t$ be a Brownian motion, let $a \in D^\times(w,\dd)$, and let $0<\dd'< \dd$. Since $f$ is analytic on the bounded set $D^\times(w,\dd')$ and continuous on $\CCC(D^\times(w,\dd'))$, we may apply It\=o's Theorem, which tells us that $f(a) = E_a[T_{D^\times(w,\dd')}]$. However, $B_t$ does not see points, so that $P_a[B_{T_{D^\times(w,\dd')}} = w] \leq P_a[B_t = w \mbox{ for some } t \geq 0] = 0$, and we conclude that $f(a) = E_a[T_{D(w,\dd')}]$. We may now follow the proof of Theorem \ref{harmseries} to conclude that $f$ may be expanded into a power series

\begin{equation} \label{}
f(a) = \bb_0 + \sum_{n=1}^{\ff} \bb_{-n} \seg{30}{(a-w)}^n + \sum_{n=1}^{\ff} \bb_{n} (a-w)^n.
\end{equation}

However, $f$ is analytic on $D^\times(w, \dd')$, so $\frac{df}{d\bar z} = \sum_{n=1}^{\ff} n \bb_{-n} \seg{30}{(a-w)}^{n-1}$ vanishes on this domain, which implies that $\bb_{-n} = 0$ for $n \geq 1$, and therefore $f(a) = \sum_{n=0}^{\ff} \bb_{n} (a-w)^n$. It is now clear that if we define $f(w) = \bb_0$, then $f$ so defined is analytic on $D(w, \dd')$. This proves $(i)$. \\

Now suppose $f(D^\times(w,\dd))$ is not dense in $\CC$. We can then find $v \in f(D^\times(w,\dd))^c$ such that $D(v, \eps) \subseteq f(D^\times(w,\dd))^c$ for some $\eps > 0$. Then $g(z) =\frac{1}{f(z) -v}$ is analytic on $D^\times(w,\dd)$ and bounded in magnitude by $\frac{1}{\eps}$, so by $(i)$ $g$ extends to be analytic on $D(w,\dd)$. If $g(w) \neq 0$, then since $f(z) = \frac{1}{g(z)} + v$ it is clear that $f$ extends to be analytic on $D(w,\dd)$, and $f(w) = \frac{1}{g(w)} + v$. If, however, $g(z) = 0$, then if we let $m$ be the order of the zero of $g$ at $w$ we can write $g(z) = (z-w)^m h(z)$, where $h$ is analytic and nonzero in a neighborhood of $w$. But then $|f(z)| = \Big| \frac{1}{g(z)}+v \Big| = \Big| \frac{1}{(z-w)^mh(z)}+v \Big| \lar \ff$ as $z \lar w$. This proves (ii). \qed

\begin{corollary} \label{}
Suppose $f$ is analytic on $D^\times(w, \dd)$ and has a pole at $w$. Then there is a unique integer $m>0$ such that $g(z) = (z-w)^m f(z)$ can be extended to be analytic on $D(w, \dd)$ with $g(w) \neq 0$. $f(z) = \frac{g(z)}{(z-w)^m}$ can therefore be expanded in a power series

\begin{equation} \label{patra}
f(z) = \sum_{n=-m}^{+\ff} \bb_n (z-w)^n,
\end{equation}

which converges on $D(w, \dd)$. Conversely, if $f$ is given by \rrr{patra} with $\bb_{-m} \neq 0$, and $\sum_{n=0}^{+\ff} \bb_n (z-w)^n$ converges in some neighborhood of $w$, then $|f(z)| \lar \ff$ as $z \lar w$, and $f$ therefore has a pole at $w$.
\end{corollary}

{\bf Proof:} The first statement is clear from the proof of Theorem \ref{casweier}, while the second holds by expanding $g$ into a power series $\sum_{n=0}^{\ff} \tilde \bb_n (z-w)^n$ on $D(w, \dd)$, then setting $\bb_n = \tilde \bb_{n+m}$ and relabeling the indices of the series. If we suppose now that $f$ is given by \rrr{patra}, then the analytic part $\sum_{n=0}^{\ff} \bb_n (z-w)^n$ of $f$ approaches $\bb_0$ as $z \lar w$, and we have

\begin{equation} \label{}
\begin{split}
\Big| \sum_{n=-m}^{-1} \bb_n (z-w)^n \Big| & \geq \Big|\frac{\bb_{-m}}{(z-w)^m}\Big| - \Big|\frac{\bb_{-m+1}}{(z-w)^{m-1}}\Big| - \ldots - \Big|\frac{\bb_{-1}}{(z-w)}\Big| \\ & = \Big| \frac{1}{(z-w)^m} \Big| |\bb_{-m} - \bb_{-m+1}(z-w) - \ldots - \bb_{-1}(z-w)^{m-1}| \lar \ff
\end{split}
\end{equation}

as $z \lar w$. Thus, $w$ is a pole of $f$. \qed

{\bf Remark:} The integer $m$ is referred to as the {\it order} of the pole at $w$, and is the same as the order of the zero of $\frac{1}{f}$ at $w$. The series \rrr{patra} is what is referred to as a {\it Laurent series}: a power series with negative powered terms. \\

We have now obtained a clear characterization of the behavior of an analytic function in the neighborhood of a pole. The singularities for which we do not have a good understanding is then the essential variety. It turns out, though we will not prove it, that a function with an essential singularity can be expanded into a series such as \rrr{patra} in a neighborhood of singularity, but with infinitely many nonzero negative powered terms. It also turns out that the Casorati-Weierstrauss Theorem can be strengthened to the following result, known as the {\it Great} or {\it Big Picard's Theorem} (what we have referred to earlier as Picard's Theorem is then demoted to the {\it Little Picard's Theorem}).

\begin{theorem} \label{}
Suppose $f$ is analytic and has an essential singularity at $w$. Then there is a point $b \in \CC$ such that $\CC \bsh \{b\} \subseteq f(U)$ for every neighborhood $U$ of $w$.
\end{theorem}

That is, $f$ attains every value in $\CC$ with at most one exception on every neighborhood of $w$. A proof exists of this fact which uses Brownian motion, but we will not present it here. 

\section{Exercises}

\ccases{idententire} a) Find all entire functions $f(z)$ such that $f(1/n) = 1/n^2$ for all $n \in \NN$.

b) Find all entire functions $f(z)$ such that $f(1/n) = 1/\sqrt{n}$ for all $n \in \NN$.

\vski

\ccases{openaxes} Suppose $f(z)$ is an analytic function on a domain $W$ with the property that $f(z)^n \in \RR$ for all $z \in W$. Show that $f$ is a constant function.

\vski

\ccases{lindprep} Suppose that $\ga$ is a {\it Jordan curve}: that is, $\ga$ can be realized as the image of a continuous and one-to-one function $\ga(t)$ from $[0,1]$ into $\CC$. $\ga$ then determines a bounded domain $W$ in its interior. Suppose that $f(z)$ is an analytic function on $W$, and suppose that there are points $w_1, \ldots , w_n \in \ga$ such that $\limsup_{\{z \lar \dd W \backslash \{w_1, \ldots ,w_n\}\}} |f(z)| \leq m < \ff$. Suppose further that there is some finite $M > m$ such that $|f(z)| \leq M$ for all $z \in W$. Show that in fact $|f(z)| \leq m$ for all $z \in W$.

\hspace{1.2in} \includegraphics[width=90mm,height=70mm]{Lindelof.pdf}

\vski

\ccases{rotatethis} Let $W$ be the square domain $\{-1<x<1,-1<y<1\}$. Suppose $f(z)$ is a bounded analytic function on $W$ such that $f(z) \lar 0$ as $z$ approaches the bottom edge of the square; that is, $\lim_{\{z \lar \{y=-1\}\}} |f(z)| = 0$. Prove that $f(z) = 0$ for all $z \in W$.

\vski

\ccases{k-roots} Suppose $p(z) = \bb_k z^k + \ldots + \bb_0$ is a degree $k$ polynomial, with $\bb_0, \ldots, \bb_k \in \CC$ and $\bb_k \neq 0$. Show that $p$ has $k$ roots, counting multiplicities. That is, show that there are values $a_1, \ldots, a_k, c \in \CC$ such that $p(z) = c(z-a_1)(z-a_2)\ldots (z-a_k)$. (Hint: If $p(a) = 0$, consider $p(z)-p(a)$).

\vski

\ccases{intprop} i) (Independence of parameterizations) Suppose $\ga_1: [a_1,b_1] \lar \CC$ and $\ga_2: [a_2,b_2] \lar \CC$ are rectifiable curves, and suppose that there is a smooth increasing function $\phi : [a_1,b_1] \lar [a_2,b_2]$ with $\phi(a_1) = a_2, \phi(b_1) = b_2$ such that $\ga_1(t) = \ga_2(\phi(t))$ for all $t \in [a_1,b_1]$. Show that

\begin{equation} \label{}
\int_{\ga_1} f(z)dz = \int_{\ga_2} f(z)dz
\end{equation}
for any analytic function $f$ defined in a neighborhood of $\ga_1=\ga_2$.

\vski

ii) (Path reversal) Let $\ga: [a_1,b_1] \lar \CC$ be a rectifiable curve. Define $\tilde \ga : [-b_1,-a_1] \lar \CC$ by $\tilde \ga (t) = \ga (-t)$. Then $\tilde \ga$ is simply the curve $\ga$ traversed in the opposite direction. Show that

\begin{equation} \label{}
\int_{\tilde \ga} f(z)dz = - \int_{\ga} f(z)dz
\end{equation}
for any analytic function $f$ defined in a neighborhood of $\ga$.

\vski

iii) (Concatenating paths) Suppose $\ga_1: [a,b] \lar \CC$ and $\ga_2: [c,d] \lar \CC$ are rectifiable curves such that $\ga_1(b) = \ga_2(c)$. Define a new curve $\tilde \ga: [a, b+d-c] \lar \CC$ by

\begin{equation} \label{17}
\tilde \ga(t) := \begin{cases} \ga_t(t) & t \in [a,b],\\
\ga_2(t-b+c) & t \in [b,b+d-c].
\end{cases}
\end{equation}
$\tilde \ga$ is the concatenation of $\ga_1$ and $\ga_2$. Show that

\begin{equation} \label{}
\int_{\tilde \ga} f(z)dz = \int_{\ga_1} f(z)dz + \int_{\ga_2} f(z)dz
\end{equation}
for any analytic function $f$ defined in a neighborhood of $\tilde \ga$.

\vski

\ccases{compdist} Let $d(z,\dd W)$ be as in Section \ref{CTsec}. Show that the function $z \lar d(z,\dd W)$ is continuous and positive on $W$. Conclude that if $K$ is a compact set in $W$ then there is some $\eps > 0$ such that $d(z,\dd W) > 0$ for all $z \in K$.

\vski

\ccases{radconvsame} Show that the radius of convergence is the same for the three series $\sum_{n=0}^{\ff} \bb_n (z-z_0)^n, \sum_{n=0}^{\ff} \frac{\bb_n}{n+1} (z-z_0)^{n+1},$ and $\sum_{n=1}^{\ff} n \bb_n (z-z_0)^{n-1}$.

\vski

\ccases{polytriang} Polygon admits triangulation.

{\bf \Large 1.} Let $\ga = \{e^{it}: 0 \leq t \leq 2\pi\}$. Evaluate the following two integrals.

\begin{equation} \label{}
\int_{\ga} \frac{1}{z}dz, \qquad \int_{\ga} \frac{1}{z^2}dz.
\end{equation}

{\bf \Large 2.} Let $\ga(t): [0, 5\pi]$ be defined as follows.

\begin{equation} \label{17}
\ga(t) = \begin{cases} (t+\frac{3\pi}{2})e^{it} & 0 \leq t \leq 2 \pi,\\
\frac{3\pi}{2} + 2\pi e^{-it} & 2\pi \leq t \leq 3 \pi,\\
(t-\frac{7 \pi}{2}) + i \cos (t-\frac{7 \pi}{2}) & 3 \pi \leq t \leq 5\pi.
\end{cases}
\end{equation}

Sketch $\ga$, and evaluate

\begin{equation} \label{}
\int_{\ga} e^{1/z} dz .
\end{equation}

{\bf \Large 3.}  Suppose $p(z) = \bb_k z^k + \ldots + \bb_0$ is a degree $k$ polynomial, with $\bb_0, \ldots, \bb_k \in \CC$ and $\bb_k \neq 0$. Show that $p$ has $k$ roots, counting multiplicities. That is, show that there are values $a_1, \ldots, a_k, c \in \CC$ such that $p(z) = c(z-a_1)(z-a_2)\ldots (z-a_k)$. (Hint: If $p(a) = 0$, consider $p(z)-p(a)$).

\ccases{exlineremove} Suppose $W$ is a domain and $C$ is a line intersecting $W$. Suppose further that $f$ is a continuous function on $W$ which is analytic on $W \bsh C$. Show that in fact $f$ is analytic on all of $W$. (Hint: Morera's Theorem)

\setcounter{cccases}{0}
\chapter{Boundary behavior of conformal maps} \label{altrings}

\section{The Three Lines Lemma and integral inequalities} \label{sec3ll}

\begin{lemma} \label{lemsegexit}
Let $a,b \in \RR$ with $a<b$, and let $W = \{a < Re(z) < b\}$. Let $C_a = \{x=a\}$ and $C_b = \{x=b\}$ be the boundary lines of $W$. Then, if $z =x+yi \in W$, we have

\begin{equation} \label{}
P_z(B_{T_W} \in C_a) = \frac{b-x}{b-a}, \qquad P_z(B_{T_W} \in C_b) = \frac{x-a}{b-a} .
\end{equation}
\end{lemma}

{\bf Proof:}  The function $h(x+yi) = \frac{b-x}{b-a}$ is harmonic on $W$, 1 on $C_a$, and 0 on $C_b$. Apply It\=o's formula to obtain

\begin{equation} \label{}
\frac{b-x}{b-a} = h(z) = E_z[h(B_{T_W})] = P_z(B_{T_W} \in C_a).
\end{equation}

The other identity follows similarly. \qed

\begin{proposition} \label{3lineslemma}
Suppose $f(z)$ is analytic and bounded on $W=\{0<Re(z)<1\}$ and continuous on $\seg{10}{W}$. Let $C_0 = \{x=0\}$ and $C_1 = \{x=1\}$ be the boundary lines of $W$, and suppose that $\sup_{z \in C_0}|f(z)| = m$, $\sup_{z \in C_1}|f(z)| = M$. Then, for any $z = x+yi \in W$, $|f(z)| \leq m^{1-x}M^x$.
\end{proposition}

{\bf Proof:} We may assume $m \leq M$, for the symmetry of $W$ allows us to merely reflect the argument if $M < m$. Let $\ZZ = \{z:|f(z)| \leq m\}$. If $z \in \ZZ$ the result holds trivially, since $|f(z)| \leq m \leq m^{1-x}M^x$, so let us assume $z \in W \bsh \ZZ$. Note that $\ln |f(z)|$ is harmonic on $W \bsh \ZZ$, and bounded above and below, so that we may start a Brownian motion at $z$ and apply It\=o's Theorem together with Lemma \ref{lemsegexit} to obtain

\begin{equation} \label{}
\begin{split}
\ln |f(z)| & = E_z[\ln |f(B_{T_{W \bsh \ZZ}})|] \\
& \leq P_z(B_{T_{W \bsh \ZZ}} \in C_1 ) \ln M + P_z(B_{T_{W \bsh \ZZ}} \in (C_0 \cup \ZZ) ) \ln m \\
& \leq P_z(B_{T_{W}} \in C_1 ) \ln M + P_z(B_{T_{W}} \in C_0) \ln m \\
& = x \ln M + (1-x) \ln m.
\end{split}
\end{equation}
Exponentiating this inequality completes the proof. \qed

\begin{theorem} \label{}
If $f,g$ are two measurable, complex-valued functions on a measure space with measure $\mu$, and $p,q > 1$ with $\frac{1}{p} + \frac{1}{q} = 1$ then

\begin{equation} \label{feyou}
\int_{A} |fg| d\mu \leq \Big(\int |f|^p d\mu\Big)^{1/p} \Big(\int |g|^q d\mu\Big)^{1/q}.
\end{equation}
\end{theorem}

{\bf Proof:} Let $F = |f|^p, G = |g|^q$. By standard approximation techniques, it is enough to prove the result when $F,G$ are simple functions; that is, we may assume that $F = \sum_{n=1}^{N} a_n 1_{A_n}$ and $G = \sum_{m=1}^{M} b_m 1_{B_m}$, where the $a_n$'s and $b_m$'s are positive constants and the $A_n$'s and $B_m$'s are disjoint measurable sets with $\mu(A_n), \mu(B_m) < \ff$. Let

\begin{equation} \label{}
\phi(z) = \int F^z G^{1-z} d \mu = \sum_{n=1}^{N} \sum_{m=1}^{M} a_n^z b_m^{1-z} \mu(A_n \cap B_m);
\end{equation}

note that a $a_n^z$ is defined to be $e^{(\ln a_n)z}$, and similarly for $b_m^{1-z}$. It is clear that $\phi$ is entire, and since $|a_n^z| = e^{(\ln a_n) Re(z)}= a_n^{Re(z)}$, $|b_n^{1-z}| = e^{(\ln b_n) (1-Re(z))} = b_n^{1-Re(z)}$ we see that $\phi$ is bounded on the $W$ defined in Proposition \ref{3lineslemma}. Note that on $C_0$, we have

\begin{equation} \label{}
|\phi(z)| \leq  \sum_{n=1}^{N} \sum_{m=1}^{M} |a_n^z| |b_m^{1-z}| \mu(A_n \cap B_m) = \sum_{n=1}^{N} \sum_{m=1}^{M} b_m  \mu(A_n \cap B_m) \leq \sum_{m=1}^{M} b_m  \mu(B_m) = \int G d\mu.
\end{equation}
A similar calculation shows that, on $C_1$, we have

\begin{equation} \label{}
|\phi(z)| \leq \int F d\mu.
\end{equation}
Proposition \ref{3lineslemma} now gives

\begin{equation} \label{}
\int F^{1/p} G^{1/q} d \mu = \phi (1/p) \leq \Big(\int F d\mu\Big)^{1/p} \Big(\int G d\mu\Big)^{1/q}.
\end{equation}
This is equivalent to \rrr{feyou}. \qed

The special case $p=q=2$ is sufficiently useful as to have its own name, the {\it Cauchy-Schwarz inequality}. We isolate it for future reference.

\begin{theorem} \label{}
If $f,g$ are two measurable, complex-valued functions on a measure space with measure $\mu$,  then

\begin{equation} \label{feyou}
\Big(\int_{A} |fg| d\mu\Big)^2 \leq \Big(\int |f|^2 d\mu\Big) \Big(\int |g|^2 d\mu\Big).
\end{equation}
\end{theorem}

\section{Carath\'eodory's Theorem for Jordan domains} \label{secara}

A {\it Jordan curve} is the homeomorphic image of a circle in $\CC$. That is, if $\ga_t$ is a continuous and injective map from $[0,1]$ to $\CC$ such that $\ga(0) = \ga(1)$, then its image, which we will refer to as $\ga$, is a Jordan curve. A domain whose boundary is a Jordan curve is called a {\it Jordan domain}. The {\it Jordan Curve Theorem} (see Appendix) tells us that any Jordan curve $\ga$ defines a unique, bounded, simply connected domain $W$ which has $\ga$ as its boundary. The Riemann Mapping Theorem tells us that $W$ is conformally equivalent to the unit disk. If $f$ is the conformal map from $\DD$ to $W$, then $f$ is, among other things, a homeomorphism between the two domains, and since we know that $\ga$ is the homeomorphic image of a circle a very natural question is whether $f$ can be extended to a homeomorphism from $\seg{10}{\DD}$ to $\bar W$. The affirmative answer to this question is known as Carath\'eodory's Theorem.

\begin{theorem} \label{}
Let $f$ be a conformal map from $\DD$ onto a Jordan domain $W$. Then $f$ extends continuously to a map sending $\dd \DD$ bijectively to $\dd W$, and $f$ so extended is therefore a homeomorphism from $\seg{10}{\DD}$ to $\seg{10}{W}$.
\end{theorem}

{\bf Proof:} Fix $z_0 \in \dd \DD$. For $0< \dd< 1/4$ let $C_\dd$ be the curve $\DD \cap \{|z_0 - z| = \dd\}$. The heart of the proof lies in showing that the length of $f(C_\dd)$ goes to 0 as $\dd \lar 0$. Let $L(\dd)$ denote the length of $C_\dd$. The magnification factor of a curve at each point $z$ is $|f'(z)|$, so we have

\begin{equation} \label{}
L(\dd) = \int_{C_\dd} |f'(z)| |dz|.
\end{equation}

It is also clear that the length of $C_\dd$ is bounded by $\pi \dd$. Thus, by the Cauchy-Schwarz inequality,

\begin{equation} \label{}
L(\dd)^2 \leq \Big( \int_{C_\dd} 1^2 |dz|\Big) \Big( \int_{C_\dd} |f'(z)|^2 |dz|\Big) \leq \pi \dd \int_{C_\dd} |f'(z)|^2 |dz|.
\end{equation}

We may therefore bound as follows.

\begin{equation} \label{}
\begin{split}
\int_{0}^{1/4} \frac{L(\dd)^2}{\dd} d\dd & \leq \pi \int_{0}^{1/4} \int_{C_\dd} |f'(z)|^2 |dz|d\dd \\
& = \pi \int \! \! \int_{\DD \cap D(z_0,1/4)} |f'(z)|^2 dA(z) \\
& = \pi Area(f(\DD \cap D(z_0,1/4))) < \ff.
\end{split}
\end{equation}

If there was an $\eps>0$ such that $L(\dd) \geq \eps$ for all $\dd$, we would have $\int_{0}^{1/4} \frac{L(\dd)^2}{\dd} d\dd \geq \int_{0}^{1/4} \frac{d\dd}{\dd} = \ff$, so we conclude that there is no such $\eps$. This does not tell us that $L(\dd) \lar 0$, but it does tell us that we can find a sequence $\dd_n \lar 0$ such that $L(\dd_n) \lar 0$ as well. Since the length of $f(C_{\dd_n})$ is finite, the curve reaches an endpoint in each direction. Let these endpoints be labeled $a_n, b_n$, and let the corresponding endpoints of $C_\dd$ be $\aa_n, \bb_n$. $a_n$ cannot lie within $W$, since if it did there would be a $z \in \DD$ such that $f(z) = a_n$, which by the open mapping theorem would imply that $f(D(z,\eps) \cap \DD) \cap f( D(\aa_n,\eps) \cap \DD) \neq \emptyset$ for any $\eps >0$, and this contradicts injectivity since $z \neq \aa_n$. Thus, $a_n \in \dd W$, and likewise for $b_n$.

\includegraphics[width=180mm,height=140mm]{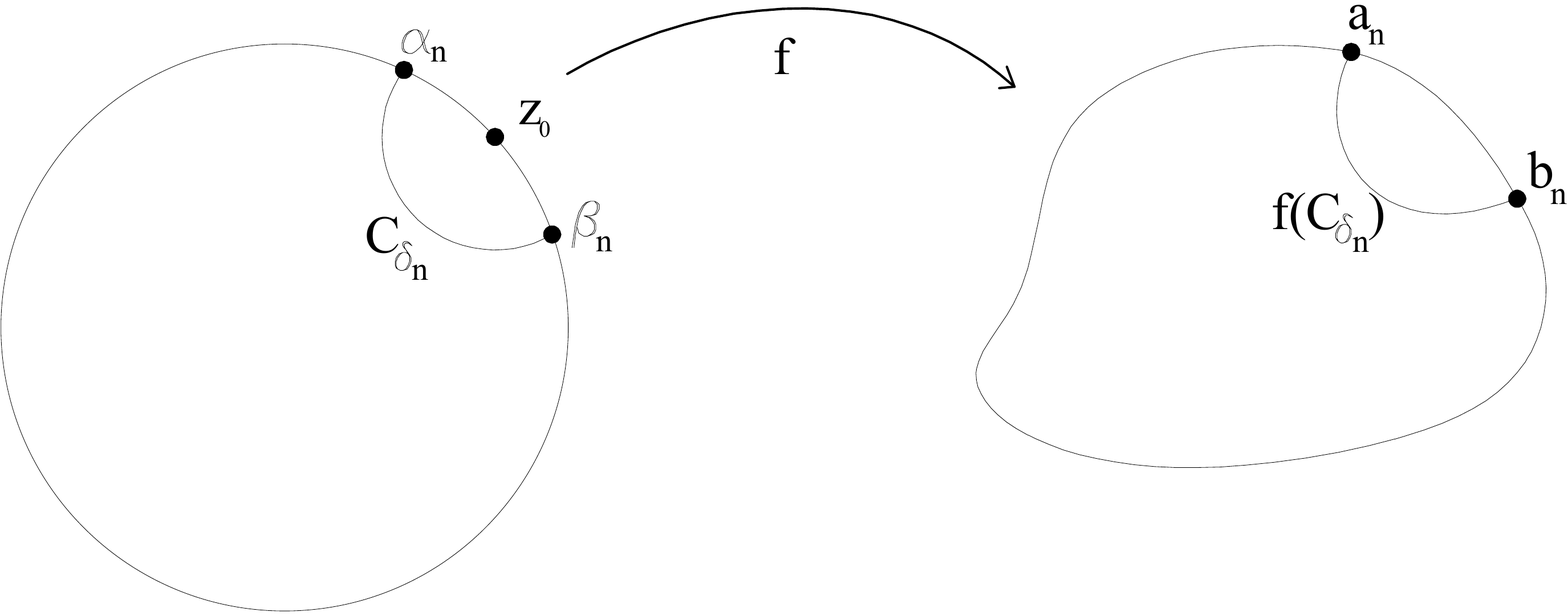}

\vspace{-1in}

Note that $|a_n - b_n| < diam(f(C_n)) \lar 0$. We can conclude from this that, if we let $D_n = D(z_0, \dd_n) \cap \DD$, we have $diam(f(D_n)) \lar 0$. The set $\seg{25}{f(D_n)}$ is compact and nonempty, and $\seg{35}{f(D_{n+1})} \subseteq \seg{25}{f(D_n)}$ for all $n$; Exercise \ref{excompint} below then implies that $\cap_{n=1}^\ff \seg{25}{f(D_n)} \neq \emptyset$. However, the diameter of this set must be less than the diameter of any $f(D_n)$, which implies that the diameter is 0. We conclude that the intersection consists of a single point, $w_0$. We define $f(z_0)$ to be $w_0$. It can now be checked that $f$, so extended to all of $\dd \DD$, defines a continuous bijection from $\seg{10}{\DD}$ to $\seg{10}{W}$. \qed

\section{Exercises}

\ccases{excompint} Suppose $K_n$ is a sequence of compact sets in $\CC$ with the property that $\cup_{n=1}^{N} K_n \neq \emptyset$ for any $N<\ff$. Show that $\cup_{n=1}^{\ff} K_n \neq \emptyset$. 

\setcounter{cccases}{0}
\chapter{Exit times and the growth of analytic functions} \label{exittimes}

\section{The expected exit time $T_W$ of Brownian motion} \label{secexpectexit}

Suppose that $W$ is a bounded domain, and as usual let $T_W = \inf\{t \geq 0: B_t \in W^c\}$ be the exit time of Brownian motion from $W$. We would like to understand the quantity $E_a[T_W]$ as a function of $a \in W$ as well as we can, and evaluate it if possible. Let us begin by noting that if we can find a real-valued function $h$ which is continuous on $\CCC(W)$, $C^2$ on $W$ with $\Del h = -2$, and 0 on $\dd W$, then we may apply It\=o's Theorem to obtain

\begin{equation} \label{}
-h(a) = E_a[h(B_{T_W})] - h(a) = \frac{1}{2} E_a \int_{0}^{T_W} \Del h(B_s) ds = - E_a \int_{0}^{T_W}ds = -E_a[T_W].
\end{equation}

Thus, $h(a) = E_a[T_W]$. However, the function $E_a[T_W]$ does not always satisfy these conditions; for example, if $W = D^{\times}(0,1) = \{0 < |z| < 1\}$, then since Brownian motion doesn't hit points we have $E_a[T_W] = E_a[T_\DD]$, and therefore $\lim_{a \lar 0} E_a[T_W] = E_0[T_\DD] = \frac{1}{2}$ even though $0 \in W^c$. On the other hand, the following holds.

\begin{proposition} \label{exitfcn}
If $W$ is a bounded simply connected domain, then $h(a) = E_a[T_W]$ is the unique continuous function on $\CCC(W)$ which is 0 on $\dd W$ and satisfies $\Del h = -2$ on $W$.
\end{proposition}

This follows from the properties of Green's function discussed earlier. There is another way to evaluate $E_a[T_W]$ which is of interest, as the next proposition shows.

\begin{proposition} \label{exitseries}
Suppose $W$ is a bounded simply connected domain, and $f : \DD \lar W$ is conformal with $f(0) = a$. Suppose further that $f$ admits the power series expansion $f(z) = a+\sum_{n=1}^{\ff} a_n z^2$. Then

\begin{equation} \label{hbp2}
E_a[T_W] = \frac{-|a|^2}{2} + \lim_{r \nearrow 1}\frac{1}{4\pi} \int_{0}^{2\pi} |f(re^{i\th})|^2 d\th = \frac{1}{2}\sum_{n=1}^{\ff}|a_n|^2.
\end{equation}
\end{proposition}

{\bf Proof:} Apply It\=o's Theorem to the function $u(z) = |z|^2$, using $\Del u = 4$, to obtain

\begin{equation} \label{hbp}
E_a[|B_{T_W}|^2] - |a|^2 = 2 E_a \int_{0}^{T_W} ds = 2 E_a[T_W].
\end{equation}

Now note that L\'evy's Theorem allows us to project a Brownian motion running on the disc onto $W$, and we obtain

\begin{equation} \label{}
E_a[|B_{T_W}|^2] = \lim_{r \nearrow 1} E_0[|f(B_{T_{D(0,r)}})|^2] = \lim_{r \nearrow 1} \frac{1}{2\pi} \int_{0}^{2\pi} |f(re^{i\th})|^2 d\th.
\end{equation}

Inserting this into \rrr{hbp} and dividing by 2 gives the first equality in \rrr{hbp2}, and the second follows from Parseval's Identity, which is proved next. \qed

\begin{theorem} \label{}
Suppose $f$ is an analytic function on $D(z_0,R) $ which admits the Taylor series expansion $f(z) = \sum_{n=0}^{\ff} a_n (z-z_0)^n$. Then

\begin{equation} \label{parid}
\frac{1}{2\pi} \int_{0}^{2\pi}|f_a(z_0 + r e^{i \th})|^2 d\th = \sum_{n=0}^{\ff} |a_n|^2 r^{2n}  \quad \mbox{for } r<R.
\end{equation}
\end{theorem}

{\bf Proof:} We have

\begin{equation} \label{}
\begin{split}
\int_{0}^{2\pi}|f_a(z_0 + r e^{i \th})|^2 d\th & = \int_{0}^{2\pi}f_a(z_0 + r e^{i \th})\seg{55}{f_a(z_0 + r e^{i \th})} d\th \\
& = \int_{0}^{2\pi} \Big(\sum_{n=0}^{\ff} a_n r^n e^{i n \th} \Big) \Big(\sum_{m=0}^{\ff} \bar a_m r^m e^{-i m \th} \Big) d\th \\
& = \sum_{n=0}^{\ff} \sum_{m=0}^{\ff} a_n \bar a_m r^{n+m} \int_{0}^{2\pi} e^{i (n-m) \th} d\th.
\end{split}
\end{equation}

\rrr{parid} now follows by noting $\int_{0}^{2\pi} e^{i (n-m) \th} d\th = 2\pi$ if $n=m$ and $\int_{0}^{2\pi} e^{i (n-m) \th} d\th = 0$ if $n \neq m$. \qed

Let us work a few examples using Proposition \ref{exitfcn}, Proposition \ref{exitseries}, or both. Note that one or the other may be substantially simpler in certain situations. \\

{\bf Example 1:} Let $W=D(0,R)$. The function $h(a) = \frac{R^2 - |a|^2}{2}$ satisfies the conditions of Proposition \ref{exitfcn}, and therefore $E_a[T_W] = \frac{r^2 - |a|^2}{2}$. This is also derivable from Proposition \ref{exitseries}, since if $f:\DD \lar W$ is the conformal map sending $0$ to $a$, then, by Carath\'eodory's Theorem, $f$ extends to $\dd \DD$ and $|f(z)| = R$ for $z \in \dd \DD$, so that

\begin{equation} \label{}
E_a[T_W] = \frac{-|a|^2}{2} + \lim_{r \nearrow 1}\frac{1}{4\pi} \int_{0}^{2\pi} |f(re^{i\th})|^2 d\th = \frac{R^2 - |a|^2}{2}.
\end{equation}

{\bf Example 2:} Let $W$ be the equilateral triangle with vertices at $1, \om, \om^2$, where $\om=e^{2\pi i/3}$. It can be shown that

\begin{equation} \label{rot}
h(a) = \frac{1}{6}(a + \bar a + 1)(\om a + \bar \om \bar a + 1)(\om^2 a + \bar \om^2 \bar a + 1)
\end{equation}

satisfies $\Del h = -2$ and $h(z) = 0$ for $z \in \dd W$: verification of the first statement is most easily accomplished by using $\triangle = 4\frac{\partial^2}{\partial a \partial \bar a}$, while the second is clear since $(a + \bar a + 1) = 0$ on $\{x=-\frac{1}{2}\}$ and $h(a) = h(\om a)$. Thus, $h(a) = E_a[T_W]$. As a simple special case, we get $E_0[T_W] = \frac{1}{6}$. \\

{\bf Example 3:} Let $W$ be the cardioid with boundary defined by the polar equation $r = 2(1+\cos \th)$.

\hspace{1.2in} \includegraphics[width=110mm,height=85mm]{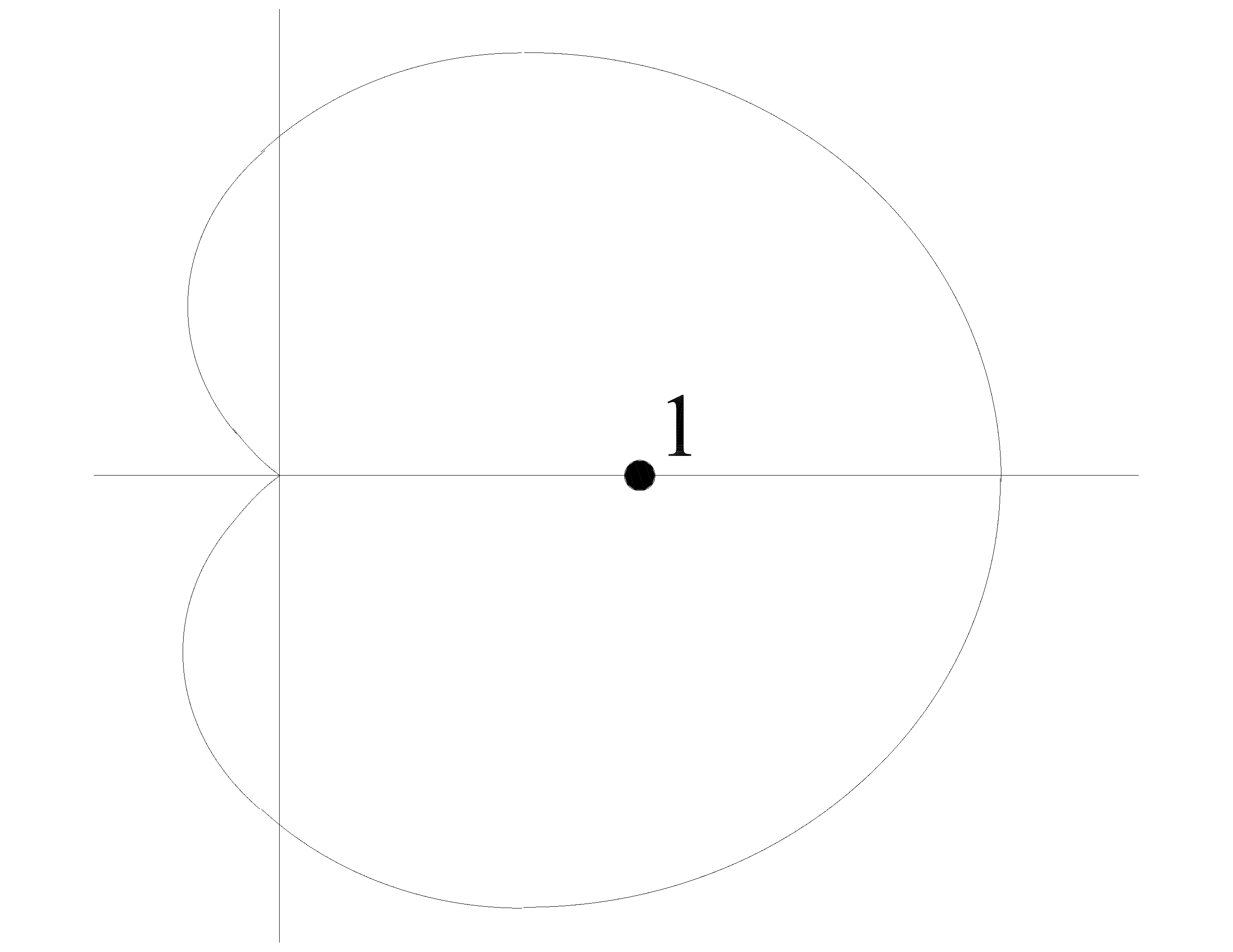}

This is the conformal image under $z^2$ of the disc $\{|z-1|=1\}$. The conformal map from $\DD$ to $W$ mapping $0$ to $1$ is therefore given by $f(z)=(z+1)^2 = 1+2z +z^2$. Applying Proposition \ref{exitseries} we get $E_1[T_W] = \frac{1}{2}(2^2 + 1^2) = \frac{5}{2}$. \\

{\bf Example 4:} Let $W$ be the infinite vertical strip $\{ \frac{-\pi}{4} < x < \frac{\pi}{4} \}$. It is easy to verify that $h(x+yi) = \frac{\pi^2}{16} - x^2$ is $0$ on $\dd W$ and satisfies $\Del h = -2$. If this truly gives $E_a[T_W]$, then we would have $E_0[T_W] = \frac{\pi^2}{16}$. On the other hand, we can find the conformal map from $\DD$ to $W$ which fixes 0 as follows. The function

\be
\tan z = \frac{\sin z}{\cos z} = -i \frac{e^{2iz} - 1}{e^{2iz} + 1}
\ee

maps $W$ conformally to $\DD$. This can be seen by noting that the function $x+iy \lar e^{2i(x+iy)} = e^{-2y + 2ix}$ maps $W$ conformally to $\{ \mbox{Re } z > 0\}$, and then that the M\"obius transformation $z \lar -i(\frac{z-1}{z+1})$ maps $\{ \mbox{Re } z > 0\}$ conformally to $\DD$. We conclude that the principal branch of $\tan^{-1}z$ maps $\DD$ conformally to $W$. $\tan^{-1}z$ admits the Taylor series expansion

\begin{equation} \label{}
\tan^{-1}z = z - \frac{z^3}{3} + \frac{z^5}{5} - \ldots = \sum_{n=1}^{\ff} \frac{(-1)^{n+1} z^{2n-1}}{2n-1}.
\end{equation}

Thus, if we may apply Proposition \ref{exitseries}, we would obtain

\begin{equation} \label{doubletrue}
E_0[T_W] = \frac{1}{2} \sum_{n=1}^{\ff} \frac{1}{(2n-1)^2}.
\end{equation}

Equating the values obtained for $E_0[T_W]$ in these two different manners yields

\begin{equation} \label{}
\sum_{n=1}^{\ff} \frac{1}{(2n-1)^2} = \frac{\pi^2}{8}.
\end{equation}

It is straightforward to verify that this is equivalent to Euler's celebrated Basel identity:

\begin{equation} \label{all}
\sum_{n=1}^\ff \frac{1}{n^2} = \frac{\pi^2}{6}.
\end{equation}

To steal a line from Chung \cite{chung2002green}, "since it works out at the end it must be justifiable"; however, it must be confessed that we have applied Propositions \ref{exitfcn} and \ref{exitseries} in a situation to which they do not directly relate, specifically to an unbounded domain $W$. To see that there really is an issue worth worrying about here, let us note that if we can find a nonconstant function $u$ harmonic on $W$ and $0$ on $\dd W$, then $h+u$ will be $0$ on $\dd W$ and satisfy $\Del(h+u) = -2$. It is also not hard to see that there always are such functions on unbounded domains. For instance, as discussed above the map $f(z) = e^{2iz}$ maps $W$ conformally to $\{ \mbox{Re } z > 0\}$, and thus the function $u(z) = Re(e^{2iz})$ is harmonic on $W$ and $0$ on $\dd W$. As it stands, we do not have a way of distinguishing between $h$ and $h+u$ in order to determine which, if either, truly gives us $E_a[T_W]$. Furthermore, the proof of Proposition \ref{exitseries} employed the identity $E_a[|B_{T_W}|^2] - |a|^2 = 2 E_a[T_W]$, but this does not necessarily hold if $W$ is unbounded. For example, let now $W = \{|z| > 1\}$, and $a \in W$. Then by recurrence we know that $T_W< \ff$ a.s., and therefore $E_a[|B_{T_W}|^2] - |a|^2 = 1-|a|^2 < 0$, but in fact it can be shown that $E_a[T_W] = \ff$! We therefore need to understand precisely when the identity $E_a[|B_{T_W}|^2] - |a|^2 = 2 E_a[T_W]$ holds. \\

The purpose of the next section is to resolve these issues.

\section{Unbounded domains}

The following lemma is basic but important for unbounded domains.

\begin{lemma} \label{}
\begin{itemize} \label{}

\item[(i)] Suppose $z \in W_1 \subseteq W_2$. Then $E_z[T_{W_1}] \leq E_z[T_{W_2}]$.

\item[(ii)] Suppose that $z \in W_1 \subseteq W_2 \subseteq \ldots$, and $W = \cup_{n=1}^\ff W_n$. Then $E_z[T_{W_n}] \nearrow E_z[T_{W}]$ as $n \lar \ff$.

\end{itemize}
\end{lemma}

{\bf Proof:} $(i)$ is obvious, and $(ii)$ is a simple consequence of the Monotone Convergence Theorem. \qed

The key to resolving the issues mentioned on unbounded domains is the following theorem, which is a special case of what we will refer to as {\it Burkholder's Theorem}, which can also be thought of as a planar form of the Davis Inequality. To state the theorem, we define the supremum process $|B^*_\tau| = \sup_{0 \leq t \leq \tau} |B_t|$ for any stopping time $\tau$. We then have

\begin{theorem} \label{burkthm}
For any stopping time $\tau$ we have

\begin{equation} \label{bteq}
2E_a[\tau] + |a|^2 \leq E_a[|B^*_\tau|^2] \leq 4(2E_a[\tau] + |a|^2).
\end{equation}
\end{theorem}

{\bf Proof:} First let us note that it is enough to prove the result in the case that $|B^*_\tau| \leq M$ a.s. This is because, if we do so, we can let $S_M = \inf\{t \geq 0:|B_t|=M\}$, apply the result to $\tau \wedge S_M$, and then let $M \nearrow \ff$, applying the Monotone Convergence Theorem to both $\tau \wedge S_M$ and $|B^*_{\tau \wedge S_M}|^2$ in order to obtain \rrr{bteq}. With this condition on $\tau$ we may apply It\=o's Theorem as in the proof of Proposition \ref{exitseries} to obtain $2E_a[\tau] + |a|^2 = E_a[|B_\tau|^2]$. It is clear that $E_a[|B_\tau|^2] \leq E_a[|B^*_\tau|^2]$, and we obtain the first inequality in \rrr{bteq}. The second follows similarly by utilizing the following important inequality.

\begin{proposition} \label{davisineq}
If $\tau$ is a stopping time with $E[\tau] < \ff$ then

\begin{equation} \label{}
E_a[|B^*_\tau|^2] \leq 4E_a[|B_\tau|^2].
\end{equation}
\end{proposition}

To prove this, we need several lemmas.

\begin{lemma} \label{}
\begin{equation} \label{}
E_a[|B_{S_M}|1_{S_M \leq \tau}] \leq E_a[|B_{\tau}|1_{S_M \leq \tau}]
\end{equation}
\end{lemma}

{\bf Proof:} Should follow from Dynkin's formula. \qed

\begin{lemma} \label{}
\begin{equation} \label{}
P_a(|B^*_\tau|\leq M) \leq \frac{1}{M} E_a[|B_\tau|1_{|B_\tau^*| \geq M}].
\end{equation}
\end{lemma}

{\bf Proof:}

\begin{equation} \label{}
P_a(|B^*_\tau|\leq M) = \frac{1}{M} E_a[|B_{S_M}|1_{S_M \leq \tau}] \leq \frac{1}{M} E_a[|B_{\tau}|1_{|B_\tau^*| \geq M}].
\end{equation}

\qed

\begin{lemma} \label{}
\begin{equation} \label{}
E[X^2] = \int_{0}^{\ff} 2M P(X > M) dM.
\end{equation}
\end{lemma}

{\bf Proof of Proposition \ref{davisineq}:}

\begin{equation} \label{}
\begin{split}
E_a[|B_\tau^*|^2 & = 2 \int_{0}^{\ff} M P(|B_\tau^*| > M) dM \\
& \leq 2 \int_{0}^{\ff} E_a[|B_\tau| 1_{|B_\tau^*| \geq M}]dM \\
& = 2 E_a[|B_\tau| \int_{0}^{|B_\tau^*|} dM] \\
& = 2 E_a[|B_\tau||B_\tau^*|] \\
& \leq E_a[|B_\tau|^2]^{1/2} E_a[|B_\tau^*|^2]^{1/2}.
\end{split}
\end{equation}

\qed

Armed with Burkholder's Theorem, we can impose a condition on harmonic functions on unbounded domains which will force them to be uniquely determined by their boundary values.

\begin{theorem} \label{phraglindharm}
Suppose $h$ is harmonic on a domain $W$ with $E[T_W] < \ff$ and continuous on $\CCC(W)$, with $h(z)= 0$ for $z \in \dd W$. If there is a constant $C>0$ such that $|h(z)| \leq C|z|^2 +C$, then $h(z) = 0$ for all $z \in W$.
\end{theorem}

{\bf Proof:} As in the proof of Theorem \ref{burkthm}, let $S_M = \inf\{t \geq 0:|B_t|=M\}$, and fix $a \in W$. By It\=o's Theorem, $h(a) = E_a[h(B_{T_W \wedge S_M})]$. We would like let $m \lar \ff$ to obtain $h(a) = E_a[h(B_{T_W})]=0$; note that the conditions on $h$ imply that $|h(B_{T_W \wedge S_M})| \leq C|B_{T_W \wedge S_M}|^2 + C \leq C|B^*_{T_W}|^2 + C$, and $E_a[C|B^*_{T_W}|^2 + C] \leq 4C (2E_a[T_W] + |a|^2) + C < \ff$ by Burkholder's Theorem. The Dominated Convergence Theorem therefore applies, and we get $h(a) = \lim_{M \nearrow \ff} E_a[h(B_{T_W \wedge S_M})] = E_a[h(B_{T_W})]=0$. \qed

This result suggests that our problem of determining the correct exit time function $h(a) = E_a[T_W]$ may be solved by imposing a growth condition on $h$, and this turns out to be correct.

\begin{theorem} \label{}
Suppose $W$ is a domain and $h$ is a nonnegative continuous function on $\CCC(W)$ such that $\Del h = -2$ on $W$ and $h = 0$ on $\dd W$. Suppose further that there is a constant $C>0$ such that $h(z) \leq C|z|^2 +C$. Then $E[T_W] < \ff$, and $h(a) = E_a[T_W]$ for all $a \in W$.
\end{theorem}

{\bf Proof:} As before, let $S_M = \inf\{t \geq 0:|B_t|=M\}$, and apply It\=o's Theorem to $-h$ to obtain

\begin{equation} \label{danp}
E_a[T_W \wedge S_M] = \frac{1}{2} \int_{0}^{T_W \wedge S_M} \Del (-h(B_s))ds = h(a) - E_a[h(B_{T_W \wedge S_M})] \leq h(a).
\end{equation}

Letting $M \nearrow \ff$, we obtain $E_a[T_W] \leq h(a) < \ff$. Once we know this, Burkholder's Theorem tells us that $E_a[|B_{T_W}^*|^2] < \ff$, which allows us to bound $h(B_{T_W \wedge S_M}) \leq C|B_{T_W \wedge S_M}|^2 + C \leq C|B^*_{T_W}|^2 + C$ and apply the Dominated Convergence Theorem in \rrr{danp} to obtain

\begin{equation} \label{}
E_a[T_W] = \lim_{M \nearrow \ff} E_a[T_W \wedge S_M] = h(a) - \lim_{M \nearrow \ff} E_a[h(B_{T_W \wedge S_M})] = h(a).
\end{equation}

\qed

We remark that the condition that $h$ be nonnegative is very important: in the example at the end of Section \ref{secexpectexit}, the function $h(a) = 1-|a|^2$ satisfies all of the other conditions in the theorem. The following theorem is a fairly general version of what is called the {\it Phragm\'en-Lindel\=of principle}. This principle is essentially a version of the maximum principle which is valid for functions which grow slowly on unbounded domains.

\begin{theorem} \label{phraglind1}
Let $W$ be a domain such that $E[T_W] < \ff$. Suppose that $f$ is an analytic function on $W$ such that $\limsup_{z \lar \dd W}|f(z)| \leq K < \ff$, and $|f(z)| \leq Ce^{C|z|^2} + C$. Then $|f(z)| \leq K$ for all $z \in W$.
\end{theorem}

We will need the following lemma. As usual, set $S_M = \inf\{t \geq 0:|B_t|=M\}$, and define

\begin{equation} \label{17dd}
\log^+ x = \begin{cases} \log x & x > 1,\\
0 & x  \leq 1.
\end{cases}
\end{equation}

\begin{lemma} \label{logpluslem}
If a function $f$ is analytic on a domain $W$ and continuous on $\CCC(W)$, then $\log^+ |f(a)| \leq E_a[\log^+|f(B_{T_W \wedge S_M})|]$.
\end{lemma}

{\bf Proof of Lemma \ref{logpluslem}:} We note first that, since $\log x$ and $g(x) = 0$ are both harmonic we have $\log |f(a)| = E_a[\log|f(B_{T_W \wedge S_M})|]$ and $g(|f(a)|) = E_a[g(|f(B_{T_W \wedge S_M})|)]$. Since $\log^+ = \max(\log, g)$, for any $a$ such that $\log^+ |f(a)| = \log |f(a)|$ we have $\log^+ |f(a)| = \log |f(a)| = E_a[\log|f(B_{T_W \wedge S_M})|] \leq E_a[\max(\log,g)(f(B_{T_W \wedge S_M}))] = E_a[\log^+(f(B_{T_W \wedge S_M}))]$. A similar calculation holds if $\log^+ |f(a)| = g(|f(a)|)$. \qed

{\bf Proof of Theorem \ref{phraglind1}:} By multiplying $f$ by a constant we may assume that $K=1$. The assumption on $f$ shows that $\log^+ |f(B_{T_W \wedge S_M})| \leq C|B_{T_W \wedge S_M}|^2 + C \leq C|B^*_{T_W \wedge S_M}|^2 + C$, and since $E[T_W] < \ff$ Burkholder's Theorem tells us that $E_a[C|B^*_{T_W \wedge S_M}|^2 + C] < \ff$. We may therefore apply the Dominated Convergence Theorem and Lemma \ref{17dd} to conclude

\begin{equation} \label{}
\log^+ |f(a)| \leq \lim_{M \nearrow \ff} E_a[\log^+|f(B_{T_W \wedge S_M})|] = E_a[\log^+|f(B_{T_W})|] = 1.
\end{equation}

Thus, $|f(a)| \leq 1$. \qed

\section{Moments of $T_W$} \label{genmomsec}

Burkholder's Theorem admits a more general formulation.

\begin{theorem} \label{burkthm}
For any stopping time $\tau$ and $p \in (0,\ff)$ there are constants $c_p, C_p > 0$ such that

\begin{equation} \label{bteq}
c_pE_a[(\tau + |a|^2)^p] \leq E_a[|B^*_\tau|^{2p}] \leq C_p E_a[(\tau + |a|^2)^p].
\end{equation}

In particular, $E_a[\tau^p] < \ff$ if, and only if, $E_a[|B^*_\tau|^{2p}]< \ff$.
\end{theorem}

The proofs of Theorems \ref{phraglindharm} and \ref{phraglind1} are easily adapted, using this more general result, and yield the following results.

\begin{theorem} \label{phraglindharm2}
Suppose $h$ is harmonic on a domain $W$ with $E[T^p_W] < \ff$ and continuous on $\CCC(W)$, with $h(z)= 0$ for $z \in \dd W$. If there is a constant $C>0$ such that $|h(z)| \leq C|z|^{2p} +C$, then $h(z) = 0$ for all $z \in W$.
\end{theorem}

\begin{theorem} \label{phraglind2}
Let $W$ be a domain such that $E[T^p_W] < \ff$. Suppose that $f$ is an analytic function on $W$ such that $\limsup_{z \lar \dd W}|f(z)| \leq K < \ff$, and $|f(z)| \leq Ce^{C|z|^{2p}} + C$. Then $|f(z)| \leq K$ for all $z \in W$.
\end{theorem}

\section{Exercises}

\ccases{idententire} a) Find all entire functions $f(z)$ such that $f(1/n) = 1/n^2$ for all $n \in \NN$.

b) Find all entire functions $f(z)$ such that $f(1/n) = 1/\sqrt{n}$ for all $n \in \NN$.

\vski

\ccases{openaxes} Suppose $f(z)$ is an analytic function on a domain $W$ with the property that $f(z)^n \in \RR$ for all $z \in W$. Show that $f$ is a constant function.

\vski

\ccases{lindprep} Suppose that $\ga$ is a {\it Jordan curve}: that is, $\ga$ can be realized as the image of a continuous and one-to-one function $\ga(t)$ from $[0,1]$ into $\CC$. $\ga$ then determines a bounded domain $W$ in its interior. Suppose that $f(z)$ is an analytic function on $W$, and suppose that there are points $w_1, \ldots , w_n \in \ga$ such that $\limsup_{\{z \lar \dd W \backslash \{w_1, \ldots ,w_n\}\}} |f(z)| \leq m < \ff$. Suppose further that there is some finite $M > m$ such that $|f(z)| \leq M$ for all $z \in W$. Show that in fact $|f(z)| \leq m$ for all $z \in W$.

\hspace{1.2in} \includegraphics[width=90mm,height=70mm]{Lindelof.pdf}

\vski

\ccases{rotatethis} Let $W$ be the square domain $\{-1<x<1,-1<y<1\}$. Suppose $f(z)$ is a bounded analytic function on $W$ such that $f(z) \lar 0$ as $z$ approaches the bottom edge of the square; that is, $\lim_{\{z \lar \{y=-1\}\}} |f(z)| = 0$. Prove that $f(z) = 0$ for all $z \in W$.

\vski

\ccases{k-roots} Suppose $p(z) = \bb_k z^k + \ldots + \bb_0$ is a degree $k$ polynomial, with $\bb_0, \ldots, \bb_k \in \CC$ and $\bb_k \neq 0$. Show that $p$ has $k$ roots, counting multiplicities. That is, show that there are values $a_1, \ldots, a_k, c \in \CC$ such that $p(z) = c(z-a_1)(z-a_2)\ldots (z-a_k)$. (Hint: If $p(a) = 0$, consider $p(z)-p(a)$).

\vski

\ccases{intprop} i) (Independence of parameterizations) Suppose $\ga_1: [a_1,b_1] \lar \CC$ and $\ga_2: [a_2,b_2] \lar \CC$ are rectifiable curves, and suppose that there is a smooth increasing function $\phi : [a_1,b_1] \lar [a_2,b_2]$ with $\phi(a_1) = a_2, \phi(b_1) = b_2$ such that $\ga_1(t) = \ga_2(\phi(t))$ for all $t \in [a_1,b_1]$. Show that

\begin{equation} \label{}
\int_{\ga_1} f(z)dz = \int_{\ga_2} f(z)dz
\end{equation}
for any analytic function $f$ defined in a neighborhood of $\ga_1=\ga_2$.

\vski

ii) (Path reversal) Let $\ga: [a_1,b_1] \lar \CC$ be a rectifiable curve. Define $\tilde \ga : [-b_1,-a_1] \lar \CC$ by $\tilde \ga (t) = \ga (-t)$. Then $\tilde \ga$ is simply the curve $\ga$ traversed in the opposite direction. Show that

\begin{equation} \label{}
\int_{\tilde \ga} f(z)dz = - \int_{\ga} f(z)dz
\end{equation}
for any analytic function $f$ defined in a neighborhood of $\ga$.

\vski

iii) (Concatenating paths) Suppose $\ga_1: [a,b] \lar \CC$ and $\ga_2: [c,d] \lar \CC$ are rectifiable curves such that $\ga_1(b) = \ga_2(c)$. Define a new curve $\tilde \ga: [a, b+d-c] \lar \CC$ by

\begin{equation} \label{17}
\tilde \ga(t) := \begin{cases} \ga_t(t) & t \in [a,b],\\
\ga_2(t-b+c) & t \in [b,b+d-c].
\end{cases}
\end{equation}
$\tilde \ga$ is the concatenation of $\ga_1$ and $\ga_2$. Show that

\begin{equation} \label{}
\int_{\tilde \ga} f(z)dz = \int_{\ga_1} f(z)dz + \int_{\ga_2} f(z)dz
\end{equation}
for any analytic function $f$ defined in a neighborhood of $\tilde \ga$.

\vski

\ccases{compdist} Let $d(z,\dd W)$ be as in Section \ref{CTsec}. Show that the function $z \lar d(z,\dd W)$ is continuous and positive on $W$. Conclude that if $K$ is a compact set in $W$ then there is some $\eps > 0$ such that $d(z,\dd W) > 0$ for all $z \in K$.

\vski

\ccases{radconvsame} Show that the radius of convergence is the same for the three series $\sum_{n=0}^{\ff} \bb_n (z-z_0)^n, \sum_{n=0}^{\ff} \frac{\bb_n}{n+1} (z-z_0)^{n+1},$ and $\sum_{n=1}^{\ff} n \bb_n (z-z_0)^{n-1}$.

\vski

\ccases{polytriang} Polygon admits triangulation.

{\bf \Large 1.} Let $\ga = \{e^{it}: 0 \leq t \leq 2\pi\}$. Evaluate the following two integrals.

\begin{equation} \label{}
\int_{\ga} \frac{1}{z}dz, \qquad \int_{\ga} \frac{1}{z^2}dz.
\end{equation}

{\bf \Large 2.} Let $\ga(t): [0, 5\pi]$ be defined as follows.

\begin{equation} \label{17}
\ga(t) = \begin{cases} (t+\frac{3\pi}{2})e^{it} & 0 \leq t \leq 2 \pi,\\
\frac{3\pi}{2} + 2\pi e^{-it} & 2\pi \leq t \leq 3 \pi,\\
(t-\frac{7 \pi}{2}) + i \cos (t-\frac{7 \pi}{2}) & 3 \pi \leq t \leq 5\pi.
\end{cases}
\end{equation}

Sketch $\ga$, and evaluate

\begin{equation} \label{}
\int_{\ga} e^{1/z} dz .
\end{equation}

\ccases{exlineremove} Suppose $W$ is a domain and $C$ is a line intersecting $W$. Suppose further that $f$ is a continuous function on $W$ which is analytic on $W \bsh C$. Show that in fact $f$ is analytic on all of $W$. (Hint: Morera's Theorem)

\setcounter{cccases}{0}
\chapter{Aspects of Stochastic Loewner Evolution} \label{exittimes}

One of the most exciting developments in probability theory in the early $21$-st century has been the discovery of Stochastic Loewner Evolution by Oded Schramm, as well as the subsequent related work by other authors.

\section{Schwarz's lemma and consequences} \label{secexpectexit}

We will now work towards a discussion of Stochastic Loewner Evolution. We begin with an important property of harmonic functions can be proved using Brownian motion.

\begin{theorem}{\rm [Mean Value Property]}  \label{}
Suppose $u$ is harmonic on a domain $W$, and the set $\{|z-z_0| \leq r\} \subseteq W$. Then

\begin{equation} \label{}
u(z_0) = \frac{1}{2\pi} \int_{0}^{2\pi} u(z_0 + re^{i\th}) d \th.
\end{equation}

\end{theorem}

That is, $u(z_0)$ is the average of the values of $u$ near $z_0$. A consequence is the following.

\begin{theorem}{\rm [Maximum/minimum principle]} \label{maxmodthm}
Suppose $\Cc$ is a simple closed curve with interior $U$. If $h(z)$ is harmonic on a domain $D$ containing $\Cc$ and $U$,
then
$$
\min_{w\in \Cc} u(w) \leq u(z) \leq \max_{w\in \Cc} u(w)
$$
for all $z\in U$.  In fact, if $u$ is nonconstant, then

$$
\min_{w\in \Cc} u(w) < u(z) < \max_{w\in \Cc} u(w)
$$
for all $z\in U$.
\end{theorem}

From this follows another important property of analytic functions.

\begin{theorem}{\rm [Maximum Modulus Theorem]} \label{maxmodthm}
Suppose $\Cc$ is a simple closed curve with interior $U$. If $f(z)$ is analytic in a domain $D$ containing $\Cc$ and $U$,
then
$$
|f(w)| \leq \max_{z\in \Cc} |f(z)|
$$
for all $w\in U$.  If $f$ is not a constant, then in fact $|f(w)| < \max_{z\in \Cc} |f(z)|$ for all $w\in U$.
\end{theorem}

Note that there is no analog of this result in real calculus. This implies the following important result, known as {\it Schwarz's Lemma}.

\begin{theorem} \label{}
Suppose $f(z)$ is an analytic function defined on $\DD$ such that $f(\DD) \subseteq \DD$, i.e. $|f(z)| < 1$ for all $z \in \DD$. Suppose also that $f(0) = 0$. Then $|f(z)| \leq |z|$ for all $z \in \DD$, and $|f'(0)| \leq 1$. Furthermore, if equality holds in either case then $f(z) = cz$, where $c$ is a complex number with $|c| = 1$.
\end{theorem}

{\bf Proof:} Consideration of the power series expansion $f(z) = \sum_{n=1}^{\ff} a_n z^n$ with $a_n = \frac{f^{(n)}(0)}{n!}$ shows that the function $g(z) = \frac{f(z)}{z}$ and $g(0) = f'(0)$ is analytic on $\DD$. Note that $\limsup_{z \lar \dd \DD} |g(z)| = \limsup_{z \lar \dd \DD} |f(z)| \leq 1$, since $f(\DD) \subseteq \DD$, so the maximum principle implies that $|g(z)| \leq 1$, with equality only when $g$ is a constant. The result follows. \qed

From Schwarz's Lemma comes the following important fact.

\begin{theorem} \label{}
Suppose $U_1, U_2$ are simply connected domains in $\CC$ and $a \in \CC$ such that $a \in U_1 \subseteq U_2$. Let $f_1(z), f_2(z)$ be conformal maps from $\DD$ onto $U_1, U_2$ respectively, such that $f_1(0) = f_2(0) = a$. Then $|f_1'(0)| \leq |f_2'(0)|$, and equality holds if and only if $U_1 = U_2$.
\end{theorem}

\section{de Branges' theorem and deterministic Loewner evolution}

In the 1800's, with the Riemann Mapping Theorem and other developments in complex analysis people began studying conformal maps deeply.  Let $\SS$ be the set of all conformal maps $f$ defined on the disk with $f(0) = 0, f'(0) = 1$.  $\SS$ is called the {\it Schlicht class}, and note that any conformal map $g$ can be transformed into a Schlicht function $f(z) = \frac{g(z) - g(0)}{g'(0)}$. Here are two examples of Schlicht functions

\begin{itemize} \label{}

\item $I(z) = z$, the identity.

\item $K(z) = \frac{z}{(1-z)^2}$, the Koebe function.

\end{itemize}

If $f$ is a Schlicht function we can write $f(z) = \sum_{n=0}^{\ff} a_n z^n = z + \sum_{n=2}^{\ff} a_n z^n$.  Note that for the Koebe function $K(z) = z + 2z^2 + 3z^3 + 4z^4 + \ldots$.  In 1916, Bieberbach proved that always $|a_2| \leq 2$ for Schlicht functions, and conjectured that always $|a_n| \leq n$.  In 1923, Loewner developed a remarkable method for proving $|a_3| \leq 3$, and his method led eventually to the full result.

\begin{theorem}[de Branges Theorem, also known as the Bieberbach Conjecture] \label{}
For any Schlicht function $f$, $|a_n| \leq n$, and equality holds for any $n$ only for the Koebe function and its rotations.
\end{theorem}

Loewner's idea was to look at only {\it slit domains}, which are domains of the form $\CC \backslash \ga$, where $\ga(t)$ is a continuous 1-to-1 function from $[0,\ff)$ to $\CC$ such that $\ga(t) \lar \ff$ as $t \lar \ff$.  For such a domain $U$ containing $0$, let $U_t$ be $\CC \backslash \ga[t,\ff]$, and let $f_t$ be the conformal map from $\DD$ to $U_t$ with $f_t(0) = 0$ and $f_t'(0) > 0$.  Then $f_t'(0)$ is increasing in $t$, and we may reparameterize $\ga$ so that $f_t'(0) = e^t$.  Let $g_t = f_t^{-1}$. We then have the following theorem.

\begin{theorem} \label{}

$f_t$ and $g_t$ satisfy the following differential equations, where $\la(t) = g_t(\ga(t))$.

\begin{equation} \label{dharry}
\begin{gathered}
\frac{\partial}{\partial t} g_t(w) = -g_t(w) \frac{\la(t) +  g_t(w)}{\la(t) - g_t(w)}; \\
\frac{\partial}{\partial t} f_t(z) = f_t'(z) z \frac{\la(t) - z}{\la(t) + z}; \\
\end{gathered}
\end{equation}

\end{theorem}

{\bf Sketch of proof:} Let $s < t$, let $\psi_{st}(z) = g_t(f_s(z))$, let $C_{st} = g_s(\ga[s,t])$, and let $J_{st} = g_t(\ga[s,t])$; see Figure 1 below for a picture of this setup. For the same reason as in Schwarz's Lemma, since $\psi(0) = 0$ the function $k(z) = \frac{\psi_{st}(z)}{z}$ for $z \neq 0$, $k(0) = \psi_{st}'(0) = e^{s-t}$, is analytic on $\DD$, and since it is never equal to 0 we can define an analytic logarithm $K(z)$ such that $e^{K(z)} = k(z)$ (this is a standard result from complex analysis theory, although we didn't prove it in class). We may also, by adding an integer multiple of $2\pi i$ if necessary, specify that $K(0) = s-t$. Now, the theory of conformal mappings shows that $\psi_{st}(z)$ extends continuously to the boundary of $\DD$, and we may therefore consider the harmonic function $Re(K(z))$ and apply the Poisson Integral Formula in the following form (see the remark immediately following Theorem \ref{PIF}):

\begin{figure}[h!] %  figure placement: here, top, bottom, or page
  \centering
  \includegraphics[width=4.2 in,height=4 in]{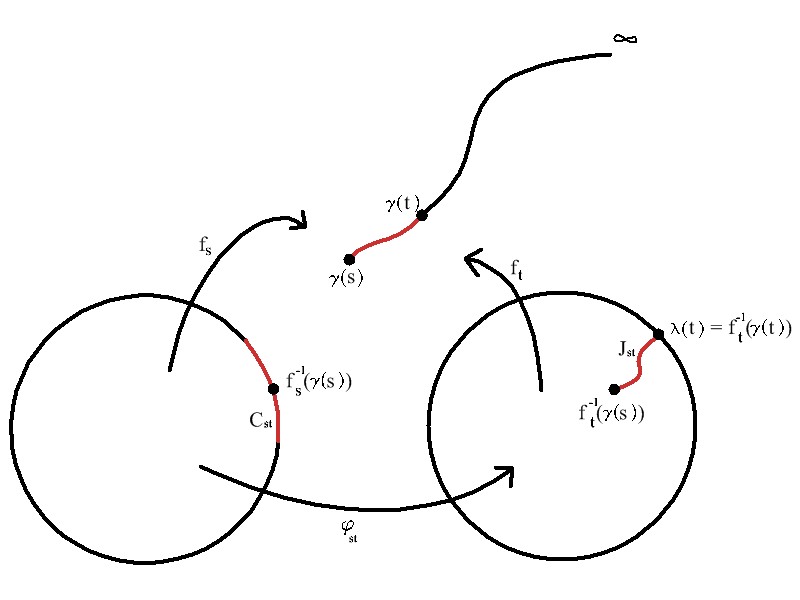}
\caption {Derivation of the Loewner equation}
\end{figure}

$$
Re(K(z)) = \frac{1}{2\pi} \int_{0}^{2\pi} Re(K(e^{i\th})) Re \Big( \frac{e^{i\th} + z}{e^{i\th} - z} \Big) d\th.
$$

Consider now the function $\frac{1}{2\pi} \int_{0}^{2\pi} Re(K(e^{i\th})) \Big( \frac{e^{i\th} + z}{e^{i\th} - z} \Big) d\th$; Leibniz's integral rule allows us to differentiate inside the integral, and since the integrand is analytic in $z$ we conclude that this function is analytic as well; since it has the same real part as $K(z)$ we see that these two analytic functions must differ only by an imaginary constant (this follows from the uniqueness of harmonic conjugates); putting $z=0$ shows that this imaginary constant is 0, and we obtain

$$
K(z) = \frac{1}{2\pi} \int_{0}^{2\pi} Re(K(e^{i\th})) \Big( \frac{e^{i\th} + z}{e^{i\th} - z} \Big) d\th.
$$

Note however that $Re(K(z)) = \ln \Big|\frac{\psi_{st}(z)}{z}\Big|$, and that $|\psi_{st}(z)|= 1$ on $\dd \DD$ except for when $z \in C_{st}$ (because $\psi_{st}(C_{st}) = J_{st}$). Thus, $Re(K(z)) = 0$ for all $z$ on $\dd \DD$ which do not lie on $C_{st}$. We obtain

$$
K(z) = \frac{1}{2\pi} \int_{\aaa}^{\bb} Re(K(e^{i\th})) \Big( \frac{e^{i\th} + z}{e^{i\th} - z} \Big) d\th,
$$

where $e^{i\aaa}$ and $e^{i\bb}$ are the endpoints of $C_{st}$. Now, as $t \searrow s$, the interval $C_{st}$ shrinks to the point $\la(s)$, and therefore the value of $\frac{e^{i\th} + z}{e^{i\th} - z}$ will be close to $\frac{la(s) + z}{\la(s) - z}$, and we will have approximately (this step is rigorously justified by the mean value theorem for integrals)

$$
K(z) \approx  \Big(\frac{\la(s) + z}{\la(s) - z}\Big) \Big( \frac{1}{2\pi}\int_{\aaa}^{\bb} Re(K(e^{i\th})) d\th\Big) = \Big(\frac{\la(s) + z}{\la(s) - z}\Big) Re(K(0)) = \Big(\frac{\la(s) + z}{\la(s) - z}\Big) (s-t),
$$

where the mean value property of harmonic functions was used. If we let $w = f_s(z)$, then it can be checked that $K(z) = \ln \frac{g_t(w)}{g_s(w)}$, and therefore

$$
\frac{\ln(g_t(w)) - \ln(g_s(w))}{t-s} = -\Big(\frac{\la(s) + g_s(w)}{\la(s) - g_s(w)}\Big).
$$

If we take the limit now as $t \searrow s$, the left side of the equation converges to the logarithmic derivative $\frac{\partial}{\partial s} \ln (g_s(w)) = \frac{\frac{\partial}{\partial s}g_s(w)}{g_s(w)}$, and the first equation in \rrr{dharry} follows; in fact, we have proved it only for the right hand derivative, but essentially the same argument shows gives proves it for the left hand derivative, and the second differential equation in \rrr{dharry} follows by inverting the first. \qed

In fact, in the proof above we used the obvious-sounding fact that $\la(t)$ is continuous in $t$; this is a true statement but it is not easy to prove. Full details can be found in the book {\it Conformal invariants: Topics in geoetric function theory} by Lars Ahlfors.

\section{Stochastic Loewner evolution} \label{genmomsec}

In 2000, Oded Schramm realized that the Loewner differential equation was also useful in studying the evolution of random processes which don't self-intersect.  That is, in the Loewner equation $\la(t)$ is defined in terms of $\ga(t)$, but the reverse works as well, and if we put a random $\la(t)$ into the Loewner equation we obtain a random output $\ga(t)$ which does not self-intersect.  The natural choice for $\la(t)$ is a 1-dimensional Brownian motion, and the resulting random curve $\ga(t)$ has been shown to be the continuous limit of a number of different discrete processes. In fact, what probabilists generally look at is a bit of a modification of the original Loewner evolution, with the curve growing into a specified domain.  A common choice is the unit disk $\DD = \{|z|<1\}$, and the same general method as was used in the plane leads to the following differential equations.

\begin{equation} \label{}
\begin{gathered}
\frac{\partial}{\partial t} g_t(w) = g_t(w) \frac{\la(t) + g_t(w)}{\la(t) - g_t(w)}; \\
\frac{\partial}{\partial t} f_t(z) = -f_t'(z) z \frac{\la(t) + z}{\la(t) - z}; \\
\end{gathered}
\end{equation}

The other common choice is the upper half-plane $\HH = \{Im(z)>0\}$, and again the same method leads to the following differential equations.

\begin{equation} \label{}
\begin{gathered}
\frac{\partial}{\partial t} g_t(w) = \frac{2}{g_t(w) - \la(t)}; \\
\frac{\partial}{\partial t} f_t(z) = \frac{2f_t'(z)}{\la(t) - z}; \\
\end{gathered}
\end{equation}

As a simple example, if $\la(t) = 0$ then

\begin{equation*} \label{}
\begin{gathered}
f_t(z) = \sqrt{z^2 - 4t} \\
g_t(z) = \sqrt{z^2 + 4t} \\
\ga(t) = 2i \sqrt{t}
\end{gathered}
\end{equation*}

For $SLE$, we take $\la(t) = B_{\kappa t} (=\sqrt{\kappa}B_t$), and call this $SLE_k$. Some basic properties of the curve for $SLE_k$ are as follows.

\begin{itemize} \label{}

\item For $0 \leq \kappa \leq 4$ the curve $\ga(t)$ is simple a.s.

\item For $4 < \kappa < 8$ the curve $\ga(t)$ intersects itself and every point is contained in a loop but the curve is not space-filling a.s.

\item For $\kappa \geq 8$ the curve $\ga(t)$ is space-filling a.s.

\end{itemize}

$SLE$ has been identified as the scaling limit of a number of discrete processes:

\begin{itemize} \label{}

\item $\kappa = 2$ is the scaling limit of the loop-erased random walk.

\item $\kappa = \frac{8}{3}$ is the outer boundary of Brownian motion.

\item $\kappa = 3$ is the limit of interfaces for the Ising model.

\item $\kappa = 4$ is the scaling limit of the path of the harmonic explorer and contour lines of the Gaussian free field.

\item $\kappa = 6$ is the scaling limit of critical percolation on the triangular lattice.

\item $\kappa = 8$ corresponds to the path separating the uniform spanning tree from its dual tree.

\end{itemize}

A number of others are conjectured, but are still unproved.
\chapter{Further reading} \label{references}

This section contains references to all papers and books I know of related to the connection between planar Brownian motion and complex analysis. I probably will update it from time to time, so if the reader knows of any I have missed (there certainly are some) please let me know.

\vski

The books \cite{bass1994probabilistic} and \cite{durBM} contain excellent overviews of the area, but are both rather less focused than this document (both devote significant time to $\RR^n$ with $n>2$, for example). Also highly recommended is the survey \cite{davis}, which is devoted exclusively to complex plane, but which of course is missing many recent results. The texts \cite{doob2012classical} and \cite{petersen1977brownian} also contain a number of interesting connections, though both are quite focused and technical. A highly enjoyable related text is \cite{chung2002green}, which benefits from being rather informal. Many general texts on Brownian motion and stochastic calculus contain sections related to conformal invariance, including \cite{bishop2017fractals, morters2010brownian, pitman2018guide, revyor,stroock2010probability}.

\vski

The solution to the Dirichlet problem presented in Chapter \ref{altrings} is from \cite{markowsky2018remark}, and the proof of the open mapping theorem in the same chapter is from \cite{markowsky2014probabilistic}. The remarkable probabilistic proof of Picard's theorem originated in \cite{davis1975picarda}, with variants appearing in \cite{davis}, \cite{durBM}, and \cite{bass1994probabilistic}. The proof of winding used in this text can be found in \cite{mckean1969stochastic}. Liouville's theorem, and thence the fundamental theorem of algebra, can be deduced easily from Levy's theorem and the bounded martingale convergence theorem, as has been mentioned in \cite{rogers2000diffusions} and doubtless in other places. Another proof of the fundamental theorem of algebra using Brownian motion can be found in \cite{pascu2005probabilistic}.

\vski

A {\it Skorokhod embedding} of Brownian motion is (loosely stated) to find a stopping time $\tau$ (satisfying certain conditions) such that $B_\tau$ has a given prescribed distribution. In one dimension this problem has generated a large amount of interest; see for instance the extensive survey \cite{obloj2004skorokhod}. Recently, it was formulated in two dimensions and solved in \cite{gross2019conformal}, with a method making extensive use of conformal mapping. This has led to the follow-up works \cite{boudabra2020new,boudabra2020remarks,mariano2020conformal}.

\vski

Planar Brownian motion has an interesting connection to the Cauchy distribution with parameters $(a,b)$, which has density function given by
\[
\frac{b}{\pi(b^2+(x-a)^{2})}
\]
This is the density of the distribution of $B_\tau$, where $B_0 = a+bi$ a.s. with $b \neq 0$ and $\tau$ is the hitting time of $\RR$; a recent proof of this using the optional stopping theorem appears in \cite{chin2019some}, but it can also be deduced by a direct calculation, using the Poisson kernel, by properties of stable distributions, or by a conformal map from the disk to the upper half-plane; see for instance \cite[Sec. 1.9]{durBM}, \cite[Ch VI.2]{feller2008introduction}, or \cite{markowskyfenn}. This fact indicates a fascinating property of the Cauchy distribution, that it must be preserved (up to a change in parameters) by an analytic function taking the upper half-plane to itself and the real line to itself. This property has been explored and exploited in a number of different places, either with or without explicit reference to Brownian motion. These sources include \cite{bourg,chin2019some,letac1977functions,MR1394988,okamura2020characterizations,pitman1967cauchy}. 

\vski

Defining Green's function using Brownian motion allows for a number of infinite product identities to be deduced, including Euler's celebrated infinite product representation for the sine function. This is shown in \cite{megreen}, while the corresponding analytic derivations of some of these identities can be found in \cite{melgreen} and \cite{melgreen2}.

\vski

It is clear that if $a \in U \subseteq V$, then $E_a[T(U)^p] \leq E_a[T(V)^p]$ for all $p \in (0,\ff)$, and in some sense the quantity $E_a[T(U)^p]$ (known as the {\it $p$-th moment}) therefore provides a measurement of the size of $U$ and the proximity of $a$ to $U^c$. On the other hand, it is a less obvious fact that the moments also tell us something about the size and shape of a domain at infinity (i.e. far from the origin). The first major work in this direction seems to have been by Burkholder in \cite{burk}, where it was proved among other things that finiteness of the $p$-th Hardy norm of $\Omega$ is equivalent to finiteness of the $\frac{p}{2}$-th moment of $T(\Omega)$. Subsequent works on this topic include \cite{mecomb,coffee,mejmaa}.

\vski

The case $p=1$ is commonly referred to as the "torsion problem" due to its connection with mechanics, and is naturally the most tractable. The function $h(a) = \mathbb{E}_a[\tau_U]$ satisfies $\Delta h = -2$, and therefore p.d.e. techniques can be employed to great effect. \cite[Ch. 6]{sperb} contains a good account of this problem and methods of attacking it in special cases, such as when the domain in question is convex. Further results along the same lines, focusing in particular on convex domains, can be found in \cite{banuelos2002torsional,carroll2001old,keady,makar,philip}.

\vski

More generally, the distribution of $T(U)$ encodes a great deal of information about $U$. This connection has been explored in a number of places, including \cite{dimitri,davis,mac,wilson2020reconceptualising}.

\vski

Other interesting related problems have been tackled by p.d.e. methods. For example, in the famous paper \cite{banuelosdrum} (see also the related works \cite{banuelos1986brownian,banuelos1995improvement, banuelosmax}) eigenvalue techniques are used to demonstrate relationships between $\mathbb{E}_a[\tau_U]$ and geometric qualities of the domain, such as the size of the hyperbolic density and the inradius (the radius of the largest disk contained in the domain). The methods developed there have been extended by other authors in a number of different directions. For example, in \cite{mendez} a number of related stochastic domination results were proved concerning convex domains in $\RR^n$ and various types of symmetrizations. These results allow conclusions to be reached concerning the comparison of $p$-th moments of the exit times from these domains. One striking consequence of the eigenvalue methods is the fact that over all domains with a given area the disk maximizes the $p$-th moment of the exit time of Brownian motion for all $p$. The recent work \cite{kimmy} contains a discussion and refinements of this result, and the preprint \cite{banuelos2020bounds} (see also the related work \cite{usexitmax}) also deals with related problems and connections to eigenvalues. Other works tackling similar issues involving eigenvalues include \cite{banuelos1989heat}

\vski

The exit distribution of Brownian motion from a domain is also known as {\it harmonic measure}, and is well known to complex (and other) analysts. The highly readable text \cite{garnett2005harmonic} provides an excellent introduction to the subject, and includes material on Brownian motion. Other sources which make explicit this connection would include \cite{betsakos1998harmonic, betsakos1998polarization, betsakos2002hitting, betsakos2003two, betsakos2008equality, betsakos2008some, betsakos2008symmetrization, doob2012classical, tao2020sendov}.

\vski

A type of domain that has yielded a particularly rich set of results in connection to the exit time and distribution of Brownian motion is the conic sections. This is doubtless due to their definition in terms of quadratic polynomials, and the importance of quadratic polynomials in the study of Brownian motion. A set of papers containing a variety of results on Brownian motion and conics is \cite{banuelos2005sharp,banuelos2001first,carroll2006brownian,mesf}.

\vski

If we start a planar Brownian motion $B_t$ at the point $1$ and let $\theta_t = arg(B_t)$, then, properly scaled, $\theta_t$ is approximately Cauchy distributed. This is known as {\it Spitzer's Theorem}, and has generated considerable interest. Many proofs and variants rely on complex analytic techniques, including \cite{baudoin2020quaternionic,bentkus2003optimal,brassesco2014density,durrett1982new}.

\vski

Interesting connections have been found between planar Brownian motion and the hyperbolic metric on simply connected domains in \cite{banuelos1993conditioned,banuelos2000extremal,banuelos1998inradius}. 

\vski
The literature on SLE is too vast for me to survey, but an interested reader can find good starting points in \cite{rohde2011oded,schramm2011conformally,schramm2005basic}.

\vski

Though containing no reference to Brownian motion, I cannot resist the chance to mention my favorite complex analysis book, \cite{schaum}.

\backmatter
Glossary
\bibliographystyle{acm} %The style you want to use for references.
\bibliography{CABMbib} %The files containing all the articles and books you ever referenced.

\begin{thebibliography}{10}

\bibitem{banuelos1986brownian}
{\sc Ba{\~n}uelos, R.}
\newblock Brownian motion and area functions.
\newblock {\em Indiana University Mathematics Journal 35}, 3 (1986), 643--668.

\bibitem{banuelos1993conditioned}
{\sc Ba{\~n}uelos, R., and Carroll, T.}
\newblock Conditioned {B}rownian motion and hyperbolic geodesics in simply
  connected domains.
\newblock {\em The Michigan Mathematical Journal 40}, 2 (1993), 321--332.

\bibitem{banuelosdrum}
{\sc Ba{\~n}uelos, R., and Carroll, T.}
\newblock Brownian motion and the fundamental frequency of a drum.
\newblock {\em Duke Mathematical Journal 75}, 3 (1994), 575--602.

\bibitem{banuelos1995improvement}
{\sc Ba{\~n}uelos, R., and Carroll, T.}
\newblock An improvement of the {O}sserman constant.
\newblock {\em Stochastic Analysis 57\/} (1995), 3.

\bibitem{banuelos2000extremal}
{\sc Ba{\~n}uelos, R., and Carroll, T.}
\newblock Extremal problems for conditioned {B}rownian motion and the
  hyperbolic metric.
\newblock In {\em Annales de l'institut Fourier\/} (2000), vol.~50,
  pp.~1507--1532.

\bibitem{banuelos2005sharp}
{\sc Ba{\~n}uelos, R., and Carroll, T.}
\newblock Sharp integrability for {B}rownian motion in parabola-shaped regions.
\newblock {\em Journal of Functional Analysis 218}, 1 (2005), 219--253.

\bibitem{banuelosmax}
{\sc Ba{\~n}uelos, R., and Carroll, T.}
\newblock The maximal expected lifetime of {B}rownian motion.
\newblock {\em Mathematical Proceedings of the Royal Irish Academy 111\/}
  (2011).

\bibitem{banuelos1998inradius}
{\sc Ba{\~n}uelos, R., Carroll, T., and Housworth, E.}
\newblock Inradius and integral means for {G}reen’s functions and conformal
  mappings.
\newblock {\em Proceedings of the American Mathematical Society 126}, 2 (1998),
  577--585.

\bibitem{banuelos1989heat}
{\sc Ba{\~n}uelos, R., and Davis, B.}
\newblock Heat kernel, eigenfunctions, and conditioned {B}rownian motion in
  planar domains.
\newblock {\em Journal of functional analysis 84}, 1 (1989), 188--200.

\bibitem{banuelos2001first}
{\sc Ba{\~n}uelos, R., DeBlassie, R.~D., and Smits, R.}
\newblock The first exit time of planar {B}rownian motion from the interior of
  a parabola.
\newblock {\em The Annals of Probability 29}, 2 (2001), 882--901.

\bibitem{banuelos2020bounds}
{\sc Ba{\~n}uelos, R., Mariano, P., and Wang, J.}
\newblock Bounds for exit times of {B}rownian motion and the first {D}irichlet
  eigenvalue for the {L}aplacian.
\newblock {\em arXiv:2003.06867\/} (2020).

\bibitem{banuelos2002torsional}
{\sc Ba{\~n}uelos, R., Van~den Berg, M., and Carroll, T.}
\newblock Torsional rigidity and expected lifetime of {B}rownian motion.
\newblock {\em Journal of the London Mathematical Society 66}, 2 (2002),
  499--512.

\bibitem{bass1994probabilistic}
{\sc Bass, R.}
\newblock {\em Probabilistic techniques in analysis}.
\newblock Springer Science \& Business Media, 1994.

\bibitem{baudoin2020quaternionic}
{\sc Baudoin, F., Demni, N., and Wang, J.}
\newblock Quaternionic {B}rownian windings.
\newblock {\em Journal of Theoretical Probability\/} (2020), 1--18.

\bibitem{bentkus2003optimal}
{\sc Bentkus, V., Pap, G., and Yor, M.}
\newblock Optimal bounds for cauchy approximations for the winding distribution
  of planar {B}rownian motion.
\newblock {\em Journal of Theoretical Probability 16}, 2 (2003), 345--361.

\bibitem{betsakos1998harmonic}
{\sc Betsakos, D.}
\newblock Harmonic measure on simply connected domains of fixed inradius.
\newblock {\em Arkiv f{\"o}r Matematik 36}, 2 (1998), 275--306.

\bibitem{betsakos1998polarization}
{\sc Betsakos, D.}
\newblock Polarization, conformal invariants, and {B}rownian motion.
\newblock {\em Annales Academiae Scientiarum Fennicae Mathematica 23\/} (1998),
  59--82.

\bibitem{betsakos2002hitting}
{\sc Betsakos, D.}
\newblock Hitting probabilities of conditional {B}rownian motion and
  polarisation.
\newblock {\em Bulletin of the Australian Mathematical Society 66}, 2 (2002),
  233--244.

\bibitem{betsakos2003two}
{\sc Betsakos, D.}
\newblock Two-point projection estimates for harmonic measure.
\newblock {\em Bulletin of the London Mathematical Society 35}, 4 (2003),
  473--478.

\bibitem{betsakos2008equality}
{\sc Betsakos, D.}
\newblock Equality cases in the symmetrization inequalities for {B}rownian
  transition functions and {D}irichlet heat kernels.
\newblock {\em Annales Academiae Scientiarum Fennicae Mathematica 33}, 2
  (2008), 413--427.

\bibitem{betsakos2008some}
{\sc Betsakos, D.}
\newblock Some properties of $\alpha$-harmonic measure.
\newblock {\em Colloquium Math 111}, 2 (2008), 297--314.

\bibitem{betsakos2008symmetrization}
{\sc Betsakos, D.}
\newblock Symmetrization and harmonic measure.
\newblock {\em Illinois Journal of Mathematics 52}, 3 (2008), 919--949.

\bibitem{dimitri}
{\sc Betsakos, D., Boudabra, M., and Markowsky, G.}
\newblock On the probability of fast exits and long stays of planar {B}rownian
  motion in simply connected domains.
\newblock {\em Journal of Mathematical Analysis and Applications\/} (2020).

\bibitem{bishop2017fractals}
{\sc Bishop, C., and Peres, Y.}
\newblock {\em Fractals in probability and analysis}, vol.~162.
\newblock Cambridge University Press, 2017.

\bibitem{usexitmax}
{\sc Boudabra, M., and Markowsky, G.}
\newblock Maximizing the $p$-th moment of exit time of planar {B}rownian motion
  from a given domain.
\newblock {\em arXiv:1907.06335\/} (2019).

\bibitem{mecomb}
{\sc Boudabra, M., and Markowsky, G.}
\newblock The $p$-th moment of exit time of planar {B}rownian motion on comb
  domains.
\newblock {\em arXiv:1907.06335\/} (2019).

\bibitem{boudabra2020new}
{\sc Boudabra, M., and Markowsky, G.}
\newblock A new solution to the conformal {S}korokhod embedding problem and
  applications to the {D}irichlet eigenvalue problem.
\newblock {\em Journal of Mathematical Analysis and Applications\/} (2020),
  124351.

\bibitem{boudabra2020remarks}
{\sc Boudabra, M., and Markowsky, G.}
\newblock Remarks on {G}ross’ technique for obtaining a conformal {S}korohod
  embedding of planar {B}rownian motion.
\newblock {\em Electronic Communications in Probability 25\/} (2020).

\bibitem{bourg}
{\sc Bourgade, P., Fujita, T., and Yor, M.}
\newblock Euler's formulae for {$\zeta(2n)$} and products of {C}auchy
  variables.
\newblock {\em Electronic Communications in Probability 12\/} (2007), 73--80.

\bibitem{brassesco2014density}
{\sc Brassesco, S., and Pire, S.}
\newblock On the density of the winding number of planar {B}rownian motion.
\newblock {\em Journal of Theoretical Probability 27}, 3 (2014), 899--914.

\bibitem{burk}
{\sc Burkholder, D.}
\newblock Exit times of {B}rownian motion, harmonic majorization, and {H}ardy
  spaces.
\newblock {\em Advances in Mathematics 26}, 2 (1977), 182--205.

\bibitem{carroll2001old}
{\sc Carroll, T.}
\newblock Old and new on the bass note, the {T}orsion function and the
  hyperbolic metric.
\newblock {\em Irish Math. Soc. Bull 47\/} (2001), 41--65.

\bibitem{carroll2006brownian}
{\sc Carroll, T.}
\newblock Brownian motion and harmonic measure in conic sections.
\newblock In {\em Potential Theory in Matsue\/} (2006), Mathematical Society of
  Japan, pp.~25--41.

\bibitem{chin2019some}
{\sc Chin, W., Jung, P., and Markowsky, G.}
\newblock Some remarks on invariant maps of the {C}auchy distribution.
\newblock {\em Statistics \& Probability Letters 158\/} (2020), 108668.

\bibitem{chung2002green}
{\sc Chung, K.~L.}
\newblock {\em Green, {B}rown, and {P}robability and {B}rownian {M}otion on the
  {L}ine}.
\newblock World Scientific Publishing Company, 2002.

\bibitem{coffee}
{\sc Coffey, M.}
\newblock Expected exit times of {B}rownian motion from planar domains:
  Complements to a paper of {M}arkowsky.
\newblock {\em {\it arXiv:1203.5142} Preprint\/}.

\bibitem{davis1975picarda}
{\sc Davis, B.}
\newblock Picard's theorem and {B}rownian motion.
\newblock {\em Transactions of the American Mathematical Society 213\/} (1975),
  353--362.

\bibitem{davis}
{\sc Davis, B.}
\newblock Brownian motion and analytic functions.
\newblock {\em Annals of Probability 7}, 6 (1979), 913--932.

\bibitem{doob2012classical}
{\sc Doob, J.}
\newblock {\em Classical potential theory and its probabilistic counterpart},
  vol.~262.
\newblock Springer Science \& Business Media, 2012.

\bibitem{durrett1982new}
{\sc Durrett, R.}
\newblock A new proof of {S}pitzer's result on the winding of two dimensional
  brownian motion.
\newblock {\em The Annals of Probability\/} (1982), 244--246.

\bibitem{durBM}
{\sc Durrett, R.}
\newblock {\em Brownian motion and martingales in analysis}.
\newblock Wadsworth Advanced Books \& Software, 1984.

\bibitem{feller2008introduction}
{\sc Feller, W.}
\newblock {\em An introduction to probability theory and its applications},
  vol.~2.
\newblock John Wiley \& Sons, 2008.

\bibitem{garnett2005harmonic}
{\sc Garnett, J., and Marshall, D.}
\newblock {\em Harmonic measure}, vol.~2.
\newblock Cambridge University Press, 2005.

\bibitem{gross2019conformal}
{\sc Gross, R.}
\newblock A conformal {S}korokhod embedding.
\newblock {\em Electronic Communications in Probability 24\/} (2019).

\bibitem{keady}
{\sc Keady, G., and McNabb, A.}
\newblock The elastic torsion problem: solutions in convex domains.
\newblock {\em New Zealand Journal of Mathematics 22\/} (1993), 43--64.

\bibitem{kimmy}
{\sc Kim, D.}
\newblock Quantitative inequalities for the expected lifetime of {B}rownian
  motion.
\newblock {\em arXiv: 1904.09565\/} (2019).

\bibitem{letac1977functions}
{\sc Letac, G.}
\newblock Which functions preserve {C}auchy laws?
\newblock {\em Proceedings of the American Mathematical Society 67}, 2 (1977),
  277--286.

\bibitem{makar}
{\sc Makar-Limanov, L.}
\newblock Solution of {D}irichlet's problem for the equation $\bigtriangleup
  u=-1$ in a convex region.
\newblock {\em Mathematical Notes of the Academy of Sciences of the USSR 9}, 1
  (1971), 52--53.

\bibitem{mariano2020conformal}
{\sc Mariano, P., and Panzo, H.}
\newblock Conformal {S}korokhod embeddings and related extremal problems.
\newblock {\em Electronic Communications in Probability 25\/} (2020).

\bibitem{mesf}
{\sc Markowsky, G.}
\newblock A method for deriving hypergeometric and related identities from the
  ${H}^2$ {H}ardy norm of conformal maps.
\newblock {\em Integral Transforms and Special Functions 24}, 4 (2013),
  302--313.

\bibitem{markowsky2014probabilistic}
{\sc Markowsky, G.}
\newblock A probabilistic proof of the open mapping theorem for analytic
  functions.
\newblock {\em Bulletin of the Australian Mathematical Society 90}, 1 (2014),
  74--76.

\bibitem{mejmaa}
{\sc Markowsky, G.}
\newblock The exit time of planar {B}rownian motion and the
  {P}hragm{\'e}n--{L}indel{\"o}f principle.
\newblock {\em Journal of Mathematical Analysis and Applications 422}, 1
  (2015), 638--645.

\bibitem{megreen}
{\sc Markowsky, G.}
\newblock On the planar {B}rownian {G}reen's function for stopping times.
\newblock {\em Journal of Mathematical Analysis and Applications 455\/} (2017),
  1221--1233.

\bibitem{markowskyfenn}
{\sc Markowsky, G.}
\newblock On the distribution of planar {B}rownian motion at stopping times.
\newblock {\em Annales Academiae Scientiarum Fennicae Mathematica 43\/} (2018),
  597--616.

\bibitem{markowsky2018remark}
{\sc Markowsky, G.}
\newblock A remark on the probabilistic solution of the {D}irichlet problem for
  simply connected domains in the plane.
\newblock {\em Journal of Mathematical Analysis and Applications 464}, 2
  (2018), 1143--1146.

\bibitem{mac}
{\sc McConnell, T.}
\newblock The size of an analytic function as measured by {L}\'evy's time
  change.
\newblock {\em The Annals of Probability 13}, 3 (1985), 1003--1005.

\bibitem{MR1394988}
{\sc McCullagh, P.}
\newblock M\"{o}bius transformation and {C}auchy parameter estimation.
\newblock {\em Ann. Statist. 24}, 2 (1996), 787--808.

\bibitem{mckean1969stochastic}
{\sc McKean, H.}
\newblock {\em Stochastic integrals}, vol.~353.
\newblock American Mathematical Soc., 1969.

\bibitem{melgreen}
{\sc Melnikov, Y.}
\newblock A new approach to the representation of trigonometric and hyperbolic
  functions by infinite products.
\newblock {\em Journal of Mathematical Analysis and Applications 344\/} (2008),
  521--534.

\bibitem{melgreen2}
{\sc Melnikov, Y.}
\newblock {\em Green's functions and infinite products: bridging the divide}.
\newblock Springer, Science \& Business Media, 2011.

\bibitem{mendez}
{\sc M{\'e}ndez-Hern{\'a}ndez, P.}
\newblock {B}rascamp-{L}ieb-{L}uttinger inequalities for convex domains of
  finite inradius.
\newblock {\em Duke Mathematical Journal 113}, 1 (2002), 93--131.

\bibitem{morters2010brownian}
{\sc M{\"o}rters, P., and Peres, Y.}
\newblock {\em Brownian motion}, vol.~30.
\newblock Cambridge University Press, 2010.

\bibitem{obloj2004skorokhod}
{\sc Ob{\l}{\'o}j, J.}
\newblock The {S}korokhod embedding problem and its offspring.
\newblock {\em Probability Surveys 1\/} (2004), 321--392.

\bibitem{okamura2020characterizations}
{\sc Okamura, K.}
\newblock Characterizations of the {C}auchy distribution associated with
  integral transforms.
\newblock {\em arXiv:2007.09707\/} (2020).

\bibitem{pascu2005probabilistic}
{\sc Pascu, M.}
\newblock A probabilistic proof of the fundamental theorem of algebra.
\newblock {\em Proceedings of the American Mathematical Society 133}, 6 (2005),
  1707--1711.

\bibitem{petersen1977brownian}
{\sc Petersen, K.}
\newblock {\em Brownian motion, {H}ardy spaces and bounded mean oscillation},
  vol.~28.
\newblock Cambridge University Press, 1977.

\bibitem{philip}
{\sc Philippin, G., and Porru, G.}
\newblock Isoperimetric inequalities and overdetermined problems for the
  {S}aint-{V}enant equation.
\newblock {\em New Zealand Journal of Mathematics 25\/} (1996), 217--227.

\bibitem{pitman1967cauchy}
{\sc Pitman, E., and Williams, E.}
\newblock Cauchy-distributed functions of {C}auchy variates.
\newblock {\em The Annals of Mathematical Statistics\/} (1967), 916--918.

\bibitem{pitman2018guide}
{\sc Pitman, J., and Yor, M.}
\newblock A guide to {B}rownian motion and related stochastic processes.
\newblock {\em arXiv preprint arXiv:1802.09679\/} (2018).

\bibitem{revyor}
{\sc Revuz, D., and Yor, M.}
\newblock {\em {Continuous martingales and Brownian motion}}.
\newblock Springer Verlag, 1999.

\bibitem{rogers2000diffusions}
{\sc Rogers, L., and Williams, D.}
\newblock Diffusions, {M}arkov processes and martingales, volume 1:
  {F}oundations.
\newblock {\em Cambridge Mathematical Library,\/} (2000).

\bibitem{rohde2011oded}
{\sc Rohde, S.}
\newblock Oded {S}chramm: from circle packing to {SLE}.
\newblock In {\em Selected Works of Oded Schramm}. Springer, 2011, pp.~3--45.

\bibitem{schramm2011conformally}
{\sc Schramm, O.}
\newblock Conformally invariant scaling limits: an overview and a collection of
  problems.
\newblock In {\em Selected Works of Oded Schramm}. Springer, 2011,
  pp.~1161--1191.

\bibitem{schramm2005basic}
{\sc Schramm, O., and Rohde, S.}
\newblock Basic properties of {SLE}.
\newblock {\em Annals of Mathematics 161}, 2 (2005), 883--924.

\bibitem{sperb}
{\sc Sperb, R.}
\newblock {\em Maximum principles and their applications}.
\newblock Elsevier, 1981.

\bibitem{schaum}
{\sc Spiegel, M.}
\newblock {\em Complex variables: {W}ith an introduction to conformal mapping
  and its applications}.
\newblock Schaum's Outline Series, 2009.

\bibitem{stroock2010probability}
{\sc Stroock, D.}
\newblock {\em Probability theory: an analytic view}.
\newblock Cambridge university press, 2010.

\bibitem{tao2020sendov}
{\sc Tao, T.}
\newblock Sendov's conjecture for sufficiently high degree polynomials.
\newblock {\em arXiv preprint arXiv:2012.04125\/} (2020).

\bibitem{wilson2020reconceptualising}
{\sc Wilson, D., Woodhouse, F., Simpson, M., and Baker, R.}
\newblock Reconceptualising transport in crowded and complex environments.
\newblock {\em arXiv:2006.16758\/} (2020).

\end{thebibliography}
\printindex %Make an index AUTOMATICALLY

\end{document}